\documentclass[preprint, 12pt]{elsarticle}
\pdfoutput=1
\usepackage{amsmath}
\usepackage{amssymb}
\usepackage{bbm}
\usepackage{graphicx}
\usepackage{wasysym}
\usepackage{subfig}
\usepackage{fullpage}
\usepackage{color}
\newdefinition{remark}{Remark}
\newdefinition{proposition}{Proposition}
\newdefinition{proof}{Proof}
\newdefinition{algorithm}{Algorithm}

\usepackage{lmodern}
\usepackage[T1]{fontenc}
\usepackage[latin9]{inputenc}
\usepackage{esint}

\usepackage{hyperref}

\numberwithin{equation}{section}
\numberwithin{figure}{section}

\newcommand{\du}{\, \mathrm{d}}
\newcommand{\myremarkend}{~\hfill$\spadesuit\/$}

\journal{...} 
\begin{document}

\begin{frontmatter}

\title{A Fourier penalty method for solving the time-dependent Maxwell's equations in domains with curved boundaries}

\author{R. Galagusz$^a$}
\ead{ryan.galagusz@mail.mcgill.ca}

\author{D. Shirokoff$^b$}
\ead{david.g.shirokoff@njit.edu}

\author{J.-C. Nave$^c$}
\ead{jcnave@math.mcgill.ca}

\address{$^a$Department of Electrical and Computer Engineering, McGill University,
         Montreal, QC, H3A0E9, CAN}

\address{$^b$Department of Mathematical Sciences, NJIT, Newark, NJ, 07102-1982, USA}         

\address{$^c$Department of Mathematics and Statistics, McGill University,
         Montreal, QC, H3A0B9, CAN}

\begin{abstract}
	We present a high order, Fourier penalty method for the Maxwell's equations in the vicinity of perfect electric conductor boundary conditions. The approach relies on extending the smooth non-periodic domain of the equations to a periodic domain by removing the exact boundary conditions and introducing an analytic forcing term in the extended domain. The forcing, or penalty term is chosen to systematically enforce the boundary conditions to high order in the penalty parameter, which then allows for higher order numerical methods. We present an efficient numerical method for constructing the penalty term, and discretize the resulting equations using a Fourier spectral method. We demonstrate convergence orders of up to $3.5$ for the one-dimensional Maxwell's equations, and show that the numerical method does not suffer from dispersion (or pollution) errors. We also illustrate the approach in two dimensions and demonstrate convergence orders of $2.5$ for transverse magnetic modes and $1.5$ for the transverse electric modes.  We conclude the paper with numerous test cases in dimensions two and three including waves traveling in a bent waveguide, and scattering off of a windmill-like geometry. 
\end{abstract}


\end{frontmatter}

\section{Introduction}

Pseudospectral and Fourier based methods \cite{TrefethenEmbree2005} provide a popular solution approach for problems involving periodic boundary conditions.  Unfortunately, pseudospectral methods which exploit the Fourier transform do not extend easily to domains with curved boundaries.  One approach for solving partial differential equations (PDEs) on domains with curved boundaries is to relax the boundary condition by introducing a forcing, or penalty term to approximately enforce the correct boundary values.  Such an approach has successfully been developed for a variety of problems in fluid dynamics \cite{AngotCaltagiron1990, ArquisCaltagirone1984, Angot2010, SarthouVincentCaltagironeAngot2008} as well as computations involving turbulent flows \cite{KevlahanGhidaglia2001}.  Other more recent applications include using penalty equations in ocean modeling \cite{ReckingerVasilyevKemper2012}, plasma physics \cite{AngotAuphanGues2014}, magneto-hydrodynamics \cite{MoralesLeroyBosSchneider2012}, and scalar advection with moving obstacles \cite{KadochKolomenskiyAngotSchneider2012}. One significant drawback with such volume based penalty methods is the introduction of analytic errors in the penalized PDE.  The resulting analytic error not only limits the accuracy of any numerical method, but also degrades the smoothness of the underlying solution.  As a result of the reduced regularity in the penalized solution, the Fourier spectral methods typically require additional filtering steps \cite{KolomenskiySchneider2009}.

In recent work \cite{ShirokoffNave2014}, a new modified penalty term was introduced to alleviate the analytic error due to the standard volume penalty method.  The approach was examined for the heat and Poisson equations to obtain a third order Fourier-based method. The method was then extended to the Navier-Stokes equations to obtain a second order Fourier scheme.  

The focus in the current paper is on hyperbolic wave equations with an emphasis on Maxwell's equations. Specifically, we focus on the time-dependent Maxwell's equations in free space in the presence of perfect electric conductors (PEC).  Perfect conductors are idealized materials that easily conduct electricity and are accompanied with corresponding boundary conditions. Mathematically, PEC boundary conditions are modeled by assuming the electric field is normal to the boundary of the conducting material. Such a condition may then be converted into an appropriate Dirichlet boundary condition on the underlying PDE.

In contrast to previous work \cite{ShirokoffNave2014} which focused primarily on elliptic and parabolic equations, here a modified approach must be applied for hyperbolic systems.  Specifically, the penalty term cannot be directly applied to a second order wave equation as it will introduce spurious oscillations in time, but rather must be introduced into the first order system so as to dampen solutions.  Even with the suitable introduction of a penalty term to a hyperbolic system, the presence of analytic errors can significantly limit the accuracy of a numerical method.  For instance we demonstrate that a conventional volume penalty method will converge at a rate of $0.5$, namely the error scales as $O(\Delta x^{1/2})$ where $\Delta x$ is the grid spacing of the scheme. Recent work by \cite{Auphan2014} suggests that an alternative penalization may yield first order methods, while other work  \cite{DymkoskiKasimovVasilyev2014} shows second order convergence rates for a class of hyperbolic systems with Neumann boundary conditions. Finally, a similar in spirit approach \cite{Carpenter1999, Hesthaven2007, CarpenterNordstromGottlieb2010, Nordstrom2014}, where an additional penalty term is prescribed to connect subdomains, or to enforce boundary conditions was developed to obtain provably stable numerical schemes.  Although that method is currently limited to low order for boundary conditions \cite{Nordstrom2014}, it is hopeful that future work may lead to the development of provably higher order penalty methods.

Another successful approach for solving wave problems with Fourier series is through Fourier extension methods \cite{BrunoLyon2010, LyonBruno2010}.  The methods have been very successful at obtaining highly accurate solutions for wave problems that do not have a divergence constraint. In particular, the Fourier method is combined with an iterative (alternate direction iteration) method to solve a sequence of elliptic problems as a means to evolve wave equations.  The methods we propose in this paper differ as they may be discretized with an explicit in time method and therefore do not require solving an elliptic problem at each time iteration.  

We emphasize that our approach is a single domain pseudospectral in space finite difference in time method. Previous single domain pseudospectral time-domain (PSTD) approaches \cite{Liu1997, Liu1999, Liu1999a, Liu1999b} cannot handle curved geometries with PEC boundary conditions. Our approach can be thought of as a new way to extend the single domain PSTD method to domains with curved geometries. In addition, our approach preserves the use of the fast Fourier transform (FFT) and does not suffer from dispersion errors. In subsequent developments of the PSTD method \cite{Fan2002} the FFT is no longer used, multiple domains must be introduced, and accuracy is lost due to subdomain coupling.

In methods such as the immersed boundary or standard penalty method, the extended solution is no longer smooth. The lack of smoothness then limits the convergence rate.  As part of our approach, we ensure that the forcing creates an extension that is smooth in a precise sense. We demonstrate that with an appropriate modification and introduction of an active penalty term, one may achieve systematically higher order methods.  Specifically, we show that for problems in one dimension, one may achieve convergence rates of up to 3.5 (the limitation currently due to time stepping), while in dimension two, one may obtain rates of 1.5 for transverse electric (TE) modes and 2.5 for transverse magnetic (TM) modes.

In the first half of the paper we introduce the Maxwell's equations with PEC boundary conditions, along with the formulation of the active penalty term. We also describe the analytic construction of the penalty term for TE and TM modes in dimension two. We then examine the analytic error in the penalty parameter for scattering of a TM mode off of a PEC wall. 
The second half of the paper focuses on the numerical implementation of solving the penalized Maxwell's equations using a Fourier pseudospectral approach.  Specifically, we provide details on how to numerically discretize the equations in both space and time using equispaced grids and Fourier series.  We then go on to outline details of stability studies in dimensions one and two and illustrate how the penalty term can be combined with PMLs to provide full time-dependent simulations of waves with PEC and radiating boundary conditions on periodic domains.  In addition, we validate the approach by performing several numerical studies.  Specifically, we show that in dimension one, the Fourier spectral method does not suffer from pollution (numerical dispersion) errors. We perform convergence studies in both one and two dimensions, showing global convergence rates of up to $3.5$ in dimension one, $1.5$ for TE modes in dimension two and $2.5$ for TM modes in dimension two.  Lastly, we illustrate the utility of the approach on some problems involving windmill shaped and waveguide geometries and demonstrate the natural extension to three dimensions.

\section{Basic approach}

In this paper we develop numerical Fourier methods for solving the time-dependent boundary value problem for Maxwell's equations.  Specifically, we focus on solving Maxwell's equations for isotropic space in the vicinity of PEC.  We denote the region of isotropic space by $\Omega_0 \subset \Omega = [0, D]^d$, for $d = 1, 2, 3$ where $[0, D]^d$ is the $d$-dimensional cube with periodic boundary conditions, and the boundary $\Gamma = \partial \Omega_0$. The Maxwell's equations then take the form
\begin{subequations} \label{MaxwellEq}
	\begin{alignat}{2}
		\frac{\partial \mathbf{H}}{\partial t} &= - \nabla \times \mathbf{E}, \quad && \textrm{in } \Omega_0 \times (0, T] \\
        		\frac{\partial \mathbf{E}}{\partial t} &= \phantom{-} \nabla \times \mathbf{H}, \quad && \textrm{in } \Omega_0 \times (0, T]\\
		\nabla \cdot \mathbf{E} \; &= \; 0, \quad && \textrm{in } \Omega_0 \times (0, T]\\
		\nabla \cdot \mathbf{H} \; &= \; 0, \quad && \textrm{in } \Omega_0 \times (0, T] \\ \label{MaxwellEq_BC}
		\mathbf{n} \times (\mathbf{E} - \mathbf{g}) &= \; 0, \quad && \textrm{in } \Gamma \times (0, T].
	\end{alignat} 
\end{subequations}

Here we work with rescaled variables $\mathbf{E}$ and $\mathbf{H}$ so that effectively $\epsilon = \mu = 1$ and $c = 1$. For instance, upon non-dimensionalizing Maxwell's equations by rescaling $(t, \mathbf{x})$ and $(\mathbf{E}, \mathbf{H})$, one arrives at equations (\ref{MaxwellEq}). 

Although our focus will largely be on PEC, we consider a more general set of boundary conditions where one prescribes a Dirichlet tangential boundary condition for $\mathbf{E}$ as some general function $\mathbf{g}$ of time. As written in the formulation (\ref{MaxwellEq}), $\mathbf{n}$ is the inward unit normal to $\Omega_0$, while $\mathbf{g}$ is the prescribed tangential component of $\mathbf{E}$ on $\Gamma$ (we assume without loss of generality that $\mathbf{n}\cdot\mathbf{g} = 0$). 

\begin{figure}[htb!]
	\centering
    \includegraphics[width = \textwidth]{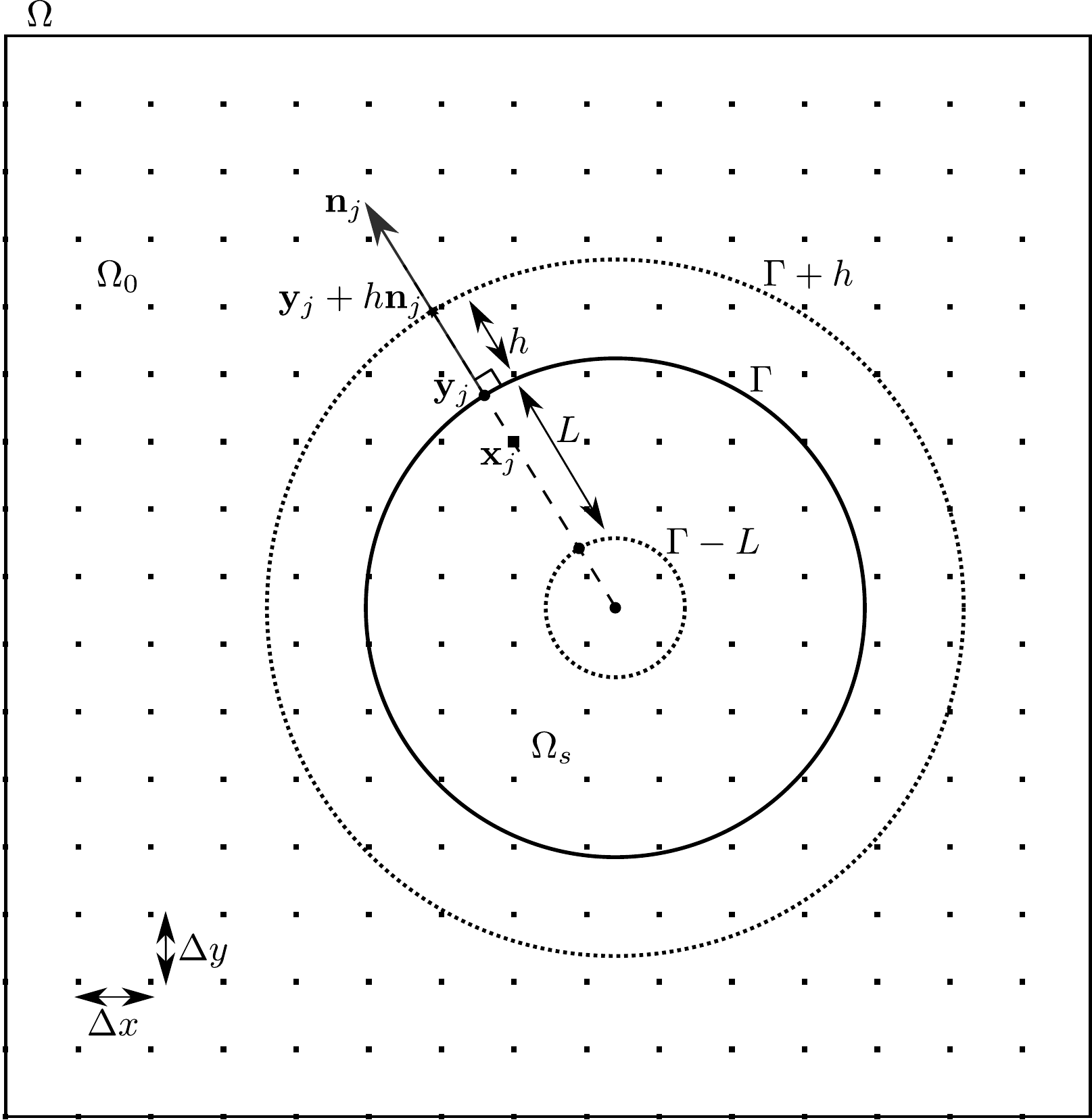} \\
    \caption{Example domain $\Omega$ with immersed PEC with boundary $\Gamma$. The plot shows a regular grid while the points $\mathbf{x}_j$, $\mathbf{y}_j$ and $\mathbf{y}_j + h\mathbf{n}_j$ are used in the construction of the extension function $\tilde{\mathbf{g}}$.} \label{Fig_Domain}
\end{figure}

A particularly practical case is that of a PEC where the complimentary domain $\Omega_s = [0, D]^d \setminus \overline{\Omega}_0$  (see Figure \ref{Fig_Domain}) is an electric conductor. In this case, one assumes $\mathbf{E} = \mathbf{H} = 0$ inside $\Omega_s$.  Due to the presence of either surface charges or currents, only the tangential component of $\mathbf{E}$ and normal component of $\mathbf{H}$ are then continuous across the interface $\Gamma$, resulting in

\begin{subequations} \label{Maxwell_BC_PEC}
	\begin{alignat}{2} 
		\mathbf{n} \times \mathbf{E} &= 0, \quad && \textrm{on } \Gamma \\
		\mathbf{n} \cdot \mathbf{H}  &= 0, \quad && \textrm{on } \Gamma.
	\end{alignat}
\end{subequations}
	
Note that the two boundary conditions (\ref{Maxwell_BC_PEC})  are equivalent. 
Given initial data $\mathbf{E}(\mathbf{x}, 0) = \mathbf{E}_0(\mathbf{x})$, $\mathbf{H}(\mathbf{x}, 0) = \mathbf{H}_0(\mathbf{x})$ satisfying the compatibility conditions $\nabla \cdot \mathbf{E}_0 = \nabla \cdot \mathbf{H}_0 = 0$ and $\mathbf{n}\times \mathbf{E}_0 = \mathbf{n} \times \mathbf{g}$, we seek a solution for (\ref{MaxwellEq}).

\subsection{Penalized Equations}

We now outline how to analytically modify the equations (\ref{MaxwellEq}) in the presence of PEC so that one may numerically solve them using time-dependent Fourier methods.  The approach relies on extending the domain $\Omega_0$ to $\Omega = \Omega_0 \bigcup \Omega_s$ and suitably modifying the equations (\ref{MaxwellEq}) inside $\Omega_s$ by the introduction of a penalty term.  For the practical implementation using Fourier methods, we take $\Omega = [0, D]^d$ to be a rectangle with periodic boundary conditions.  We then solve the full penalized equations on $\Omega$ with the understanding that the restriction of the solution to $\Omega_0$ represents the physical solution.  Meanwhile the solution on $\Omega_s$ is fictitious and used only to aid in the numerical computation. 

The modified penalty equations take the form
\begin{subequations} \label{Penalized_Maxwell}
   \begin{alignat}{2}
     \frac{\partial \mathbf{H}_{\eta}}{\partial t} & =  - \nabla \times \mathbf{E}_{\eta}, \quad & \textrm{on } \Omega
        \\ 
     \frac{\partial \mathbf{E}_{\eta}}{\partial t}  & = \phantom{-} \nabla \times \mathbf{H}_{\eta} - \eta^{-1} \chi_h(\mathbf{x})(\mathbf{E}_{\eta} - \tilde{\mathbf{g}} ), \quad & \textrm{on } \Omega
   \end{alignat}
 \end{subequations}

\begin{subequations}
	\begin{alignat}{2}  \label{PenalizedDivergenceConstraint}
		\nabla \cdot \mathbf{E}_{\eta} &= 0, \quad & \textrm{on } \Omega_0 \phantom{.} \\
		\nabla \cdot \mathbf{H}_{\eta} &= 0, \quad & \textrm{on } \Omega_0.
	\end{alignat}
\end{subequations}
Here $\tilde{\mathbf{g}}$ is an active penalty function, and $\chi_h(\mathbf{x})$ is a characteristic function defined by
\begin{eqnarray}
	\chi_h(\mathbf{x}) = \begin{cases}
		0, & \text{for } \text{dist}(\mathbf{x}, \Omega_s) > h \\
		1, & \text{for } \text{dist}(\mathbf{x}, \Omega_s) \leq h
	\end{cases}.
\end{eqnarray}
In other words, $\chi_h(\mathbf{x}) = 1$ if $\mathbf{x} \in \Omega_s$ or within a distance\footnote{Here $\text{dist}(\mathbf{x}, \Omega_s) = \inf_{\mathbf{y} \in \Omega_s} |\mathbf{x} - \mathbf{y}|$ is the distance of the point $\mathbf{x}$ to the set $\Omega_s$.} $h$ to the set $\Omega_s$.  

The goal is to choose $\tilde{\mathbf{g}}$ so that the penalized solution $(\mathbf{E}_{\eta}, \mathbf{H}_{\eta})$ with the same initial data converges rapidly to the exact solution $(\mathbf{E}, \mathbf{H})$
\begin{eqnarray} \label{ConvergenceLimit_E}
	\lim_{h, \eta \rightarrow 0} \mathbf{E}_{\eta} \rightarrow \mathbf{E}, \quad \quad \lim_{h, \eta \rightarrow 0} \mathbf{H}_{\eta} \rightarrow \mathbf{H}.
\end{eqnarray}
In such a case, solving the penalized equations (\ref{Penalized_Maxwell}) with small $h$ and $\eta$ provide accurate approximations to the true fields.  

In the following subsections, we outline how to construct $\tilde{\mathbf{g}}$ to satisfy (\ref{ConvergenceLimit_E}).   Although there are some similarities with the procedure outlined in \cite{ShirokoffNave2014}, the new method described here differs in the sense that (at the level of a continuum PDE) the penalty term $\chi_h(\mathbf{x})(\mathbf{E} - \tilde{\mathbf{g}})$ is continuous. In other words, we choose $\tilde{\mathbf{g}}$ to continuously match $\mathbf{E}_{\eta}$ at the jump discontinuity in $\chi_h$.  In fact, when $\tilde{\mathbf{g}}$ matches $m$ derivatives at the jump, then for a fixed $\eta$, $h$, and smooth enough boundary $\Gamma$, the forcing term $\eta^{-1} \chi_h(\mathbf{x})(\mathbf{E}_{\eta} - \tilde{\mathbf{g}} )$ is $C^m$.  The construction then turns out to be simpler and more accurate than in \cite{ShirokoffNave2014}.

\subsection{Penalty function $\tilde{\mathbf{g}}$ in one dimension}\label{Sec_PenaltyFunction_dim1}

\begin{figure}[htb!]
	\centering
    \includegraphics[width = \textwidth]{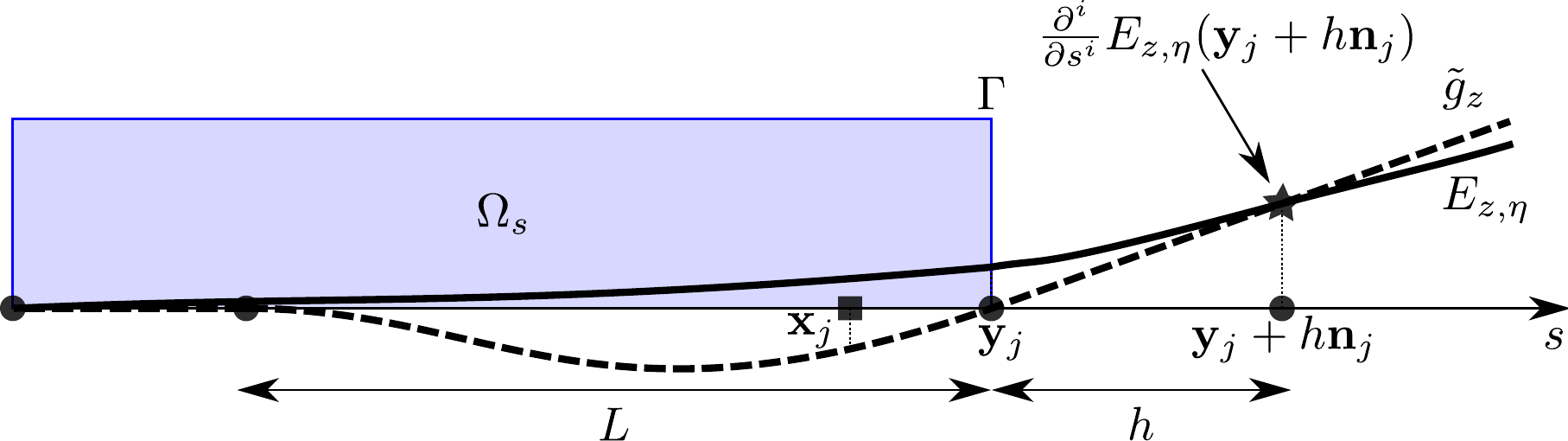} \\
    \caption{Construction of $\mathbf{\tilde{g}}$ along the direction normal to the interface $\Gamma$.} \label{Fig_OneD_Construction}
\end{figure}

We start by explicitly presenting the construction of $\tilde{\mathbf{g}}$ in one dimension.  For this construction, we assume that the boundary of the domain is located at $x = 0$, so that $\Omega_s = \{ x < 0\}$ and $\Omega_0 = \{ x > 0\}$. In addition we assume that the fields take the form $\mathbf{E}_{\eta} = (0, 0, E_{z, \eta}(x, t))$, $\mathbf{H}_{\eta} = (0, H_{y, \eta}(x, t), 0)$, and that the exact solution $E_z(0, t) = g$ satisfies a Dirichlet boundary condition at $x = 0$. In this case we take $\tilde{\mathbf{g}} = (0, 0, \tilde{g}_z(x))$ to have only one component. The Maxwell's equations then become
\begin{subequations} \label{Maxwell_Penalty_dim1}
	\begin{align}
	\frac{\partial H_{y, \eta}}{\partial t} &= \frac{\partial E_{z, \eta}}{\partial x} \\
	\frac{\partial E_{z, \eta}}{\partial t} &= \frac{\partial H_{y, \eta}}{\partial x}  - \eta^{-1} \chi_h(x) (E_{z, \eta} - \tilde{g}_z), 
	\end{align}
\end{subequations}
where 
\begin{eqnarray*}
	\chi_h(x) = \begin{cases}
		0, & x > h \\
		1, & x \leq h
	\end{cases}.
\end{eqnarray*}
The prescription is now to build $\tilde{g}_z(x)$ as a smooth extension of $E_{z, \eta}$ that also goes through the exact boundary condition $\tilde{g}_z(0) = g$.  To obtain an extension, we match $m$ derivatives of $E_{z, \eta}$ at $x = h$, for instance see Figure \ref{Fig_OneD_Construction}.  We choose for $\tilde{g}_z(x)$ to be supported on the interval $[-L, h]$, where $L > 0$ is an order 1 parameter and $h > 0$ is a small parameter. Eventually $L$ and $h$ will be fixed by numerical considerations. Explicitly, we take $\tilde{g}_z(x)$ as a polynomial extension of degree $(2m+3)$ to smoothly extend $E_{z, \eta}$ at $x = h$ and decay to $0$ at $x = -L$:
\begin{enumerate}
	\item[A.] Matching $m = 0$ derivatives at $x = h$, which will yield a $(\Delta x)^{1.5}$ order scheme
		\begin{align} \label{Extension_dim1_m0}
			\tilde{g}_z(x) = g \; P_{0, 0}(x) + E_{z, \eta}(h)\; P_{0, 1}(x) 
		\end{align}
		where
		\begin{align*}
			P_{0,0}(x) &= -\frac{(x-h)(x+L)^2}{h L^2}, &
			P_{0,1}(x) = \frac{x(x+L)^2}{h(h+L)^2}.
		\end{align*}	
		
	\item[B.] Matching $m = 1$ derivatives at $x = h$, which will yield a $(\Delta x)^{2.5}$ order scheme
		\begin{align} \label{Extension_dim1_m1}
			\tilde{g}_z(x) = g \; P_{1, 0}(x) + E_{z, \eta}(h)\; P_{1, 1}(x) + E_{z, \eta}'(h) \; P_{1, 2}(x)
		\end{align}
		where
		\begin{align*}
			P_{1,0}(x) &= \frac{(x+L)^3}{h^2 L^3}(x-h)^2, 	& P_{1,1}(x) = \frac{x(x+L)^3}{h(h+L)^3}\Big[1 - \frac{(4h+L)}{h(h+L)}(x-h) \Big], \\
			P_{1,2}(x) &= \frac{x (x-h) (x + L)^3}{h(h+L)^3}.
		\end{align*}
	\item[C.] Matching $m = 2$ derivatives at $x = h$, which will yield a $(\Delta x)^{3.5}$ order scheme
		\begin{align} \label{Extension_dim1_m2}
			\tilde{g}_z(x) = g \; P_{2, 0}(x) + E_{z, \eta}(h)\; P_{2, 1}(x) + E_{z, \eta}'(h) \; P_{2, 2}(x) + E_{z, \eta}''(h) \; P_{2, 3}(x)
		\end{align}
		where
		\begin{align*}
			P_{2,0}(x) &= -\frac{1}{h^3 L^4}(x + L)^4 (x-h)^3, \\
			P_{2,1}(x) &= \frac{x (x+L)^4}{h (h + L)^4}\Big[1 - \frac{(5h+L)}{h(h+L)}(x - h) + \frac{(15h^2 + 6hL + L^2)}{h^2(h+L)^2}(x-h)^2 \Big], \\
			P_{2,2}(x) &= \frac{x(x-h)(x+L)^4}{h(h+L)^4} \Big[ 1 - \frac{(5h+L)}{h(h+L)}(x-h)\Big], \\
			P_{2,3}(x) &= \frac{x(x-h)^2(x+L)^4}{2h(h+L)^4}.
		\end{align*}
\end{enumerate}

\begin{remark}
	The important ingredient in constructing $\tilde{g}_z$ is to build a smooth extension of $E_{z, \eta}$ that also satisfies the exact boundary condition $g$.  As a result, the polynomial prescription described here is not unique. In fact, other constructions -- such as using an exponentially decaying basis  \cite{ShirokoffNave2014}, or solving a minimization problem -- are also feasible.  Future research involves understanding the stability properties for different extension constructions. \myremarkend 
\end{remark}

\begin{remark}
	In practice, when using the high order extension [C] with a Fourier method, one only approximately computes the derivatives $E_{z, \eta}'$, $E_{z, \eta}''$. Refer to details in the numerical implementation regarding the Fourier method. \myremarkend
\end{remark}

\begin{remark}
	The analytic convergence of the penalized solution to the underlying solution does not depend on the exact details of $\tilde{g}_z(x)$ away from the interface $\Gamma$. However, we explicitly choose $\tilde{g}_z(x)$ to decay to $0$ at $x = -L$ with a polynomial degree $m + 2$.  Such a rate ensures that the solution $E_{z, \eta}$ is smoother at $x = -L$ than at the point $x = h$. \myremarkend
\end{remark}

\subsection{Penalty function $\tilde{\mathbf{g}}$ for a TM mode}\label{Sec_PenaltyFunction_TMmode}

In the case when the initial data $\mathbf{E}_0$ and $\mathbf{H}_0$, and the subsequent solutions do not depend on the $z$ coordinate, the components of the magnetic field decouple into a transverse magnetic $(\mathrm{TM}_z)$ mode consisting of $(H_x, H_y, E_z)$ and a transverse electric $(\mathrm{TE}_z)$ mode consisting of $(E_x, E_y, H_z)$.

In such a case, we prescribe the penalized $\mathrm{TM}_z$ mode to be
\begin{subequations} \label{TM_z_Mode}
	\begin{align}
	\frac{\partial H_{x, \eta}}{\partial t}  &=	-\frac{\partial E_{z,\eta}}{\partial y}, \\
	\frac{\partial H_{y, \eta}}{\partial t} &= \phantom{-}\frac{\partial E_{z, \eta}}{\partial x}, \\
	\frac{\partial E_{z, \eta}}{\partial t} &= \phantom{-}\frac{\partial H_{y, \eta}}{\partial x} - \frac{\partial H_{x,\eta}}{\partial y} - \eta^{-1} \chi_h(\mathbf{x})(E_{z, \eta} - \tilde{g}_z).
	\end{align}
\end{subequations}
Here the penalty term can be taken to be $\tilde{\mathbf{g}} = (0, 0, \tilde{g}_z(\mathbf{x}))$, where $\tilde{g}_z(\mathbf{x})$ depends only on $E_{z, \eta}$.  
We also note that since the $\mathrm{TM}_z$ mode only contains an $E_z$ component, equations $(\ref{TM_z_Mode})$ imply that $\nabla \cdot \mathbf{H}_{\eta} = 0$ and $\nabla \cdot \mathbf{E}_{\eta} = 0$ for all time.

The primary difference between the two-dimensional $\mathrm{TM}_z$ mode (\ref{TM_z_Mode}), and the one-dimensional equations (\ref{Maxwell_Penalty_dim1}) is that $\tilde{g}_z$ is now an extension of a two-dimensional function.  To efficiently construct $\tilde{g}_z$, we follow a similar approach to \cite{ShirokoffNave2014} where we build $\tilde{g}_z$ along rays from the boundary $\Gamma$. We note that in the current formulation for a fast construction of $\tilde{g}_z$, we require that $\Gamma \in C^2$.

Again, we choose $\tilde{g}_z$ to be a continuous extension of $E_{z, \eta}$ satisfying the exact boundary conditions $g_z$ on $\mathbf{x} \in \Gamma$.  To describe the construction, we make use of the following sets of points which are located a distance $h \ll 1$ from $\Gamma$ inside $\Omega_0$, and $L$ away from $\Gamma$ inside $\Omega_s$:
\begin{align}
	\Gamma + h &= \{\mathbf{x} \in \mathbbm{R}^d: \mathbf{x} \in \Omega_0, \textrm{dist}(\mathbf{x}, \Gamma) = h\} \\
	\Gamma - L &= \{\mathbf{x} \in \mathbbm{R}^d: \mathbf{x} \in \Omega_s, \textrm{dist}(\mathbf{x}, \Gamma) = L\} .
\end{align}
We then choose $\tilde{g}_z(\mathbf{x})$ to
	\begin{enumerate}
		\item[(a)] Match $m = 0$ or $m = 1$ normal derivatives of $E_{z, \eta}$ at $\Gamma + h$, 
		\item[(b)] Go through the exact boundary condition $g(\mathbf{y})$ for any $\mathbf{y} \in \Gamma$,
		\item[(c)] Decay smoothly to 0 at $\Gamma - L$.
	\end{enumerate}	

\begin{figure}[htb!]
	\centering
    \includegraphics[width = \textwidth]{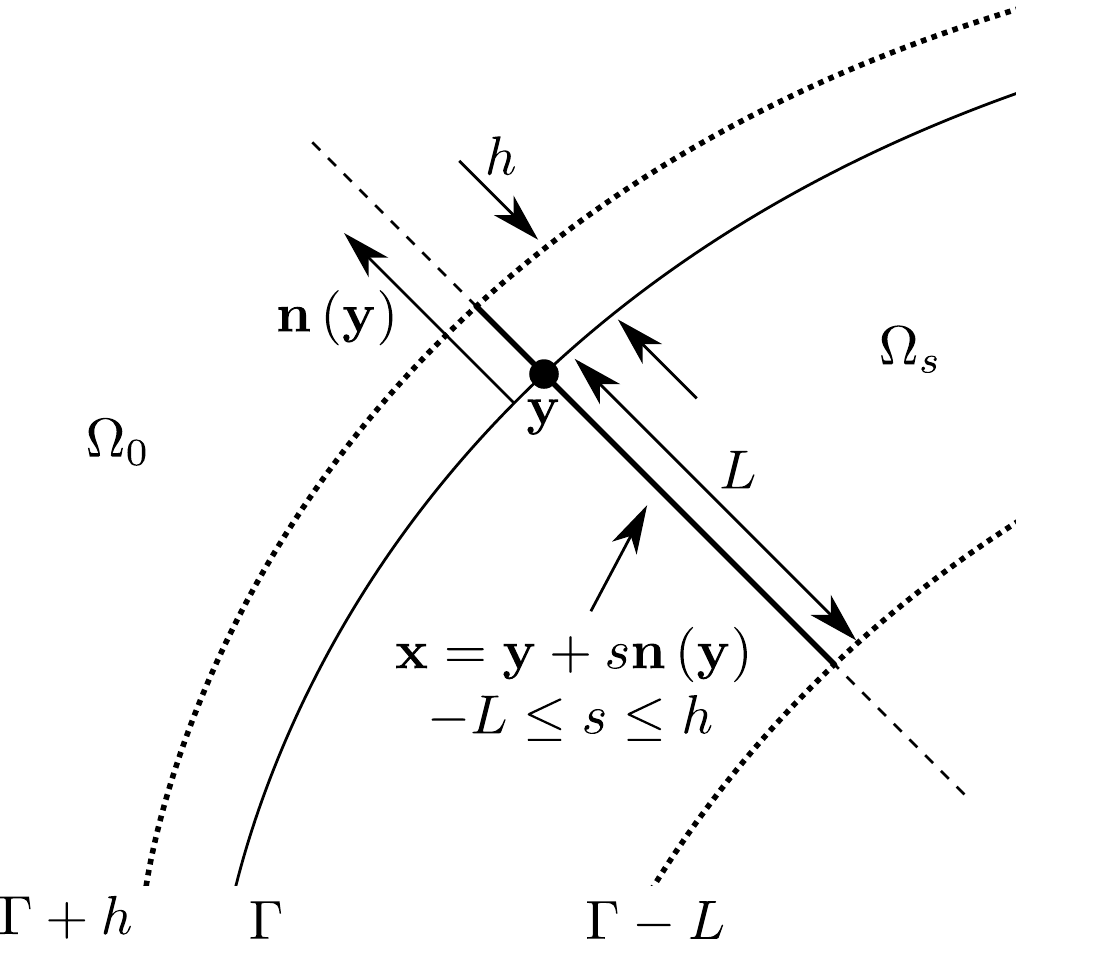} \\
    \caption{The local construction of a ray in the vicinity of $\Gamma$. The function $\tilde{\mathbf{g}}$ is built as a polynomial along each ray.} \label{Fig_TwoD_Construction}
\end{figure}

The extension is constructed as follows:
\begin{description}
\item [Step 1] \hfill \\
Build a local coordinate system surrounding the interface $\Gamma$ in a region between $\Gamma - L$ and $\Gamma + h$ (see Figure \ref{Fig_TwoD_Construction}). 
Suppose $\mathbf{n}$ is the outward normal at $\mathbf{y} \in \Gamma$. Then one can write a local system $(\mathbf{y}, s)$ defined implicitly by
	\begin{align} \label{LocalCoordinates}
		\mathbf{x} &= \mathbf{y} + s \mathbf{n}(\mathbf{y}) 
	\end{align}
	where $-L \leq s \leq h$ and $\mathbf{y} \in \Gamma$. 	For $\Gamma \in C^2$ and sufficiently small $L$ and $h$, one can always invert (\ref{LocalCoordinates}) so that $s = s(\mathbf{x})$ and $\mathbf{y} = \mathbf{y}(\mathbf{x})$ are functions of the coordinates $\mathbf{x}$. 

\item [Step 2] \hfill \\
Build $\tilde{g}_z(\mathbf{x})$ using one-dimensional polynomials along rays. 
	Given a point $\mathbf{x}$ between $\Gamma - L$ and $\Gamma + h$, along with the corresponding point $\mathbf{y} = \mathbf{y}(\mathbf{x})$ and distance $s = s(\mathbf{x})$ from Step 1, the extension is
	\begin{enumerate}
		\item [A.] Matching $m = 0$ derivatives at $\mathbf{x} \in \Gamma + h$, which will yield a $(\Delta x)^{1.5}$ order scheme
	\begin{align}	
		\tilde{g}_z(\mathbf{x}) &= g(\mathbf{y}) \; P_{0,0}(s) + E_{z, \eta}(\mathbf{y} + h \mathbf{n}(\mathbf{y})) \; P_{0, 1}(s)
	\end{align}
		\item [B.] Matching $m = 1$ derivatives at $\mathbf{x} \in \Gamma + h$, which will yield a $(\Delta x)^{2.5}$ order scheme
	\begin{align}	
		\tilde{g}_z(\mathbf{x}) &= g(\mathbf{y}) \; P_{1,0}(s) + E_{z, \eta}(\mathbf{y} + h \mathbf{n}(\mathbf{y})) \; P_{1, 1}(s) + \Big[\frac{\partial E_{z, \eta}}{\partial s} (\mathbf{y} + h \mathbf{n}(\mathbf{y}))\Big] \; P_{1, 2}(s)
	\end{align}
\end{enumerate}
	Note that in the constructions $[A]$ and $[B]$, the point $\mathbf{y} + h \mathbf{n}(\mathbf{y})$ is on $\Gamma + h$. Meanwhile, in construction [B], for small $h\ll 1$ the expression $\frac{\partial E_{z, \eta}}{\partial s}= (\mathbf{n}\cdot \nabla) E_{z, \eta}$ is the derivative of $E_{z, \eta}$ in the normal direction $\mathbf{n}$.
\end{description}

\begin{remark}
	In the case where the interface $\Gamma = \{\mathbf{x} \in \mathbbm{R}^d: \psi(\mathbf{x}) = 0\}$ is described by a level set $\psi$ with $|\nabla \psi| = 1$ and
\begin{align}
	\Omega_0 &= \{\mathbf{x} \in \mathbbm{R}^d : \psi(\mathbf{x}) > 0 \} \\
   \Omega_s &= \{\mathbf{x} \in \mathbbm{R}^d : \psi(\mathbf{x}) < 0 \},
\end{align}
then
\begin{align}
	\Gamma + h &= \{\mathbf{x} \in \mathbbm{R}^d: \psi(\mathbf{x}) = h\} \\
	\Gamma - L &= \{\mathbf{x} \in \mathbbm{R}^d: \psi(\mathbf{x}) = -L\}.
\end{align}
In addition, $\mathbf{n} = \nabla \psi$ represents a local normal to the level sets. 
\myremarkend
\end{remark}
\begin{remark}
	In simple geometries, such as a circular arc, one can explicitly solve equation (\ref{LocalCoordinates}) to recover $(\mathbf{y}, s)$ from $\mathbf{x}$. In cases where the interface $\Gamma$ is described as the zero level set of a function $\psi \in C^2$ so that $\psi(\mathbf{y}) = 0$ for all $\mathbf{y} \in \Gamma$, then one can easily recover $s(\mathbf{x})$ and $\mathbf{y}(\mathbf{x})$ using a Newton iteration method. In this case, $\psi$ does not need to have unit normal ($|\nabla \psi| \neq 1$).  We provide further numerical details in Section \ref{Sec_NumSolveCoordinates}. \myremarkend
\end{remark}
\begin{remark}
	The boundary curvature cannot be infinite. In practice, the curvature of the boundary $\Gamma$ should be small enough to be resolved by the grid spacing of the spatial discretization. \myremarkend
\end{remark}

\subsection{Penalty function $\tilde{\mathbf{g}}$ for a $TE$ mode} \label{Sec_PenaltyFunction_TEmode}

When the initial data and boundary data does not depend on $z$, one also obtains a decoupled $\mathrm{TE}_z$ mode consisting of components $(E_x, E_y, H_z)$. Here we write the boundary data as $\mathbf{g} = (g_x, g_y, 0)$, so that if $\mathbf{n} = (n_x, n_y, 0)$ is a unit normal at any point $\mathbf{y} \in \Gamma$ on the boundary, the boundary condition (\ref{MaxwellEq_BC}) reads
		\[ \big( \mathbf{n} \times (\mathbf{E} - \mathbf{g})\big)_z = n_x(E_y - g_y) - n_y(E_x - g_x) = 0. \]
Here the penalized equations take the form
\begin{subequations} \label{TE_z_Mode}
	\begin{align}
	\frac{\partial E_{x, \eta}}{\partial t} &=	\phantom{-}\frac{\partial H_{z, \eta}}{\partial y}   - \eta^{-1} \chi_h(\mathbf{x})(E_{x, \eta} - \tilde{g}_x), \\
	\frac{\partial E_{y, \eta}}{\partial t} &= -\frac{\partial H_{z, \eta}}{\partial x} - \eta^{-1} \chi_h(\mathbf{x})(E_{y, \eta} - \tilde{g}_y), \\
	\frac{\partial H_{z, \eta}}{\partial t} &= \phantom{-}\frac{\partial E_{x, \eta}}{\partial y} - \frac{\partial E_{y, \eta}}{\partial x}.
	\end{align}
\end{subequations}

\begin{remark}
	It is also possible to penalize only the $H_{z,\eta}$ component of the $\textrm{TE}_z$ mode using the equivalent Neuman PEC boundary condition $(\mathbf{n}\cdot \nabla) H_{z, \eta} = 0$.
\myremarkend
\end{remark}

\begin{remark}
Note that in this case $\nabla\cdot \mathbf{H}_{\eta} = 0$ by virtue of the fact that $\mathbf{H}_{\eta}$ only depends on $z$.  Meanwhile, for a point $\mathbf{x} \in \Omega_0$ (or more precisely outside $(\Gamma + h)$), we may take the divergence of (\ref{TE_z_Mode}) to obtain 
\[ \frac{\partial (\nabla \cdot \mathbf{E}_{\eta})}{\partial t} = 0, \quad \textrm{for } \mathbf{x} \textrm{ outside } (\Gamma + h). \]
Therefore if $\nabla\cdot \mathbf{E}_0 = 0$, then the divergence is preserved to be zero.  In the case where numerical spectral derivatives are used, additional care must be taken to ensure that $\nabla\cdot \mathbf{E}_{\eta}=0$ remains zero. \myremarkend
\end{remark}

The goal is to choose $\tilde{\mathbf{g}}$ to penalize the tangential component of $\mathbf{E}_{\eta} = (E_{x, \eta}, E_{y, \eta}, 0)$ in exactly the same fashion that $\tilde{g}_z$ penalized $E_{z,\eta}$ in the $\mathrm{TM}_z$ mode.  However, since there are now two components of $\tilde{\mathbf{g}}$, we choose a second condition to ensure that $\tilde{\mathbf{g}}$ does not affect the normal component of $\mathbf{E}_{\eta}$ at the boundary. Namely, we choose $\tilde{\mathbf{g}}$ so that the penalty term penalizes only the tangential component of $\mathbf{E}_{\eta}$:
	\begin{align*}
		\mathbf{n} \cdot (\mathbf{E}_{\eta} - \tilde{\mathbf{g}} ) &= 0, \quad \quad \quad \quad \textrm{for } \mathbf{y} \in \Gamma,  \\
		\mathbf{n} \times (\tilde{\mathbf{g}} - \mathbf{g}) &= 0\phantom{,}  \quad \quad \quad \quad \textrm{for } \mathbf{y} \in \Gamma.
	\end{align*}
	The two conditions can be guaranteed provided we take $\tilde{\mathbf{g}}$ to be 
	\begin{align*}
		\tilde{\mathbf{g}}	&= (\mathbf{E}_{\eta} \cdot \mathbf{n})\mathbf{n} + (\mathbf{g} - (\mathbf{g} \cdot \mathbf{n})\mathbf{n} ), \quad \textrm{for } \mathbf{y} \in \Gamma.
	\end{align*}

 To make the construction explicit, let $\mathbf{n} = (n_x, n_y, 0)$ be the normal at any point $\mathbf{y} \in \Gamma$ on the boundary.  Then for any point $\mathbf{x}$ between $-L$ and $h$ of $\Gamma$, we solve (\ref{LocalCoordinates}) to find $\mathbf{y} = \mathbf{y}(\mathbf{x})$, $s = s(\mathbf{x})$ and the corresponding normal $\mathbf{n}(\mathbf{y})$. The components of $\tilde{\mathbf{g}}$ are then constructed in a very similar fashion to $\tilde{g}_z$ for the $\mathrm{TM}_z$ mode.
	\begin{enumerate}
		\item [A.] Matching $m = 0$ derivatives at $\mathbf{x} \in \Gamma + h$, which will yield a $(\Delta x)^{1.5}$ order scheme
		\begin{align*}
		\tilde{\mathbf{g}}(\mathbf{x}) = \big[ (\mathbf{E}_{\eta}(\mathbf{y}) \cdot \mathbf{n})\mathbf{n} + (\mathbf{g}(\mathbf{y}) - (\mathbf{g}(\mathbf{y}) \cdot \mathbf{n})\mathbf{n})
 \big] \ P_{0,0}(s) +  \mathbf{E}_{\eta}(\mathbf{y} + h \mathbf{n}) \; P_{0, 1}(s),
		\end{align*}
		or explicitly in components
		\begin{align} \nonumber
			\tilde{g}_x(\mathbf{x}) &= \big[(E_{x, \eta}(\mathbf{y}) n_x + E_{y, \eta}(\mathbf{y}) n_y) n_x + g_{x}(\mathbf{y}) - (g_{x}(\mathbf{y}) n_x + g_{y}(\mathbf{y}) n_y) n_x \big] \; P_{0,0}(s) \\
		&+ E_{x,\eta}(\mathbf{y} + h \mathbf{n}) \; P_{0, 1}(s), \\ \nonumber
		\tilde{g}_y(\mathbf{x}) &= \big[ (E_{x, \eta}(\mathbf{y})  n_x + E_{y, \eta}(\mathbf{y})  n_y) n_y + g_{y}(\mathbf{y}) - (g_{x}(\mathbf{y}) n_x + g_{y}(\mathbf{y}) n_y) n_y \big] \; P_{0,0}(s) \\
		&+ E_{y, \eta}(\mathbf{y} + h \mathbf{n}) \; P_{0, 1}(s).
		\end{align}
	\end{enumerate}
	It is important to note that the constructions for $\tilde{g}_x$ and $\tilde{g}_y$ are very similar to the construction for $\tilde{g}_z$ in the $\mathrm{TM}_z$ mode (with the exception of having different coefficients for the $P_{0,0}$ term), and are done independently for each of the two components.	 Moreover, in the case of a PEC boundary condition, $\mathbf{g} = 0$ and the penalty term $\tilde{\mathbf{g}}$ only depends on $\mathbf{n}\cdot \mathbf{E}_{\eta}$ at the boundary $\Gamma$.

\subsection{Construction of $\tilde{\mathbf{g}}$ in the general case} \label{Sec_PenaltyFunction_General}

The more general case of constructing the extension $\tilde{\mathbf{g}}$ in higher dimensions builds on the general prescription described in the previous Section \ref{Sec_PenaltyFunction_TEmode} for the $\mathrm{TE}_z$ mode.  In particular, the penalty function $\tilde{\mathbf{g}}$ is chosen to penalize the tangential component of the field $\mathbf{E}_{\eta}$ and to approximately enforce $\mathbf{n}\times(\mathbf{E}_{\eta} - \mathbf{g}) = 0$ on $\Gamma$.  As a result, we take $\tilde{\mathbf{g}}$ to match the exact value of $\mathbf{E}_{\eta}$ at $\Gamma + h$, and also satisfy the tangential component of the boundary condition at $\Gamma$
	\begin{alignat*}{2} \label{ExtensionConditionGeneral}
			\tilde{\mathbf{g}}	&= (\mathbf{E}_{\eta} \cdot \mathbf{n})\mathbf{n} + (\mathbf{g} - (\mathbf{g} \cdot \mathbf{n})\mathbf{n} ), \quad && \textrm{for } \mathbf{y} \in \Gamma, \\
			\tilde{\mathbf{g}} &= \mathbf{E}_{\eta}, \quad && \textrm{for } \mathbf{y} \in \Gamma + h.
	\end{alignat*}
	Again, we may make the construction explicit. First, given any $\mathbf{x}$ within $-L$ and $h$ of $\Gamma$, solve (\ref{LocalCoordinates}) for $\mathbf{y} = \mathbf{y}(\mathbf{x})$ and $s = s(\mathbf{x})$. The extension is then written as
	\begin{enumerate}
		\item [A.] Matching $m = 0$ derivatives at $\mathbf{x} \in \Gamma + h$, which will yield a $(\Delta x)^{1.5}$ order scheme
		\begin{align*}
		\tilde{\mathbf{g}}(\mathbf{x}) = \big[ (\mathbf{E}_{\eta}(\mathbf{y}) \cdot \mathbf{n})\mathbf{n} + (\mathbf{g}(\mathbf{y}) - (\mathbf{g}(\mathbf{y}) \cdot \mathbf{n})\mathbf{n}) \big] \ P_{0,0}(s) +  \mathbf{E}_{\eta}(\mathbf{y} + h \mathbf{n}) \; P_{0, 1}(s).
		\end{align*}
	\end{enumerate}
	\[ \]

\begin{remark}
	While we have only provided the explicit construction for a $(\Delta x)^{1.5}$ scheme in the $\textrm{TE}_z$ and general cases, it is possible to obtain systematically higher rates of convergence analytically by including additional normal derivatives in the construction of $\tilde{\mathbf{g}}$. However, there are difficulties associated with obtaining stable numerical schemes in these cases. \myremarkend
\end{remark}

\section{$\mathrm{TM}_z$ Mode Convergence Analysis} \label{AnalyticTM_Analysis}
In this section we examine the analytic convergence rate in $\eta$ and $h$ for plane wave scattering solutions off of a flat wall for the two-dimensional Maxwell's equations.  The problem of a plane wave $\mathrm{TM}_z$ mode defined on $x > 0$ scattering off a flat wall at $x = 0$ has a solution of the form $\mathbf{E} = (0, 0, E_{z}(x,y,t))$, $\mathbf{H} = (H_{x}(x, y, t), H_{y}(x, y, t), 0)$. Introducing complex notation, $E_{z} = \Re\{ \hat{E}_z \}$, $H_{x} = \Re\{ \hat{H}_x \}$, $H_{y} = \Re\{ \hat{H}_y \}$, the solution has components
\begin{subequations} \label{Incident_TM_Wave}
	\begin{align} 
	\hat{E}_{z} &= \phantom{-} E_i e^{\imath \omega t} \Big( e^{\imath(k_x, k_y)\cdot(x, y)} - e^{\imath(-k_x, k_y)\cdot(x, y)} \Big) \\
	\hat{H}_{x} &= -E_i \frac{k_y}{\omega} e^{\imath \omega t} \Big( e^{\imath(k_x, k_y)\cdot(x, y)} - e^{\imath(-k_x, k_y)\cdot(x, y)} \Big) \\
	\hat{H}_{y} &= \phantom{-} E_i \frac{k_x}{\omega} e^{\imath \omega t} \Big( e^{\imath(k_x, k_y)\cdot(x, y)} + e^{\imath(-k_x, k_y)\cdot(x, y)} \Big), 
\end{align}
\end{subequations}
where $E_i$ is the amplitude of the incoming wave.  The PEC boundary condition $\hat{E}_z(0, y, t) = 0$ at $x = 0$ forces all of the incoming wave to be reflected back. Here we have introduced $\omega^2 = k_x^2 + k_y^2 \neq 0$ as the dispersion relation. We now examine the error associated with a solution to the penalized equations (\ref{TM_z_Mode}) containing the same incident wave as (\ref{Incident_TM_Wave}) with
\begin{eqnarray}
	\chi_h(x,y) = \begin{cases}
		0, & x > h \\
		1, & x \leq h
	\end{cases}.
\end{eqnarray}
In this model problem, we take $\tilde{g}_z(x, y)$ to match the function value of $E_{z, \eta}$ at $x = h$ with, for analytic simplicity, a lower order polynomial than in (\ref{Extension_dim1_m0}) given by
\begin{displaymath}
	\tilde{g}_z(x, y) = \frac{ x(x+1)}{h(h+1)} \; E_{z, \eta}(h, y), \hspace{10mm} x \leq h.
\end{displaymath}

\begin{remark}
	For simplicity, in this example we are interested in quantifying the analytic error induced by the improved penalty term. As a result, we take the simplest function $\tilde{g}_z(x,y)$ to be a low order polynomial which vanishes at $x = -1$ inside $\Omega_s$. In practice, numerical implementations for $\tilde{g}_z(x,y)$ require $\tilde{g}_z$ to vanish more smoothly at $x = -1$ inside the obstacle region, as to avoid oscillations in the Fourier representation of $E_{z, \eta}$.
	\myremarkend
\end{remark}

We now solve the penalized equations for an incoming wave with amplitude $E_i$ and determine the error in the reflection due to the penalty term.  For $x > h$, the penalized equations reduce to the Maxwell's equations in free space and we can write a general solution as
\begin{align} \label{E_Incident}
	\hat{E}_{z, \eta} &= \phantom{-} E_i e^{\imath \omega t} \Big( e^{\imath(k_x, k_y)\cdot(x-h, y)} - R e^{\imath(-k_x, k_y)\cdot(x-h, y)} \Big) \\
	\hat{H}_{x, \eta} &= -E_i \frac{k_y}{\omega} e^{\imath \omega t} \Big( e^{\imath(k_x, k_y)\cdot(x-h, y)} - R e^{\imath(-k_x, k_y)\cdot(x-h, y)} \Big) \\
	\hat{H}_{y, \eta} &= \phantom{-} E_i \frac{k_x}{\omega} e^{\imath \omega t} \Big( e^{\imath(k_x, k_y)\cdot(x-h, y)} + R e^{\imath(-k_x, k_y)\cdot(x-h, y)} \Big) 
\end{align}
where again $E_i$ is the amplitude of the incoming wave, while $R$ is a reflection coefficient to be determined by matching the solution across the penalty region.  We note that in the exact unpenalized problem, $R_{exact} = e^{-2\imath h k_x}$.

For $x < h$, we may use separation of variables and write $E_{z,\eta} = \Re\{ e^{\imath \omega t + \imath k_y y} \hat{E}(x) \}$. Using the ansatz for $E_{z,\eta}$, we obtain 
\[ H_{x, \eta} = \Re\Big\{ -\frac{k_y}{\omega} \hat{E}(x) e^{\imath \omega t + \imath k_y y} \Big\}, \hspace{10mm} H_{y, \eta} = \Re \Big\{ \frac{1}{\imath \omega} \; \frac{d \hat{E}}{d x} \; e^{\imath \omega t + \imath k_y y} \Big\}, \]
along with an ODE obtained from (\ref{TM_z_Mode}) for $\hat{E}(x)$
\begin{eqnarray} \label{ODE_TM_Penalty}
	\frac{d^2 \hat{E}}{d x^2} + (\omega^2 - k_y^2) \hat{E} - \frac{\imath \omega}{\eta}(\hat{E} - \hat{g}(x) ) = 0, \hspace{10mm} x \leq h
\end{eqnarray}
where
\[ \hat{g}(x) = \frac{x(x+1)}{h(h+1)} \hat{E}(h). \]
The ODE can then be simplified into the following form
\begin{equation} \label{simplified_ODE}
	\frac{d^2 \hat{E}}{dx^2} + \gamma^2 \hat{E} + A (x + x^2) = 0
\end{equation}
where
\[ \gamma^2 = \omega^2 - k_y^2 - \imath \omega \eta^{-1}\]
and 
\[ A = \imath \omega \eta^{-1} \frac{1}{h(h+1)} \hat{E}(h). \]
On $x \leq h$, the ODE (\ref{simplified_ODE}) has the solution
\begin{equation}
	\hat{E}(x) = E_0 e^{\imath \, \gamma (x - h)} + A[ 2\gamma^{-4} - \gamma^{-2} x - \gamma^{-2} x^2]
\end{equation}
where $\gamma$ is the unique root with $\Re \{\imath \gamma\} \geq 0$ chosen to satisfy the radiation condition (exponential decay) for $x \rightarrow -\infty$, and $E_0$ is the constant of integration.  

We now solve for the two unknowns $E_0$ and $R$ by imposing continuity of the solution $E_{z, \eta}$ and $H_{x, \eta}$, $H_{y, \eta}$ at $x = h$ for all $y$. Continuity of $\hat{E}_{z, \eta}$ at $x = h$ yields
\begin{eqnarray}
	E_0 + A[2\gamma^{-4} - \gamma^{-2} h - \gamma^{-2} h^2] = E_i (1 - R) 
\end{eqnarray}
while continuity of $\hat{H}_{y, \eta}$ yields
\begin{eqnarray}
	\imath \gamma E_0 + A[-\gamma^{-2} - 2 \gamma^{-2} h] = \imath k_x E_i (1 + R)
\end{eqnarray}
where $A = \imath \omega \eta^{-1} \frac{E_i}{h(h+1)}(1 - R)$.  The two equations can be used to find $E_0$ in terms of $E_i$ and $R$ as functions of $h$ and $\eta$. Specifically, we can eliminate $E_0$ and write $E_r = R E_i$ for some reflection coefficient $R$. Solving for $R$ (via MATLAB's symbolic package) and expanding in powers of $\eta$ and $h$ yields
\begin{align}
	R &\approx (1 - 2\imath h k_x) - [2\sqrt{2}(1+\imath)k_x \omega^{-1/2}] \sqrt{\eta}h + [2 k_x \omega^{-1} (k_x^2 - 4)]\eta h + \mathcal{O}(h^2).
\end{align}

The leading term $(1-2\imath h k_x)$ is exactly the first order term in the reflection coefficient $R_{exact}$. Therefore, $R$ and $R_{exact}$ differ by order $\sqrt{\eta}h$.  Fixing $\eta = O(\Delta x)$ and $h = O(\Delta x)$ yields an error of $O(\Delta x^{1.5})$ in both the amplitude $|R|$ and phase $\angle R$. Hence, we have a global error of order $1.5$. 

\begin{remark}
	Errors at both order $O(\Delta x^{1.5})$ and $O(\Delta x^2)$ appear in the expansion for $|R|$ and $\angle R$.  Hence, one may initially see $2$nd order convergence before observing the asymptotic convergence rate of $1.5$. 
	\myremarkend
\end{remark}

\begin{remark}
	One can repeat the calculation in this section by taking a static, non-active volume penalty term of the form $\eta^{-1} \chi_0(\mathbf{x}) \, \mathbf{E}_{\eta}$, where $h = 0$ and $\tilde{\mathbf{g}} = 0$.  Such a choice for a non-active penalty term recovers the PEC boundary conditions, however results in a slow analytic convergence rate (\ref{ConvergenceLimit_E}) of $O(\eta^{1/2})$.  For numerical purposes, such an analytic convergence rate translates into a numerical scheme with global convergence $O(\Delta x^{1/2})$.    \myremarkend
\end{remark}

\section{Perfectly Matched Layers (PML)} \label{sec_perfectly_matched_layers}
In our current approach using Fourier methods, we work on a rectangular domain with periodic boundary conditions.  In many applications, however, one is not interested in solving Maxwell's equations in a periodic domain, but rather on an infinite one.  One major difficulty which arises when using a periodic computational domain to compute solutions on an infinite one is the artificial wrapping of traveling waves.  Namely, waves which should radiate out on an infinite domain simply wrap back into the computational domain as a result of the periodic boundary conditions.  In this section, we outline how to eliminate the artificial wrapping so that one may compute time-domain radiating solutions, such as those arising from scattering problems, on an effective infinite domain.  The approach is through the introduction of a perfectly matched layer (PML) \cite{Berenger1994}.  Although PMLs were originally introduced to eliminate artificial reflections which arise from a finite truncation of a computational domain, they are easily modified to the case of a periodic domain.

Here we outline how to modify the PML from a square domain with Dirichlet boundary conditions, to a periodic one.  We do so for the case of a $\mathrm{TE}_z$ mode and note that the modification closely follows the formulation originally proposed in \cite{Berenger1994}.  

As a first step, we decompose the field $H_{z, \eta} = H_{zx, \eta} + H_{zy, \eta}$ into two components.  In the absence of a PML, we choose the decomposition so that the two components evolve according to
\begin{subequations}
	\label{PML_Hz} 
	\begin{align}
		\frac{\partial H_{zx, \eta}}{\partial t} &= -\frac{\partial E_{y,\eta}}{\partial x}  \\
		\frac{\partial H_{zy, \eta}}{\partial t} &= \phantom{-}\frac{\partial E_{x,\eta}}{\partial y}.
	\end{align}
\end{subequations}
Although the addition of an extra equation appears redundant, the decomposition simplifies the resulting implementation of a PML.  To add a PML we further modify the extended $\mathrm{TE}_z$ equations (\ref{PML_Hz}) to contain an absorbing layer.  For the absorbing layer, we let $(\sigma_x, \sigma_y)$ denote two effective material parameters.  Since the domain is periodic, we simply choose the PML to have two bands, one vertical and one horizontal.  For example, Figure \ref{Fig_TE_scattering} shows a periodic domain with two such strips outlined by dashed lines.  Outside of each strip, we take $(\sigma_x, \sigma_y) = (0,0)$ as a physical domain which allows for the normal propagation of the $\mathrm{TE}_z$ mode.  In such a region, one may have curved obstacles.  In the PML region, we take $\sigma_x, \sigma_y > 0$ and modify the $\mathrm{TE}_z$ equations as follows:
\begin{subequations}
	\label{perfectly_matched_layers}
	\begin{align} \label{TE_PML_1}
	\frac{\partial E_{x, \eta}}{\partial t} &=  \phantom{-}\frac{\partial}{\partial y} ( H_{zx, \eta} + H_{zy, \eta}) - \sigma_y E_{x, \eta} - \eta^{-1} \chi(\mathbf{x}) (E_{x, \eta} - g_x)\\ \label{TE_PML_2}
	\frac{\partial E_{y, \eta}}{\partial t}  &= -	\frac{\partial}{\partial x} ( H_{zx, \eta} + H_{zy, \eta}) - \sigma_x E_{y, \eta} - \eta^{-1} \chi(\mathbf{x}) (E_{y, \eta} - g_y) \\ \label{TE_PML_3}
	\frac{\partial H_{zx, \eta}}{\partial t} &= -	\frac{\partial}{\partial x} E_{y,\eta} - \sigma_x H_{zx,\eta} \\ \label{TE_PML_4}
	\frac{\partial H_{zy, \eta}}{\partial t} &= \phantom{-}\frac{\partial}{\partial y}E_{x,\eta} - \sigma_y H_{zy, \eta}.
	\end{align}
\end{subequations}
Through direct calculation, \cite{Berenger1994} showed that such a modification\footnote{In the case where $\epsilon_0$ and $\mu_0$ are not $1$, one must rescale the coefficients $\sigma_x$ and $\sigma_y$ in equations (\ref{perfectly_matched_layers}). For example, see equations $(2)$--$(3)$ in \cite{Berenger1994}.} results in a perfectly matched layer.  Specifically, a wave traveling from the region where $(\sigma_x, \sigma_y) = (0, 0)$ does not reflect off the region where $\sigma_x$ (or $\sigma_y$) is non-negative regardless of the incident angle or frequency.  Although discontinuous jumps in $\sigma_x$ (or $\sigma_y$), do not theoretically reflect waves in a PML, they can result in numerical reflections when computing a numerical solution.  As a result, in practice, we choose $\sigma_x$ (or $\sigma_y$) to grow linearly up to a maximum value $\sigma_{x, max}$ (or  $\sigma_{y, max}$).  Here, the slope and maximum value may depend on the exact problem. In practice, one can ramp up to the maximum value over a few wavelengths.

\section{A Numerical Fourier Algorithm} \label{Sec_NumFourierAlgorithm}

In this section we outline the numerical method, and details we use when solving the penalized Maxwell's equations.

Let $\Omega = [0, D]^d$ be the domain. Then in our scheme, we use an equispaced grid with $N$ (even) points, and spacing $\Delta x = D/N$. In two dimensions we take $\Delta y = \Delta x$, however one does not in general require equal grid spacing.  Grid points are denoted as 
	\begin{align*}
		x_j &= j \Delta x, \quad \quad y_j = j \Delta y, \quad \quad \textrm{for } j = 0, \ldots, N-1, 
	\end{align*}
	and variables evaluated at gridpoints as $u_j = u(x_j)$.
	
	In the numerics, we also make use of the discrete Fourier transform of a function $u(x)$ (on a domain of length $D$)
\begin{equation}
	\hat{u}_l = \mathcal{F}\left\{ u \right\} = \sum_{j = 0}^{N-1} u_j \; e^{-\imath k_l x_j},
\end{equation}
where 
\begin{alignat*}{2}
	k_l &= \frac{2\pi l}{D} \quad && \textrm{for } 0 \leq l \leq N/2 \\
	k_l &= \frac{2\pi (l - N)}{D} \quad && \textrm{for } N/2 + 1 \leq l \leq N-1
\end{alignat*}	
	are the wavenumbers.
The inverse is then taken as 
\begin{equation}
	u_j = \mathcal{F}^{-1}\{ \hat{u} \} = \frac{1}{N} \sum_{l = 0}^{N-1} \hat{u}_l \; e^{\imath k_l x_j}.
\end{equation}
The discrete Fourier transform pairs also have natural extensions to higher dimensions.

\subsection{Time-stepping Details}

Linear wave equations, such as Maxwell's equations, have an evolution governed by operators with purely imaginary eigenvalues.  As a result, explicit time-stepping schemes may not be stable if the stability region does not incorporate a sufficient portion of the imaginary axis. The purpose of this section is to present the stability results for standard Runge-Kutta time stepping schemes using Fourier spectral differentiation in space in the abscence of penalization. 

As an example, we consider a Fourier method for the one-dimensional Maxwell's equations on a periodic domain of $D = 2\pi$, given by
\begin{subequations}	\label{Maxwell_Fourier}
\begin{align} 
	\label{OneD_Fourier1}
	\frac{\partial E_z}{\partial t} &= \mathcal{F}^{-1} \left\{ \, \imath k \, \mathcal{F} \left\{ H_y \right\} \right\}, \\ \label{OneD_Fourier2}
	\frac{\partial H_y}{\partial t} &= \mathcal{F}^{-1} \left\{ \, \imath k \, \mathcal{F}\left\{ E_z \right\} \right\}.
\end{align}
\end{subequations}

We report stability requirements for common spectral time-stepping schemes to (\ref{Maxwell_Fourier}) by listing the eigenvalues $\lambda$ to the discrete linear time evolution in Table \ref{Table_TimeStepStability}.  Here the eigenvalue amplitude $|\lambda|^2 < 1$ is required for stability.  We denote 
\begin{eqnarray}
	r = \Delta t \, k \in \mathbbm{R}
\end{eqnarray}
as the real parameter which combines the time step $\Delta t$ and wavenumber $k$.  
\begin{table}
	\centering
	\caption{Stability for standard time-stepping schemes}
	\label{Table_TimeStepStability}
   \begin{tabular}{ |l | c | r |} 
     \hline
     Integration of (\ref{OneD_Fourier1})--(\ref{OneD_Fourier2}) & Eigenvalue amplitude $|\lambda|^2$ & Stability \\ \hline \hline
     Simple Euler & $1 + r^2$ & Unstable \\ \hline
     Modified Euler (RK2) & $1 + \frac{1}{4}r^4$ & Unstable \\
     \hline
    $4$th order Runge-Kutta (RK4) & $1 - \frac{1}{72}r^6 + \frac{1}{576}r^8$ & Stable for $r < 2.83$ \\
	\hline
	 Implicit Euler & $(1 + r^2)^{-1}$ & Unconditionally stable \\  \hline
	 Integrating factor \cite{MilewskiTabak1999} & 1 & Unconditionally stable \\ \hline
   \end{tabular}
\end{table}
As outlined in Table \ref{Table_TimeStepStability}, the simple Euler and Modified Euler schemes are always unstable.  Meanwhile, RK4 is stable provided $r < \sqrt{576/72} \sim 2.83$.  For a $d$-dimensional periodic square with side length $D$, one then has  $k_{max} =  \pi \sqrt{d} (N / D)$ with
\begin{align}
	\Delta t &< 2.83 k_{max}^{-1} \\
	&= \frac{2.83}{\pi \sqrt{d}}\frac{D}{N}.
\end{align}

\subsection{Solving equation (\ref{LocalCoordinates}) for the local coordinates} \label{Sec_NumSolveCoordinates}

In many applications one describes the boundary or interface as the zero level set of a function $\hat{\psi}(\mathbf{x})$.  In this section we provide some brief numerical details on how one can numerically use the level set (which may not have unit norm) to build the local coordinate system $\mathbf{y}=\mathbf{y}(\mathbf{x})$ and $s = s(\mathbf{x})$ as the solution to equation (\ref{LocalCoordinates}).

 We accomplish this numerically using a damped Newton method, as described in \cite{Persson2005}.
That is, for every grid point $\mathbf{x}_j$ between $(\Gamma - L)$ and $(\Gamma + h)$, we seek the point $\mathbf{y}_j$ on the zero level set of $\hat{\psi}$ (i.e., the interface) such
that $\mathbf{y}_j - \mathbf{x}_j$ is parallel to the normal direction $\mathbf{n}(\mathbf{y}_j) \propto \nabla\hat{\psi}\left(\mathbf{y}_j\right)$.
That is, we would like 
\begin{equation} \label{NewtonRoot}
f\left(\mathbf{y}_j\right)=\left[\begin{array}{c}
\hat{\psi}\left(\mathbf{y}_j\right)\\
(\mathbf{y}_j - \mathbf{x}_j)\times\nabla\hat{\psi}\left(\mathbf{y}_j\right)
\end{array}\right]=0.
\end{equation}
Here one could also arrive at equation (\ref{NewtonRoot}) by dotting and crossing (\ref{LocalCoordinates}) through with $\nabla \hat{\psi}$ since it is proportional to $\mathbf{n}$. We note that for a Newton iteration to work, we assume that $\hat{\psi} \in C^2$ locally near the interface $\Gamma$ so that one may compute the Jacobian of (\ref{NewtonRoot}).

Once we have the point $\mathbf{y}_j$ corresponding to each grid point $\mathbf{x}_j$, we compute
\begin{equation}
\psi\left(\mathbf{x}_j\right)=\mathrm{sign}(\hat{\psi}\left(\mathbf{x}_j\right))\left\Vert \mathbf{y}_j - \mathbf{x}_j \right\Vert _{2}
\end{equation}
on the grid which is now a level set function with unit norm $|\nabla \psi| = 1$.  In addition we take $s(\mathbf{x}_j) = \psi(\mathbf{x}_j)$ and the normal used in the local ray construction at each point is simply $\mathbf{n}(\mathbf{y}_j) = \nabla\hat{\psi}\left(\mathbf{y}_j\right)$.

\subsection{Main algorithm and details} \label{sec_main_algorithm}
	\begin{itemize}
		\item Discretize the spatial derivatives in (\ref{Penalized_Maxwell}) using pseudospectral differentiation so that
		\begin{subequations} \label{SpectralMaxwell}
		   \begin{align} 
			     \frac{\partial \mathbf{H}_{\eta}}{\partial t}  &=  - \mathcal{F}^{-1}\{ \imath \mathbf{k} \times \mathcal{F} \{\mathbf{E}_{\eta}\} \}, \\		     
     \frac{\partial \mathbf{E}_{\eta}}{\partial t}  &= \phantom{-}  \mathcal{F}^{-1}\{ \imath \mathbf{k} \times \mathcal{F} \{ \mathbf{H}_{\eta} \}\} - \eta^{-1} \chi_h(\mathbf{x})(\mathbf{E}_{\eta} - \tilde{\mathbf{g}} ).
   			\end{align}
   		\end{subequations}
   		Here, the FFT is used to compute the discrete Fourier transform (and its inverse) on the right hand side (RHS) of (\ref{SpectralMaxwell}).
		\item Compute $\tilde{\mathbf{g}}$ on the RHS of equation (\ref{SpectralMaxwell}):
   			\begin{enumerate}
   				\item[(a)]  Build and store the local coordinate system $s_j = s(\mathbf{x}_j)$, $\mathbf{y}_j = \mathbf{y}(\mathbf{x}_j)$ and normal $\mathbf{n}_j = \mathbf{n}(\mathbf{y}_j)$. \\
Do so for all grid points $\mathbf{x}_j$ between $\Gamma - L$ and $\Gamma + h$. This may typically be done only once.  If required, a Newton iteration with a level set may be used to solve equations (\ref{LocalCoordinates}) for each $\mathbf{x}_j$. 
				\item[(b)]  Compute $\mathbf{E}_{\eta}$ at $\Gamma$ and $\Gamma + h$, i.e., at the points $\mathbf{y}_j$ and $\mathbf{y}_j + h \mathbf{n}_j$, respectively. \\
				Interpolate (via cubic interpolation) the values of $\mathbf{E}_{\eta}$ at the points $\mathbf{y}_j$ (only required for the $\mathrm{TE}_z$ or full three-dimensional cases) and $\mathbf{y}_j + h \mathbf{n}_j$ using the values of $\mathbf{E}_{\eta}$ at the equispaced gridpoints $\mathbf{x}_j$
					\begin{align*}
						\mathbf{E}_{\eta}(\mathbf{y}_j) &\leftarrow \textrm{Interpolate}(\mathbf{E}_{\eta}) \\
						\mathbf{E}_{\eta}(\mathbf{y}_j + h \mathbf{n}_j) &\leftarrow \textrm{Interpolate}(\mathbf{E}_{\eta})
					\end{align*}
				\item[(c)] Compute derivatives of $\mathbf{E}_{\eta}$ at $\Gamma + h$, i.e., at the points $\mathbf{y}_j + h \mathbf{n}_j$.		
			\begin{enumerate}
				\item[(i)] For dimension one: obtain approximate derivatives $E'$, $E''$ on the regular grid as
					\begin{align*}
						E_{z, \eta}' &= \mathcal{F}^{-1} \{ e^{-c_f k^2 /N^2} \imath k \mathcal{F}\{ E_{z, \eta} \}\} \\
						E_{z, \eta}'' &= \mathcal{F}^{-1} \{ e^{-c_f k^2 /N^2} ( \imath k)^2 \mathcal{F}\{ E_{z, \eta} \}\}
					\end{align*}
					where $c_f = 16$ is a high frequency filtering parameter, followed by interpolation to $y_j + h \mathbf{n}_j$ (where $\mathbf{n}_j = \pm 1$ in one dimension)
					\begin{align*}
					  E_{z, \eta}'(y_j+ h \mathbf{n}_j) &\leftarrow \textrm{Interpolate}(E_{z, \eta}')\\
					  E_{z, \eta}''(y_j+ h \mathbf{n}_j) &\leftarrow \textrm{Interpolate}(E_{z, \eta}'').
					 \end{align*}					  
				\item[(ii)]For dimension two: obtain the required derivatives of $\mathbf{E}_{\eta}$,  $\frac{\partial \mathbf{E}_{\eta}}{\partial x}$ and  $\frac{\partial \mathbf{E}_{\eta}}{\partial y}$ on the grid
					\begin{align*}
						\frac{\partial \mathbf{E}_{\eta}}{\partial x} &= \mathcal{F}^{-1} \{ e^{-c_f k^2 /N^2} \imath k_x \mathcal{F}\{ \mathbf{E}_{\eta} \}\} \\
						\frac{\partial \mathbf{E}_{\eta}}{\partial y} &= \mathcal{F}^{-1} \{ e^{-c_f k^2 /N^2} \imath k_y \mathcal{F}\{ \mathbf{E}_{\eta} \}\}  
					\end{align*} 	
					where $c_f = 16$ is a high frequency filtering parameter. Interpolate the derivatives to the points $\mathbf{y}_j + h \mathbf{n}_j \in \Gamma + h$ as
					\begin{align*}
					  \frac{\partial \mathbf{E}_{\eta}}{\partial x}(\mathbf{y}_j + h\mathbf{n}_j) &\leftarrow \textrm{Interpolate}\Big(\frac{\partial \mathbf{E}_{\eta}}{\partial x}\Big)\\
					  \frac{\partial \mathbf{E}_{\eta}}{\partial y}(\mathbf{y}_j + h\mathbf{n}_j) &\leftarrow \textrm{Interpolate}\Big(\frac{\partial \mathbf{E}_{\eta}}{\partial y}\Big).
					\end{align*}		
			\end{enumerate}
	\item[(d)] Using the interpolated values at $\mathbf{y}_j$ and $\mathbf{y}_j + h \mathbf{n}_j$, use the formula [A], [B], or [C] from Sections \ref{Sec_PenaltyFunction_dim1}--\ref{Sec_PenaltyFunction_General} to build $\tilde{\mathbf{g}}(\mathbf{x}_j)$. 
\end{enumerate} 
		\item Evolve (\ref{SpectralMaxwell}) forward in time by $\Delta t$ using RK4 time stepping.
   		\item Due to the spectral derivatives in (\ref{SpectralMaxwell}), when solving either the $\mathrm{TE}_z$ mode or the full equations (\ref{Penalized_Maxwell}), a small non-zero amplitude for $\nabla \cdot \mathbf{E}_{\eta}$ may arise after the 4 stages of RK4 (this does not occur in dimension one or for the $\mathrm{TM}_z$ mode when only $E_z$ appears in the equations). Thus, after the 4 stages of RK4, project out the small divergence of $\mathbf{E}_{\eta}$ by computing
			\begin{align}
				\nabla \cdot \mathbf{E}_{\eta} &= \mathcal{F}^{-1} \{ \imath \mathbf{k} \cdot \mathcal{F} \{ \mathbf{E}_{\eta} \} \}, \\
				\hat{p}(\mathbf{k}) &= -\frac{1}{k^2} \mathcal{F}\{ \; (\nabla\cdot \mathbf{E}_{\eta}) (1 - \chi_h(\mathbf{x})) \; \}, \quad \textrm{for } k = |\mathbf{k}| \neq 0, \\ 
				\mathbf{E}_{\eta} &\leftarrow \mathbf{E}_{\eta} - \mathcal{F}^{-1}\{\imath \mathbf{k} \; \hat{p}(\mathbf{k}) \}.
			\end{align}   		
	\end{itemize}

	\begin{remark}
		Note that in part (c) of our numerical algorithm, we apply a high frequency filter to obtain approximate derivatives for $\mathbf{E}_{\eta}$ in the construction of $\tilde{\mathbf{g}}$.  It is important to note that at no point do we filter the actual solution $\mathbf{E}_{\eta}$ as such a procedure would destroy the accuracy of the algorithm.  Instead, filtering $\mathbf{E}_{\eta}$ in the construction of $\tilde{\mathbf{g}}$ only slightly modifies the penalty forcing term.  The parameter $c_f = 16$ was chosen to ensure stability of the numerical algorithm, while remaining small enough to preserve the overall accuracy of the method. \myremarkend
	\end{remark}
	
\begin{remark}
		Note that each step of our approach makes use of well established algorithms. When the interface $\Gamma$ is described by a level set, the method only requires the FFT, Newton's method, bicubic interpolation, and RK4 time stepping so that implementing the method is straightforward given standard robust numerical packages.
\end{remark}

\section{Stability}

In this section we discuss the stability of the numerical method in Section \ref{Sec_NumFourierAlgorithm}, as well as the  stability of the underlying penalty PDE.  Specifically, we note that there are two separate stability issues to consider. The first is the analytic stability effects that the penalty term has on the underlying solution, while the second is the conventional numerical stability of the Fourier algorithm.  

\subsection{Energy and analytic stability}

\begin{figure}[tb!] 
	\centering
    \includegraphics[width = 1\textwidth]{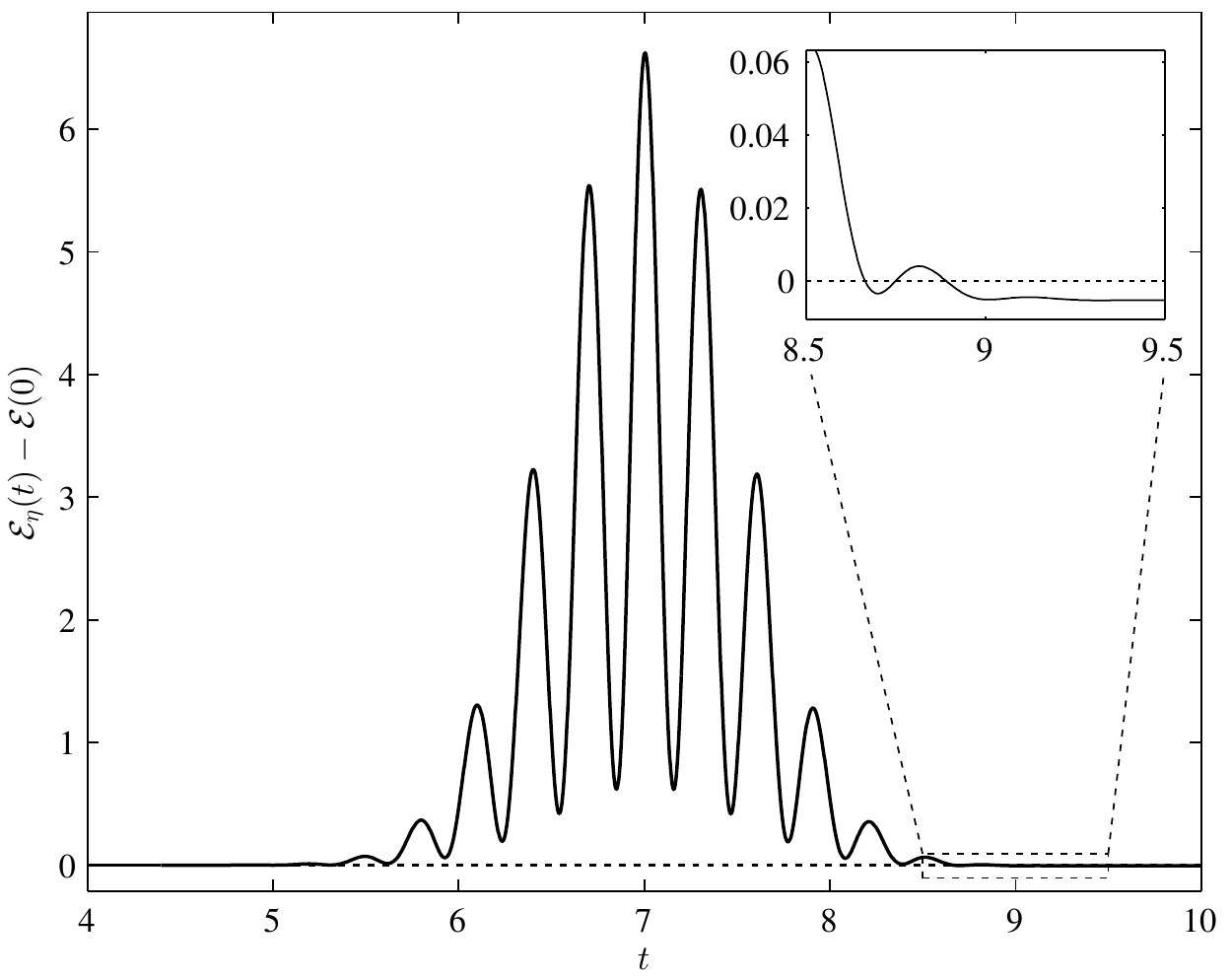} \\
    \caption{Computation of the total energy $\frac{1}{2}\int_{\Omega} |\mathbf{E}_{\eta}|^2 + |\mathbf{B}_{\eta}|^2 \du \mathbf{x} - \mathcal{E}(0) $ for a modulated Gaussian wave in dimension one.  The penalty term forces the solution in the extended region, thereby adding energy to the total system.  The energy in the physical domain $\Omega_0$, remains close to the initial energy in the system. After the wave reflects, the total energy is slightly less than the initial energy. See Section \ref{sec_GaussianScattering1d} for details associated with the specific scattering problem whose energy is depicted.}
 \label{AddedEnergy}
\end{figure}

In domains with PEC boundary conditions, the underlying Maxwell's equations (\ref{MaxwellEq}) conserve the quadratic energy 
\[ \mathcal{E} = \frac{1}{2} \int_{\Omega_0} |\mathbf{E}|^2 + |\mathbf{B}|^2 \du \mathbf{x}. \]
In other words, $\mathcal{E}$ does not depend on time. In the case of the penalized equations, the associated energy of the penalty field 
\[ \mathcal{E}_{\eta} = \frac{1}{2} \int_{\Omega_0} |\mathbf{E}_{\eta}|^2 + |\mathbf{B}_{\eta}|^2 \du \mathbf{x} \] 
is not exactly conserved due to the fact that $\mathbf{E}_{\eta}$ no longer satisfies the exact PEC boundary conditions.  Since $\mathbf{E}_{\eta}$ is close to the exact field $\mathbf{E}$, the energy $\mathcal{E}_{\eta} = \mathcal{E} + O(\eta^{\gamma})$ for the appropriate $\gamma$ corresponding to the convergence rate of the method.  It is important to note that in general $\mathcal{E}_{\eta}$ could be larger (or smaller) than $\mathcal{E}$, which has the interpretation of the penalty term pumping (or removing) a small energy into the reflected fields (see Figure \ref{AddedEnergy}). As a result of the small increase in energy, there can be an associated weakly unstable eigenvalue to the penalized equations.  Note that this eigenvalue can occur at the analytic level and is independent of numerical implementation details.  Numerically the small increase in energy is not problematic since the errors are on the order of the numerical method.  We also note that PMLs have similar behavior reported in the literature \cite{AbarbanelGottlieb1997, AbarbanelGottliebHesthaven2002, BecachePetropoulosGednev2004, DuruKreiss2014} whose study is ongoing.

\begin{remark}
Some numerical experiments based on varying $\tilde{\mathbf{g}}$ suggest that different formulations of the penalty term may act to increase or decrease the small energy difference in the penalized energy $\mathcal{E}_{\eta}$ with respect to $\mathcal{E}$.  We intend to investigate the differences in future work.  \myremarkend
\end{remark}

\subsection{Stability of the numerical scheme}

\begin{figure}[htb!] 
	\centering
    \includegraphics[width = 0.45\textwidth]{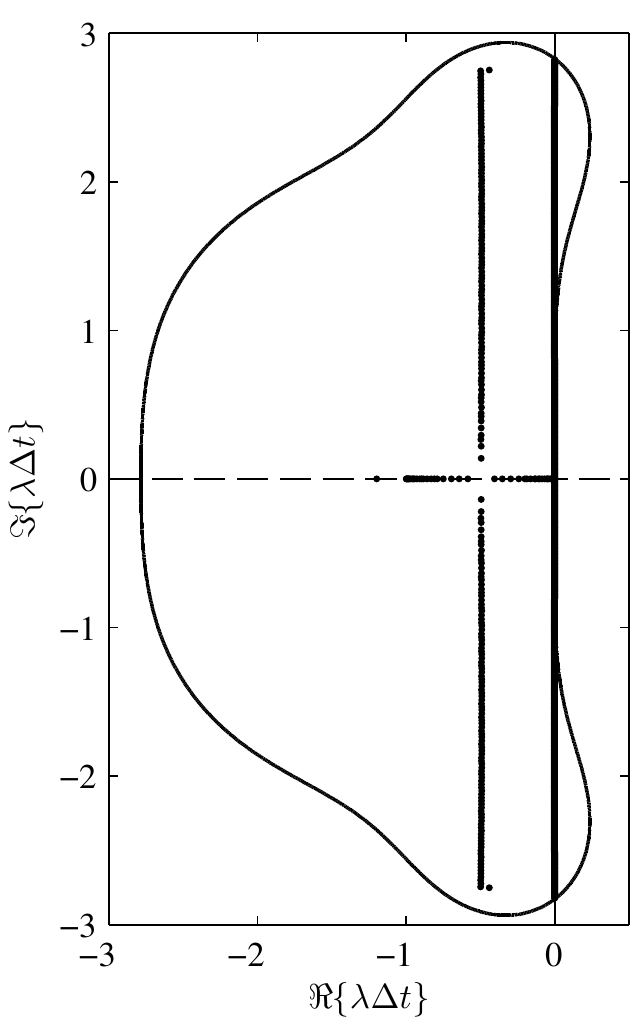} 
    \includegraphics[width = 0.45\textwidth]{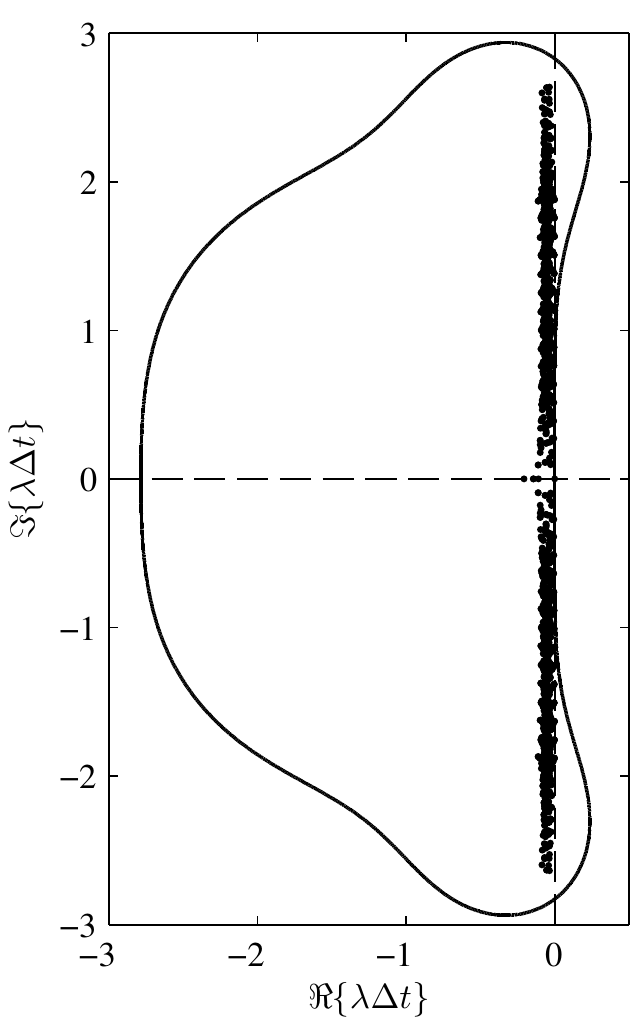} \\
    \caption{Left: Eigenvalues (dots) for the evolution operator $\mathbf{A}$ in one dimension ($m = 0$, $N = 1024$, with geometry described in Section \ref{sec_GaussianScattering1d}), along with the RK4 stability region (line). Right: Eigenvalues for the discrete evolution operator $\mathbf{A}$ for the two-dimensional $\textrm{TM}_z$ mode ($m = 0$, $Nx = Ny = 32$, with geometry described in Section \ref{sec_manufactured2d}). }
 \label{Eigenvalues}
\end{figure}

Once the numerical scheme is discretized according to the algorithm in Section \ref{sec_main_algorithm}, we examine stability by numerically computing the eigenvalues of the associated linear operators.  One should note that the penalized equations are the sum of two operators (the wave operator and the penalty operator) whose eigenvalues can be independently, analytically computed. Unfortunately the penalty term is non-normal, so that stability is not determined by the eigenvalues of the penalty term alone.  Alternatively, one can use energy arguments to show that for sufficiently small $L$ in the penalty term $\tilde{\mathbf{g}}$, one guarantees a strong stability preserving (SSP) scheme (in the $L^2(\Omega)$ norm), however for such $L$ one loses the global accuracy of the method. Therefore, to show stability for the current method, we compute the associated eigenvalues. To compute the eigenvalues, we write the numerical scheme in the form
\begin{eqnarray}
	\mathbf{u} \approx (\mathbf{E}_{\eta}, \mathbf{H}_{\eta})^T
\end{eqnarray}
and introduce the discretized operator
\begin{eqnarray}
	\mathbf{A} &=& \begin{pmatrix}
				-\eta^{-1} \chi(\mathbf{x}) (\mathbf{I} - \mathbf{G}) & \nabla \times \\
				-\nabla \times & 0
	\end{pmatrix}	
\end{eqnarray}
where $\mathbf{G}$ is the discrete operator that approximates the penalty term with zero boundary condition $\mathbf{\tilde{g}} \approx \mathbf{G} \mathbf{E}_{\eta}$.  The Maxwell's equations are then approximated by
\begin{align}
	\frac{\partial \mathbf{u}}{\partial t} = \mathbf{A} \mathbf{u}.
\end{align}
We compute the eigenvalues of $\mathbf{A}$ for numerous test cases and compare them to the RK4 stability region. Specifically, we compute the eigenvalues for dimension one matching $m = 0, 1, 2$ and for the two-dimensional $\textrm{TM}_z$ mode with $m = 0, 1$ for the domain with a hole removed.  Although we varied different values of the parameters $\eta$ and $L$, we show two typical eigenvalue plots in Figure \ref{Eigenvalues} indicating that the scheme is numerically stable.  Finally we remark that when matching higher derivatives in the numerical algorithm, one needs to add the extra filtering step outlined in Section \ref{sec_main_algorithm} part (c) for stability.  Mathematically, this filtering step modifies the matrix $\mathbf{G}$ of the penalty term to make the scheme stable without affecting accuracy.

\section{Numerical test cases}
\subsection{Test 1: One-dimensional Gaussian scattering} \label{sec_GaussianScattering1d}

In this section, we perform a numerical convergence study for the active penalty method.  We do so for the one-dimensional scattering of an incident Gaussian wave packet. Specifically, we seek solutions of the form $\mathbf{E}=\left(0, 0, E_{z}\left(x,t\right)\right)$ and $\mathbf{H}=\left(0,H_{y}\left(x,t\right),0\right)$ to Maxwell's
equations on the domain $x\in\left[0,\infty\right)$ such that
\begin{subequations} \label{Maxwell_1D}
\begin{align} \label{Maxwell_1D_H}
\frac{\partial H_{y}}{\partial t} & =\frac{\partial E_{z}}{\partial x},  \\ \label{Maxwell_1D_E}
\frac{\partial E_{z}}{\partial t} & =\frac{\partial H_{y}}{\partial x}.
\end{align}
\end{subequations}
In addition, we impose initial conditions
\begin{subequations} \label{InitialData_1D}
\begin{align} \label{InitialData_1D_H}
H_{y}\left(x,0\right) & =E_{0} \left\{ f(x - x_0) + f( -x - x_0 ) \right\}\chi_0(x) \\ \label{InitialData_1D_E}
E_{z}\left(x,0\right) & =E_{0} \left\{ f(x - x_0) - f( -x - x_0 ) \right\} \chi_0(x) \\
	f(\alpha) &= e^{-\frac{1}{2}\left(\frac{\alpha}{\sigma}\right)^{2}}\sin\left(\omega_{0} \alpha \right)
\end{align}
\end{subequations}
along with the boundary condition $E_{z}\left(0,t\right)=0$ and 
\begin{equation}
\chi_{0}\left(x\right)=\begin{cases}
0 & x < 0\\
1 & x \geq 0
\end{cases}.
\end{equation}
Provided $x_0/\sigma$ is large enough, the initial conditions simplify to a single incident Gaussian packet $E_{z}\left(x,0\right) \approx E_{0} f(x - x_0)$.  The Maxwell's equations then have the solution
\begin{subequations} \label{Exact_1D_Solution}
\begin{align} \label{Exact_1D_SolutionH}
H_{y}\left(x,t\right) & =E_{0}\left\{ f(t + x - x_0) + f( t - x - x_0) \right\}  \chi_0(x) \\ \label{Exact_1D_SolutionE}
E_{z}\left(x,t\right) & = E_{0}\left\{ f(t + x - x_0) - f( t - x - x_0 ) \right\} \chi_0(x).
\end{align}
\end{subequations}
 Here the first terms correspond to the incident wave, while the second terms correspond to the reflected wave. We compare the exact solution (\ref{Exact_1D_Solution})
to the numerical solution of the penalized equations
\begin{subequations} \label{Maxwell_1D_discretize}
\begin{align}
\frac{\partial H_{y,\eta}}{\partial t} & =\mathcal{F}^{-1}\left\{ \imath k_{x}\mathcal{F}\left\{ E_{z,\eta}\right\} \right\} \\
\frac{\partial E_{z,\eta}}{\partial t} & =\mathcal{F}^{-1}\left\{ \imath k_{x}\mathcal{F}\left\{ H_{y,\eta}\right\} \right\} -\eta^{-1}\chi_{h}\left(x\right)\left(E_{z,\eta}- \tilde{g}_{z}\right)
\end{align}
\end{subequations}
with initial data (\ref{InitialData_1D}), where
\begin{equation}
\chi_{h}\left(x\right)=\begin{cases}
0 & x > h\\
1 & x \leq h
\end{cases}.
\end{equation}
For our test, we take $E_0 = 1$, $x_0 = 7$, $\sigma = 1/\sqrt{2}$ and $\omega_0 = 10$ to be the parameters of the Gaussian wave packet. Meanwhile we take $\Omega_s = [-2, 0)$ and the physical domain $\Omega_0 = [0, 14)$ so that the box size is $D = 16$.  We then integrate equations (\ref{Maxwell_1D_discretize}) to a final time $T = 12$ using a 4th order Runge-Kutta (RK4) scheme with $\Delta t = 0.2 \Delta x$, and $\eta = h = \Delta x$. The extension function $\tilde{g}_{z}$ is constructed such that $\tilde{g}_{z}(0)=0$.

\begin{figure}[htb!]
	\centering
    \includegraphics[width = \textwidth]{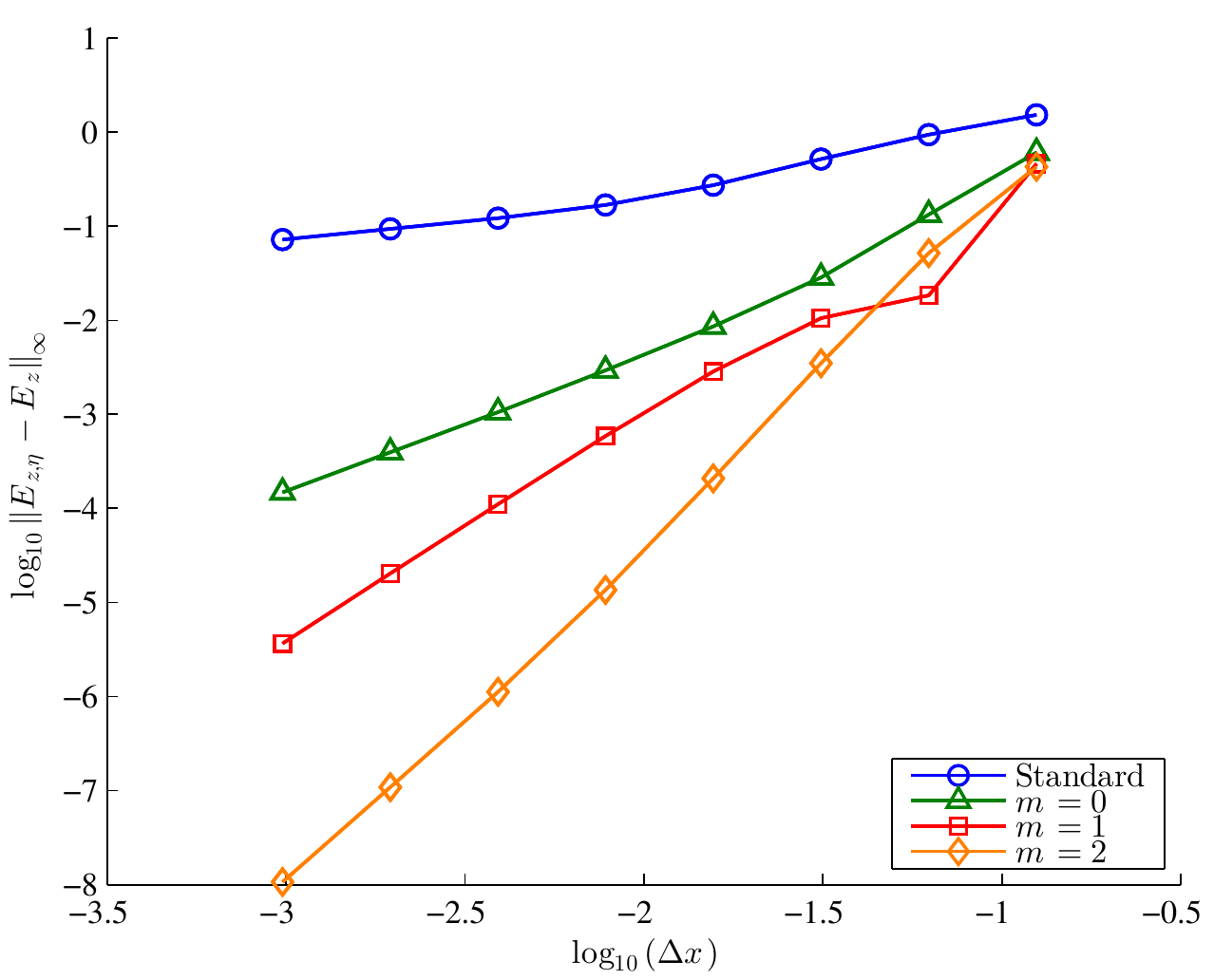} \\
    \caption{Convergence rates for a scattered Guassian (\ref{InitialData_1D}) with RK4 time-stepping.  The plot shows a non-active penalty method (blue-circles), and active methods matching $0$ (green-triangles), $1$ (red-squares), and $2$ (orange-diamonds) derivatives in constructing $\tilde{g}_z$. The global convergence rates in $L^{\infty}$ are approximately $0.38, 1.42, 2.48$ and $3.34$ respectively. } \label{GuassianScattering_CvgRates}
\end{figure}

\begin{figure}[htb!]
	\centering
    \includegraphics[width = \textwidth]{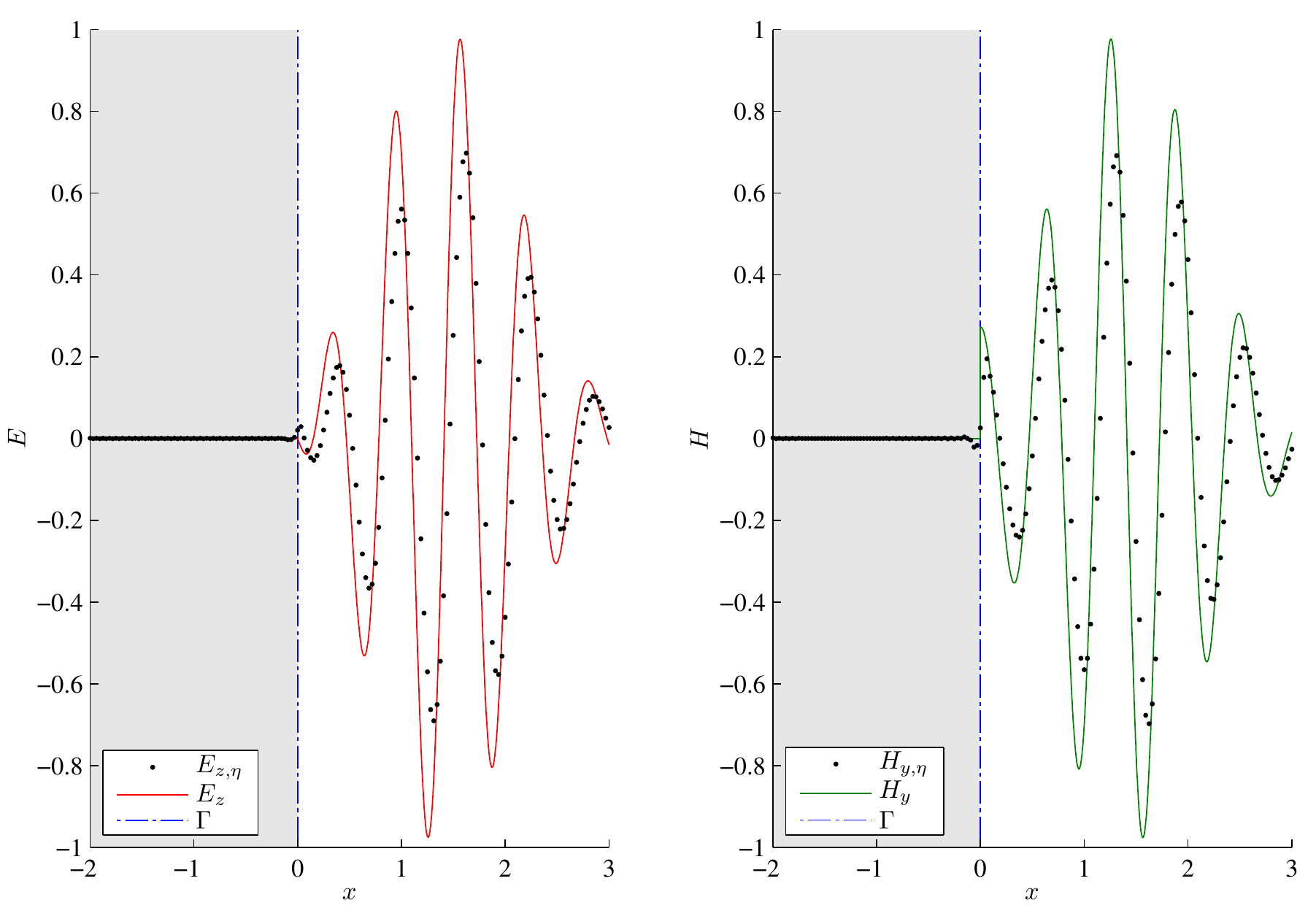} \\
    \caption{Comparison of the exact solution of a reflected Gaussian wave packet to that of a numerically computed approximation obtained using a standard non-active penalty method where $\tilde{g}_z = 0$.} \label{Maxwell_1D_standard}
\end{figure}

Figure \ref{GuassianScattering_CvgRates} compares the error for a scattered Gaussian using a non-active penalty method (where $\tilde{g}_z = 0$) to the proposed active penalty method matching $m = 0$, $m = 1$, and $m = 2$ derivatives at the interface $\Gamma = 0$. For each method, we compute the asymptotic convergence rate and report them to be $O(\Delta x^{\gamma})$ where $\gamma = 0.38, 1.42, 2.48$ and $3.34$ respectively.  We note that the rate of approximately $1.42$ when matching $m = 0$ derivatives is quite close to the predicted analytic rate of $1.5$ derived in Section \ref{AnalyticTM_Analysis}.

\begin{figure}[htb!]
	\centering
    \includegraphics[width = \textwidth]{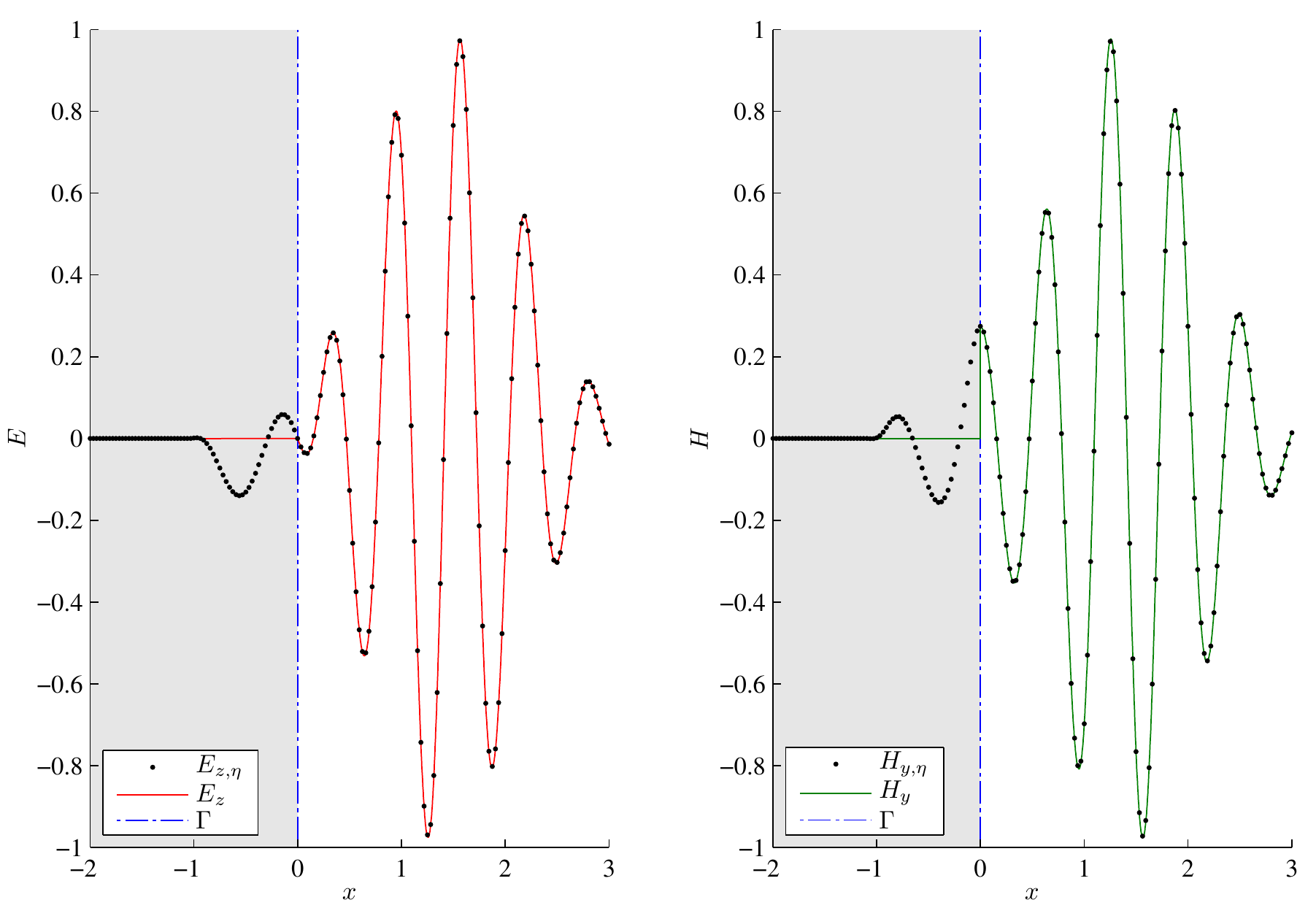} \\
    \caption{Comparison of the exact solution of a reflected Gaussian wave packet to that of a numerically computed approximation obtained using an active penalty method matching $m = 2$ derivatives. Note that as a result of the active penalty term $\tilde{g}_z$, the solution $E_{z,\eta}$ oscillates inside the PEC. }
\label{Maxwell_1D_k2}
\end{figure}

To illustrate the role of the active penalty term, we plot the penalized solution against the exact solution in the vicinity of the interface $\Gamma$.  Figure \ref{Maxwell_1D_standard} shows a standard non-active penalty method where $\tilde{g}_z = 0$ for all time.  Note that the poor convergence rate leads to a large error after the wave has reflected from the interface.  Meanwhile, Figure \ref{Maxwell_1D_k2} shows the penalized solution when matching $m = 2$ derivatives.  Here, the penalty term is a smooth extension which matches the boundary condition at $x = 0$. This results in a significant increase in accuracy for the same number of grid points.

\subsubsection{One-dimensional dispersion errors}

One difficulty which arises when using finite difference methods for solving a wave equation is the introduction of numerical dispersion errors.  Specifically, the numerical discrete dispersion relation can differ from the exact analytic one at large wavenumbers.  As a result, one must increase the resolution of the scheme, i.e., the number of grid points per wavelength (ppwl), with the wavenumber of the initial data. 

In contrast, provided one fixes the ppwl resolution, Fourier methods have been shown \cite{LyonBruno2010} to maintain a constant error over a wide range of wavelengths.  In this subsection we examine the pollution error for the proposed active penalty method.  Here we perform the same test as in the previous section using the initial data (\ref{Maxwell_1D}) ($E_0 = 1$, $x_0 = 7$, $\sigma = 1/\sqrt{2}$), however we vary $\omega_0 \in [10, 500]$.  In the test, we fix the ppwl at either $15$ or $20$ so that the total number of grid points increases with the frequency $\omega_0$ (or number of wavelengths) of the initial data.  As in the previous test cases, we take $\Delta t = C \Delta x$ where $C = 0.5$ for $m = 0$ and $m = 1$ derivatives and $C = 0.2$ for $m = 2$ derivatives.  We also set the integration time $T = 15$, $\eta = \Delta t$ and $h = 1.002 \Delta x$. Here the factor $1.002$ is taken slighly larger than $1$ to ensure that $h$ is at least one gridpoint away from $\Gamma$.  In all test cases the RK4 time stepping scheme is used to guarantee that time discretization error is smaller than the error associated with the introduction of the penalty term. 

\begin{figure}[htbp!]
	\centering
    \includegraphics[width = .79\textwidth]{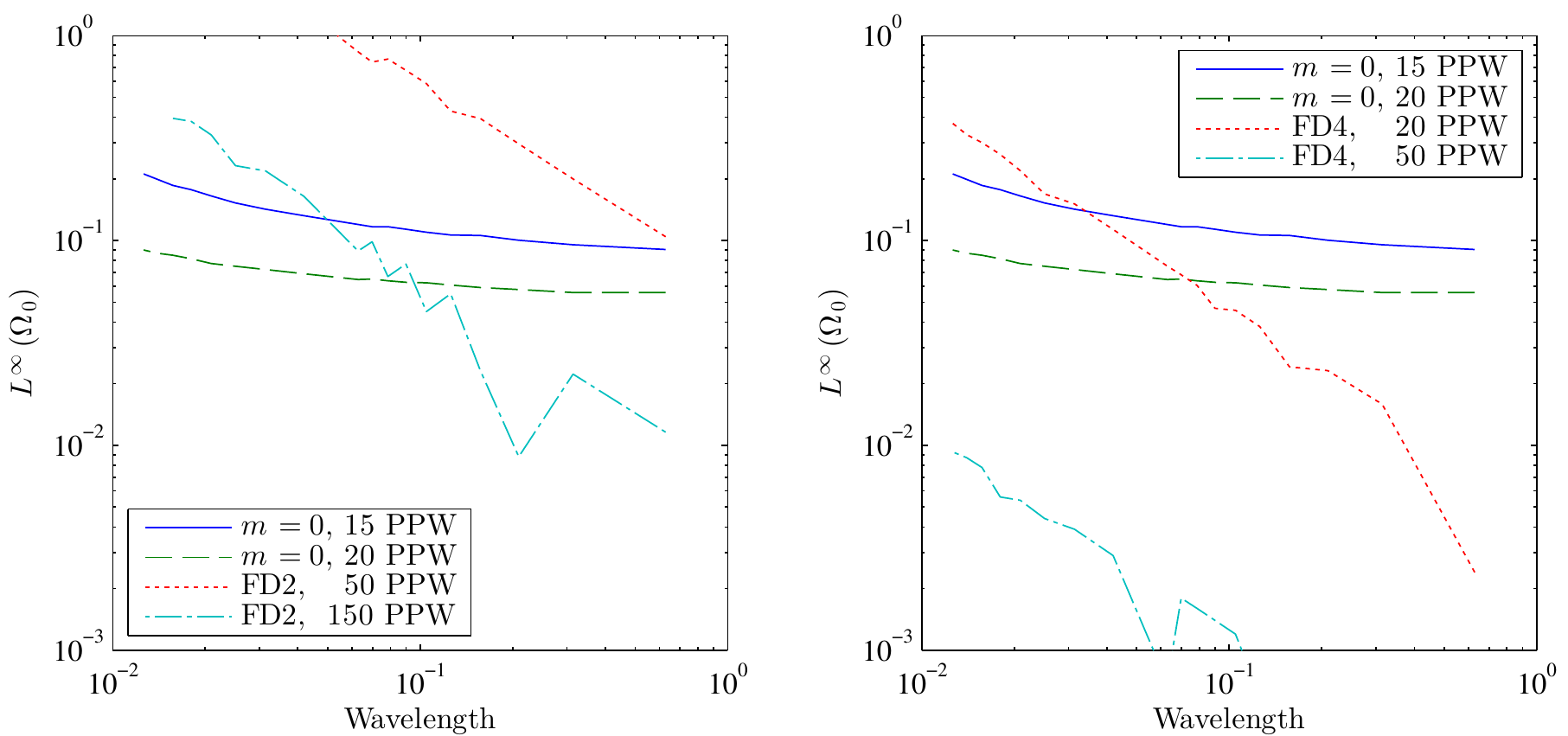} \\
    \includegraphics[width = .79\textwidth]{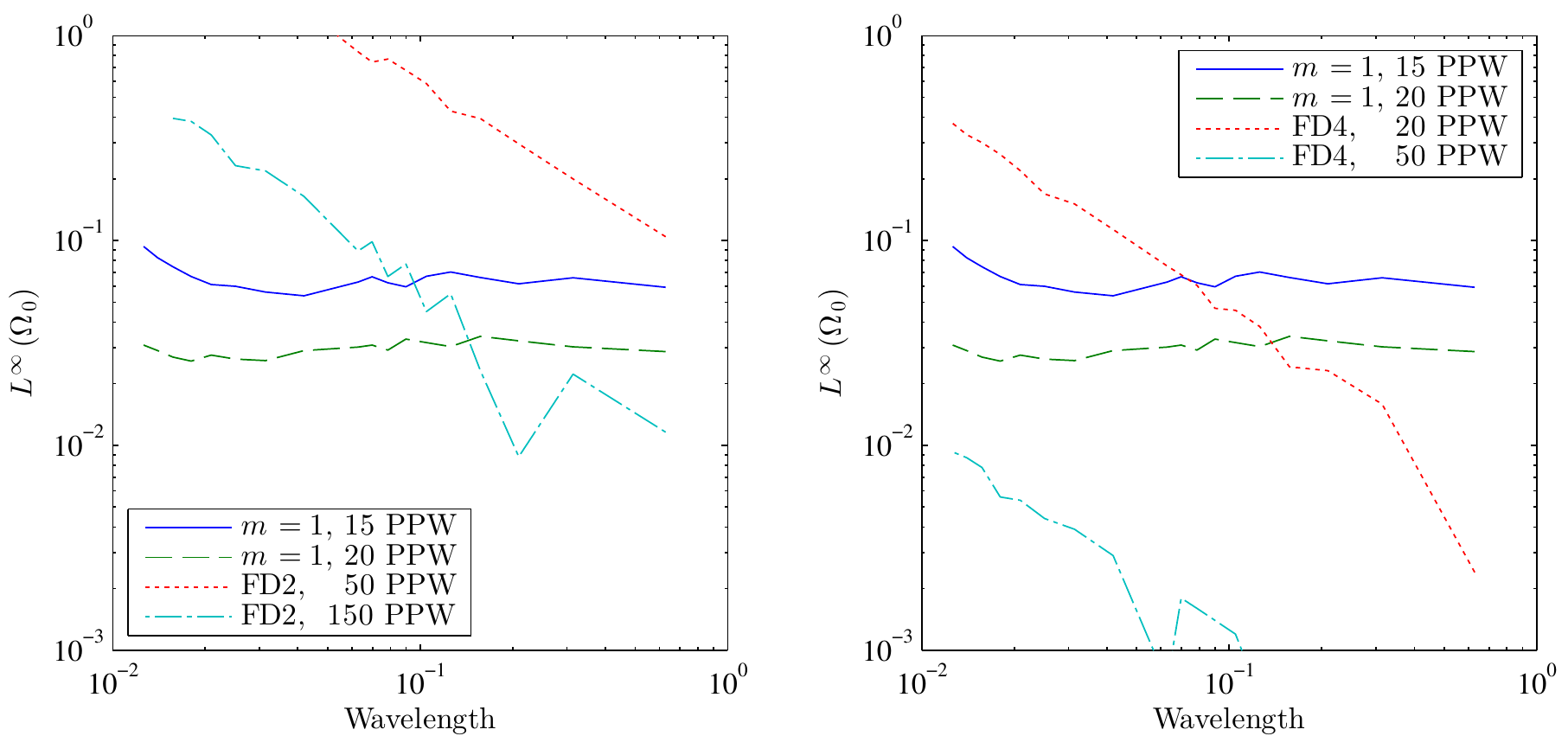} \\
    \includegraphics[width = .79\textwidth]{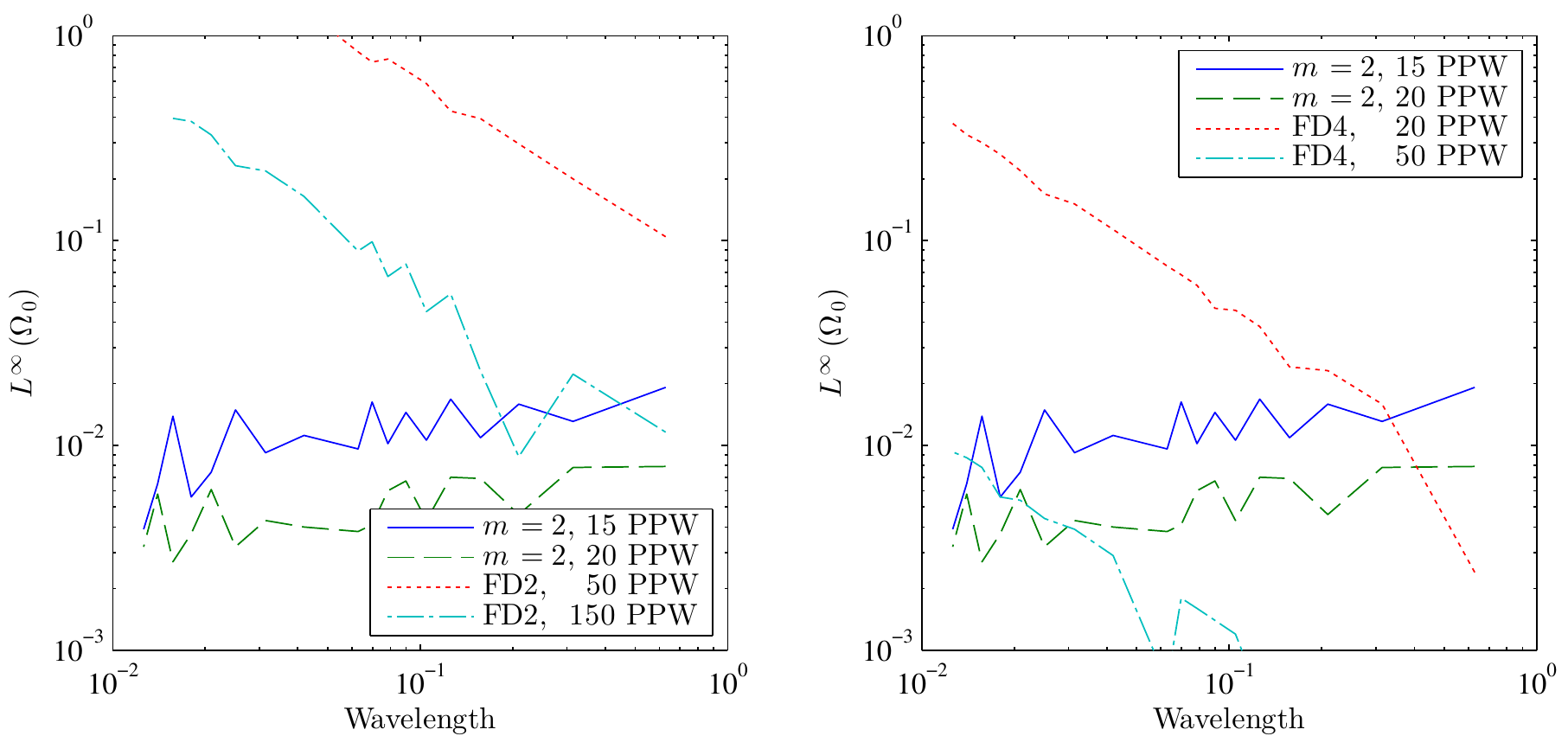} \\
    \caption{The pollution error (in $L^{\infty}(\Omega_0) )$ for various active penalty methods and finite difference schemes (FD2 and FD4 are 2nd and 4th order methods). For a fixed ppwl, the active penalty error remains flat over a wide range of wavenumbers. An active penalty method matching one derivative ($m = 1$) surpasses a 4th order finite difference scheme at moderate wave numbers $\sim 10^{-1}$.} \label{PollutionTest_1D}
\end{figure}

As shown in Figure \ref{PollutionTest_1D}, the error (in $L^{\infty}(\Omega_0)$) for active penalty methods remains relatively flat over a wide range of wavelengths.  The plots also show 2nd and 4th order finite difference schemes.  As expected, both finite difference schemes show an increase in error as the wavenumber increases.  

\subsection{Test 2: Two-dimensional manufactured solutions for a domain with a circular hole} \label{sec_manufactured2d}

In the following section, we test the Fourier penalty method (Section \ref{sec_main_algorithm}) using a manufactured solution approach on a periodic domain $\Omega=\left[0,2\pi\right]^2$, with a circular hole removed.  Specifically, the boundary $\Gamma$ of the hole is given 
 by the zero level set of the signed distance function
\begin{equation}
\psi\left(x,y\right)=\sqrt{\left(x-x_{0}\right)^{2}+\left(y-y_{0}\right)^{2}}-a. \label{eq:levelset}
\end{equation}
The zero level set of (\ref{eq:levelset}) is circular with radius $a$ and center $(x_0,y_0)$. In our tests, we fix $a=2$ with center $x_{0}=y_{0}=\pi$.

Here the manufactured solution approach allows for the direct convergence test of the penalized Maxwell's equations.  Two test problems are chosen to verify two independent modes of propagation
supported by the two-dimensional Maxwell's equations. We treat
the $\mathrm{TM}_{z}$ mode followed by the $\mathrm{TE}_{z}$ mode. In each case, we impose a boundary condition $\mathbf{n}\times ( \mathbf{E} - \mathbf{g}) = 0$. Recall that when $\mathbf{g} = \mathbf{0}$, we have the boundary condition for a PEC.

\subsubsection{$TM_{z}$ mode} \label{section_TM_manufacture}

We seek solutions to the forced Maxwell's equations of the form $\mathbf{E}=\left(0,0,E_{z}\left(x,y,t\right)\right)$
and $\mathbf{H}=\left(H_{x}\left(x,y,t\right),H_{y}\left(x,y,t\right),0\right)$
which satisfy
\begin{subequations} \label{Maxwell_2D_TM}
\begin{align} \label{Maxwell_2D_TM_begin}
\frac{\partial H_{x}}{\partial t} & =-\frac{\partial E_{z}}{\partial y}\\
\frac{\partial H_{y}}{\partial t} & = \phantom{-} \frac{\partial E_{z}}{\partial x}\\ \label{Maxwell_2D_TM_end}
\frac{\partial E_{z}}{\partial t} & = \phantom{-} \frac{\partial H_{y}}{\partial x}-\frac{\partial H_{x}}{\partial y}+F
\end{align}
\end{subequations}
with forcing function 
\begin{equation} \label{forcing_function_TM}
F= \sin\left(x\right)\cos\left(y\right)\sin\left(t\right),
\end{equation}
and initial conditions
\begin{subequations} \label{initial_conditions_TM}
\begin{align} \label{initial_conditions_TM_begin}
H_{x}\left(x,y,0\right) & =0\\
H_{y}\left(x,y,0\right) & =0\\ \label{initial_conditions_TM_end}
E_{z}\left(x,y,0\right) & =\sin\left(x\right)\cos\left(y\right).
\end{align}
\end{subequations}
The solution to (\ref{Maxwell_2D_TM}) with forcing function (\ref{forcing_function_TM}), subject to initial conditions (\ref{initial_conditions_TM}) is given by
\begin{subequations}
\begin{align}
H_{x}\left(x,y,t\right) & =\sin\left(x\right)\sin\left(y\right)\sin\left(t\right)\\
H_{y}\left(x,y,t\right) & =\cos\left(x\right)\cos\left(y\right)\sin\left(t\right)\\ \label{maxwell_2D_TM_exact}
E_{z}\left(x,y,t\right) & =\sin\left(x\right)\cos\left(y\right)\cos\left(t\right),
\end{align}
\end{subequations}
which one can verify satisfies the divergence-free criteria.

\begin{figure}[htb!]
	\centering
    \includegraphics[width = \textwidth]{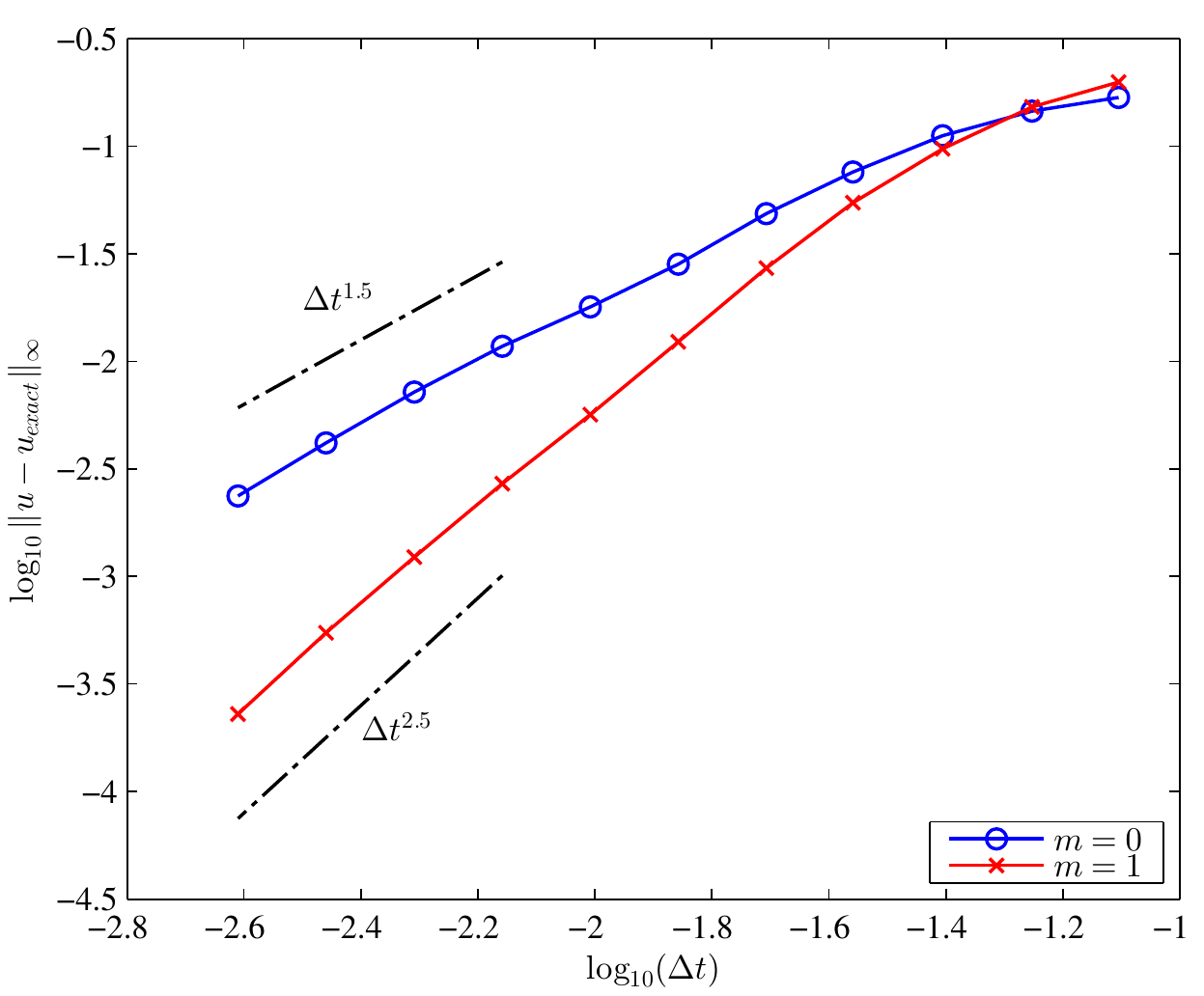} \\
    \caption{Convergence study for the manufactured $\textrm{TM}_z$ example. The global convergence rates in $L^{\infty}(\Omega_0)$ are approximately 1.5 and 2.5 for $m=0$ and $m=1$ respectively. Here, $u=(H_{x,\eta},H_{y,\eta},E_{z,\eta})$, and a smoothing parameter $c_f = 16$ was used for the $m=1$ case. No smoothing was required for $m=0$.} \label{TM_manufacture_convergence_rates}
\end{figure}

We then compare the exact solution to the solution of the penalized equations
\begin{subequations} \label{maxwell_2D_TM_disc}
\begin{align} \label{maxwell_2D_TM_disc_begin}
\frac{\partial H_{x,\eta}}{\partial t} & =- \mathcal{F}^{-1}\left\{ \imath k_{y}\mathcal{F}\left\{ E_{z,\eta}\right\} \right\} \\
\frac{\partial H_{y,\eta}}{\partial t} & = \phantom{-} \mathcal{F}^{-1}\left\{ \imath k_{x}\mathcal{F}\left\{ E_{z,\eta}\right\} \right\} \\ \label{maxwell_2D_TM_disc_end}
\frac{\partial E_{z,\eta}}{\partial t} & = \phantom{-} \mathcal{F}^{-1}\left\{ \imath k_{x}\mathcal{F}\left\{ H_{y,\eta}\right\} \right\} - \mathcal{F}^{-1}\left\{ \imath k_{y}\mathcal{F}\left\{ H_{x,\eta}\right\} \right\} -\eta^{-1}\chi_{h}\left(x,y\right)\left(E_{z,\eta}-\tilde{g}_{z}\right)+F
\end{align}
\end{subequations}
at grid points belonging to the physical domain $\Omega_0$.
We note that in this instance, $\tilde{g}_{z}$ is constructed to handle the non-zero boundary condition in the exact solution (\ref{maxwell_2D_TM_exact}).
We integrate (\ref{maxwell_2D_TM_disc}) to a final time $T=1.1\pi$ using RK4 with $\Delta t =0.4 \Delta x$, $h=2 \Delta x$, $\eta = 4 \Delta t$, and $L=1$. Figure \ref{TM_manufacture_convergence_rates} illustrates the convergence rates for the $m=0$ and $m=1$ cases. The results agree with the expected rates of 1.5 and 2.5.

\subsubsection{$TE_{z}$ mode}

\begin{figure}[htb!]
	\centering
    \includegraphics[width = \textwidth]{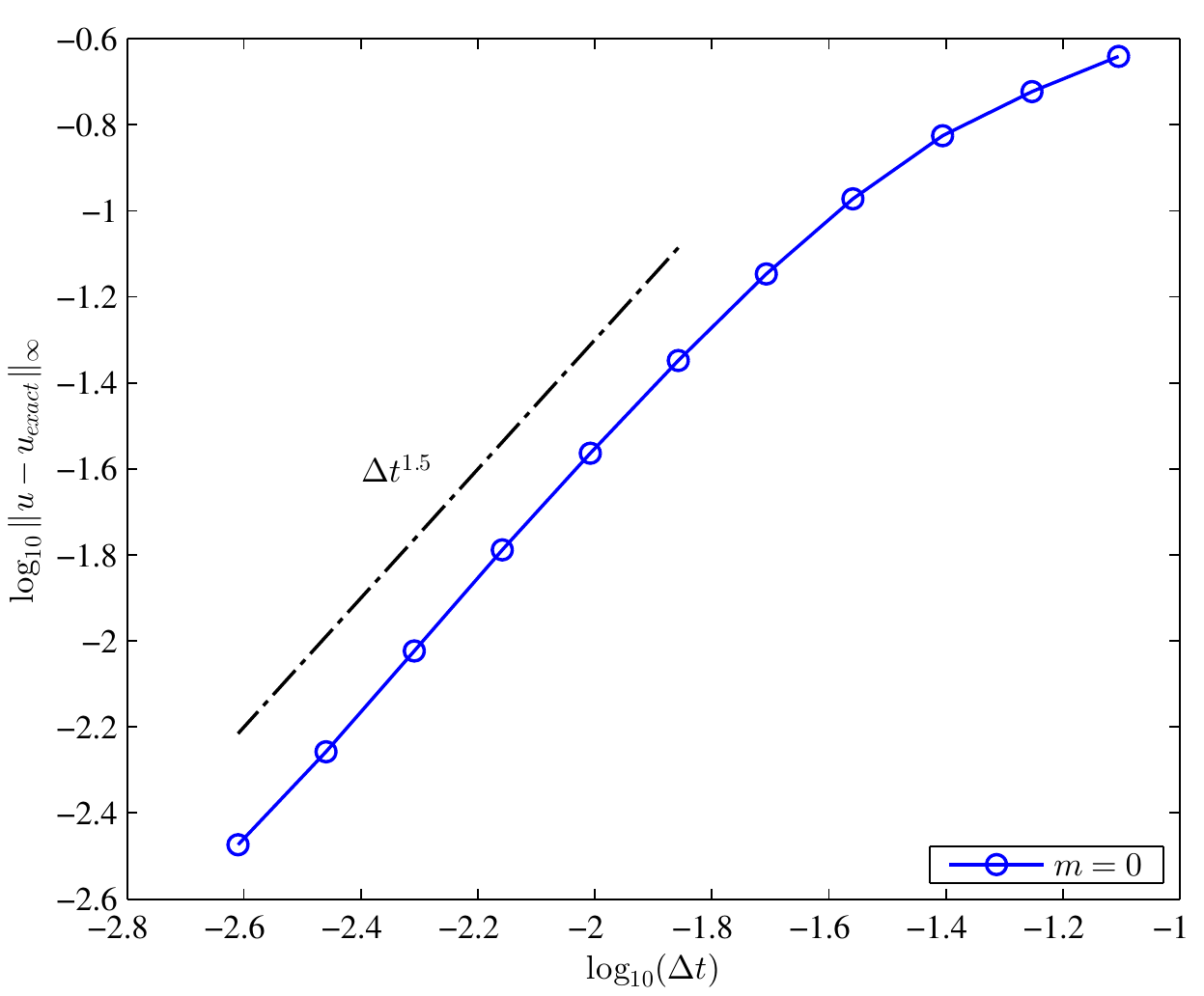} \\
    \caption{Convergence study for the manufactured $\textrm{TE}_z$ example. The global convergence rate in $L^{\infty}(\Omega_0)$ is approximately 1.5 for $m=0$. Here, $u=(E_{x,\eta},E_{y,\eta},H_{z,\eta})$.} \label{TE_manufacture_convergence_rates}
\end{figure}

Similarly, we may also seek solutions of the form $\mathbf{E}=\left(E_{x}\left(x,y,t\right),E_{y}\left(x,y,t\right),0\right)$
and $\mathbf{H}=\left(0,0,H_{z}\left(x,y,t\right)\right)$ over $\Omega = [0,2\pi]^2$
such that
\begin{subequations} \label{Maxwell_2D_TE}
\begin{align} \label{Maxwell_2D_TE_begin}
\frac{\partial E_{x}}{\partial t} & =\phantom{-} \frac{\partial H_{z}}{\partial y}\\
\frac{\partial E_{y}}{\partial t} & =-\frac{\partial H_{z}}{\partial x}\\  \label{Maxwell_2D_TE_end}
\frac{\partial H_{z}}{\partial t} & =\phantom{-} \frac{\partial E_{x}}{\partial y}-\frac{\partial E_{y}}{\partial x}+F
\end{align}
\end{subequations}
with forcing function 
\begin{equation} \label{forcing_function_TE}
F=\sin\left(x\right)\cos\left(y\right)\sin\left(t\right),
\end{equation}
and initial conditions
\begin{subequations} \label{initial_conditions_TE}
\begin{align} \label{initial_conditions_TE_begin}
E_{x}\left(x,y,0\right) & =0\\
E_{y}\left(x,y,0\right) & =0\\ \label{initial_conditions_TE_end}
H_{z}\left(x,y,0\right) & =\sin\left(x\right)\cos\left(y\right).
\end{align}
\end{subequations}
The solution to (\ref{Maxwell_2D_TE}) with forcing function (\ref{forcing_function_TE}), subject to initial conditions (\ref{initial_conditions_TE}) is given by
\begin{subequations}
\begin{align}
E_{x}\left(x,y,t\right) & =-\sin\left(x\right)\sin\left(y\right)\sin\left(t\right)\\
E_{y}\left(x,y,t\right) & =-\cos\left(x\right)\cos\left(y\right)\sin\left(t\right)\\
H_{z}\left(x,y,t\right) & = \phantom{-} \sin\left(x\right)\cos\left(y\right)\cos\left(t\right)
\end{align}
\end{subequations}
which one can verify satisfies the divergence-free criteria. Using the same circular obstacle as in Section \ref{section_TM_manufacture},
the penalized solution is computed by integration in time (using RK4) of 
\begin{subequations}
\begin{align}
\frac{\partial E_{x,\eta}}{\partial t} & =\phantom{-} \mathcal{F}^{-1}\left\{ \imath k_{y}\mathcal{F}\left\{ H_{z,\eta}\right\} \right\} -\eta^{-1}\chi_{h}\left(x,y\right)\left(E_{x,\eta}-\tilde{g}_{x}\right)\\
\frac{\partial E_{y,\eta}}{\partial t} & =-\mathcal{F}^{-1}\left\{ \imath k_{x}\mathcal{F}\left\{ H_{z,\eta}\right\} \right\} -\eta^{-1}\chi_{h}\left(x,y\right)\left(E_{y,\eta}-\tilde{g}_{y}\right)\\
\frac{\partial H_{z,\eta}}{\partial t} & = \phantom{-} \mathcal{F}^{-1}\left\{ \imath k_{y}\mathcal{F}\left\{ E_{x,\eta}\right\} \right\} -\mathcal{F}^{-1}\left\{ \imath k_{x}\mathcal{F}\left\{ E_{y,\eta}\right\} \right\} +F,
\end{align}
\end{subequations}
where $\tilde{\mathbf{g}}$ is constructed using $\mathbf{g} = \mathbf{E}$. That is, we penalize the electric field such that the tangential component
at the boundary $\Gamma$ is equal to the tangential component of
the exact solution. Figure \ref{TE_manufacture_convergence_rates} illustrates the convergence rate for the $m=0$ case. The parameter values for $T$, $\Delta t$, $h$, $\eta$, and $L$  are unchanged from Section \ref{section_TM_manufacture}.

\subsection{Test 3: Solution inside a circular cavity}
We may also examine a problem similar to test case 2 where we solve the penalized equations on the interior of a circular cavity domain that is embedded in a periodic domain.  In this case, we use the level set
\begin{equation}
\psi\left(x,y\right)=1-\sqrt{\left(x-x_{0}\right)^{2}+\left(y-y_{0}\right)^{2}}
\end{equation}
to construct the local coordinate system for $\tilde{g}_{z}$.  The physical, circular domain then corresponds to the region where $\psi > 0$.

\begin{figure}[tb!]
    \centering
    \includegraphics[width = .325\textwidth]{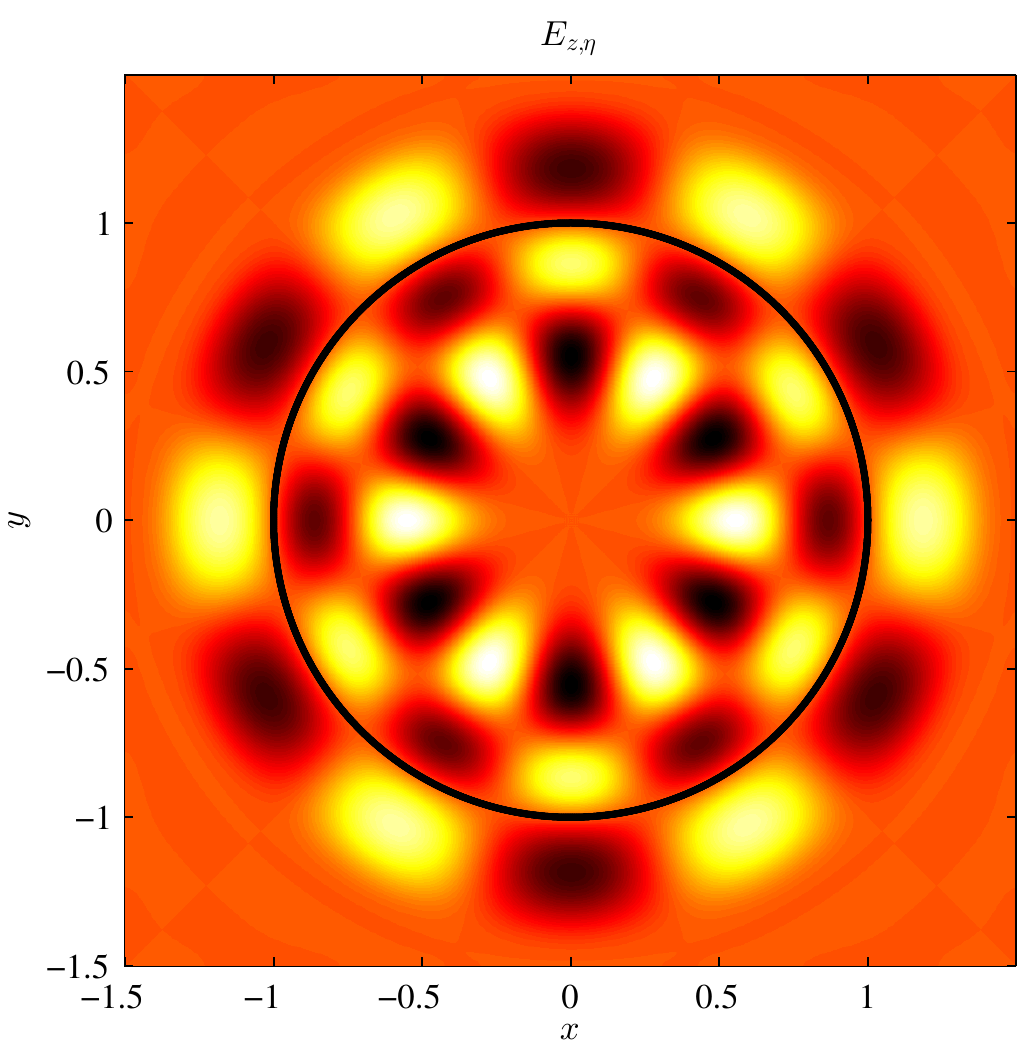} 
    \includegraphics[width = .325\textwidth]{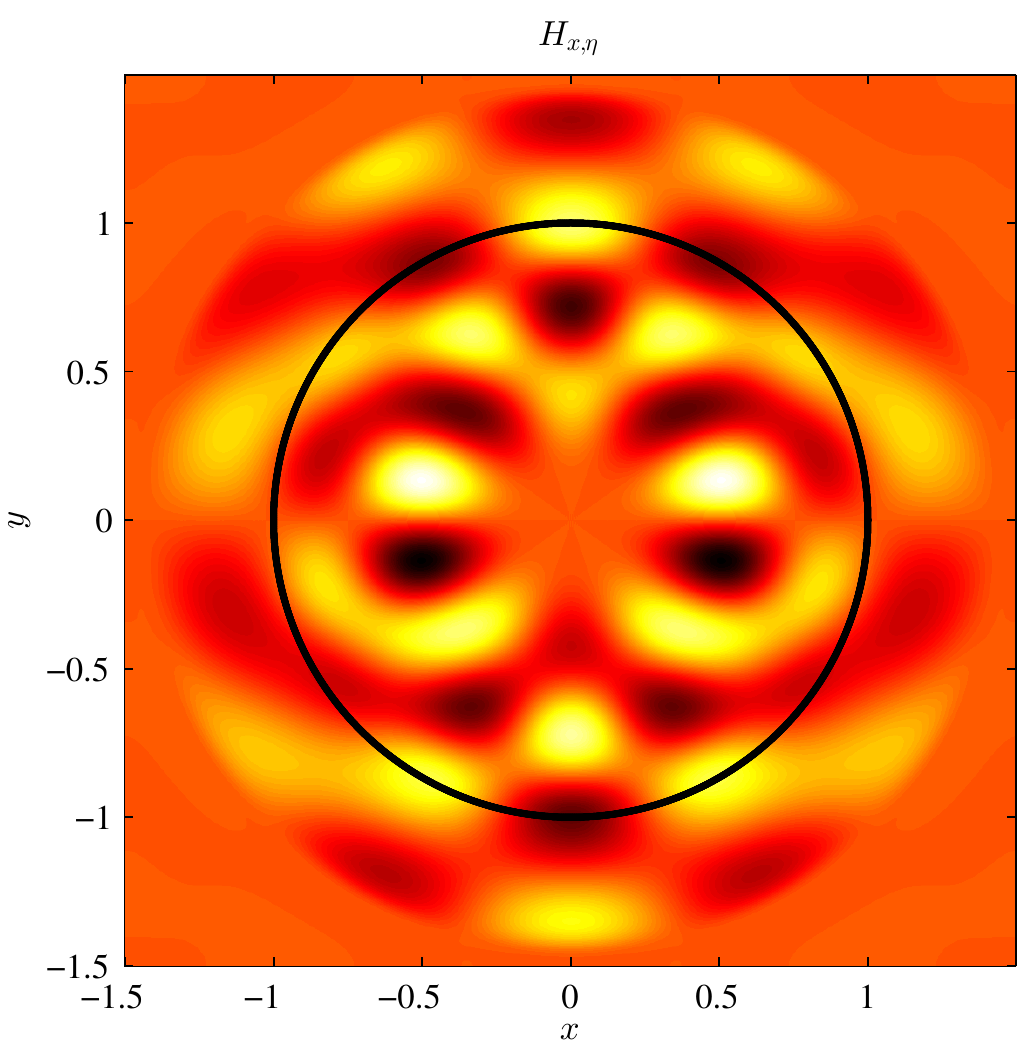} 
    \includegraphics[width = .325\textwidth]{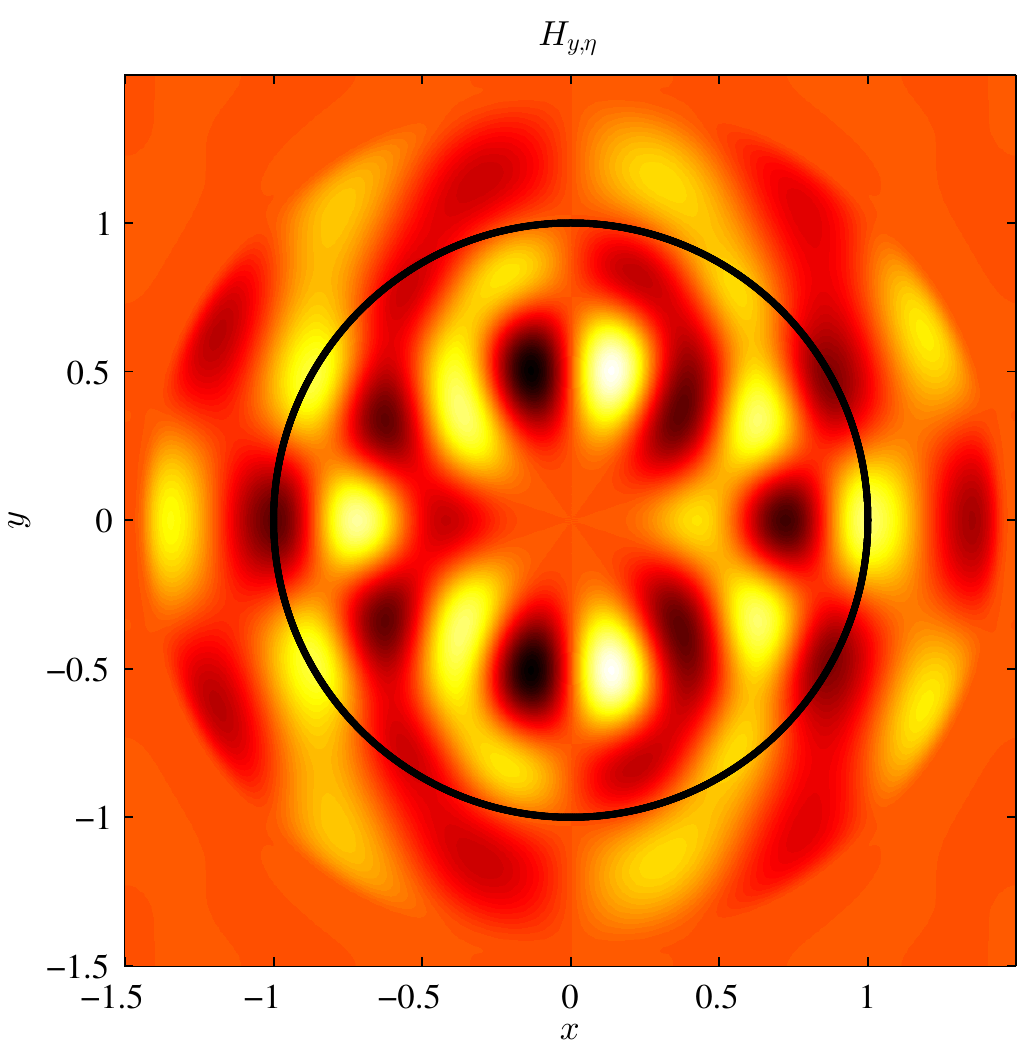} 
    \caption{Plots of the $\textrm{TM}_z$ mode cavity problem solution using the penalization method. Note how the $E_{z,\eta}$ component tracks the penalty function $\tilde{g}_z$ outside of $\Omega_0$.} \label{Fig_TM_eigenfunction}
\end{figure}

Now consider solving (\ref{TM_z_Mode}) with initial conditions
\begin{align}
H_{x,\eta}\left(\rho,\phi,0\right) & =0\\
H_{y,\eta}\left(\rho,\phi,0\right) & =0\\
E_{z,\eta}\left(\rho,\phi,0\right) & =J_{i}\left(\alpha_{i,j}\rho\right)\cos\left(i\phi\right),
\end{align}
where $J_{i}$ is the Bessel function of the first kind of integer
order $i$ and $\alpha_{i,j}$ is the $j^{\textrm{th}}$ positive
real root of the order $i$ Bessel function. The solution to the unpenalized
$\textrm{TM}_{z}$ equations is 
\begin{subequations} \label{TM_eigenfunction_exact_solution}
\begin{align}
H_{\rho}\left(\rho,\phi,t\right) & =\frac{i}{\alpha_{i,j}\rho}J_{i}\left(\alpha_{i,j}\rho\right)\sin\left(i\phi\right)\sin\left(\alpha_{i,j}t\right)\\
H_{\phi}\left(\rho,\phi,t\right) & =\frac{1}{2}\left[J_{i-1}\left(\alpha_{i,j}\rho\right)-J_{i+1}\left(\alpha_{i,j}\rho\right)\right]\cos\left(i\phi\right)\sin\left(\alpha_{i,j}t\right)\\
E_{z}\left(\rho,\phi,t\right) & =J_{i}\left(\alpha_{i,j}\rho\right)\cos\left(i\phi\right)\cos\left(\alpha_{i,j}t\right).
\end{align}
\end{subequations}
Figure \ref{Fig_TM_eigenfunction} illustrates the $E_{z}$, $H_{x}$, and $H_{y}$
components of the solution computed using the proposed penalization
method at time $T=0.3$ with $i=6$, $j=2$, and $\alpha_{6,2}\approx13.5892$.
For this example, $\Delta t=0.4\Delta x$, $\eta=4\Delta t$, $h=2\Delta x$,
$m=1$, and $L=0.45$. 

\begin{figure}[htb!]
\centering
\includegraphics[width = \textwidth]{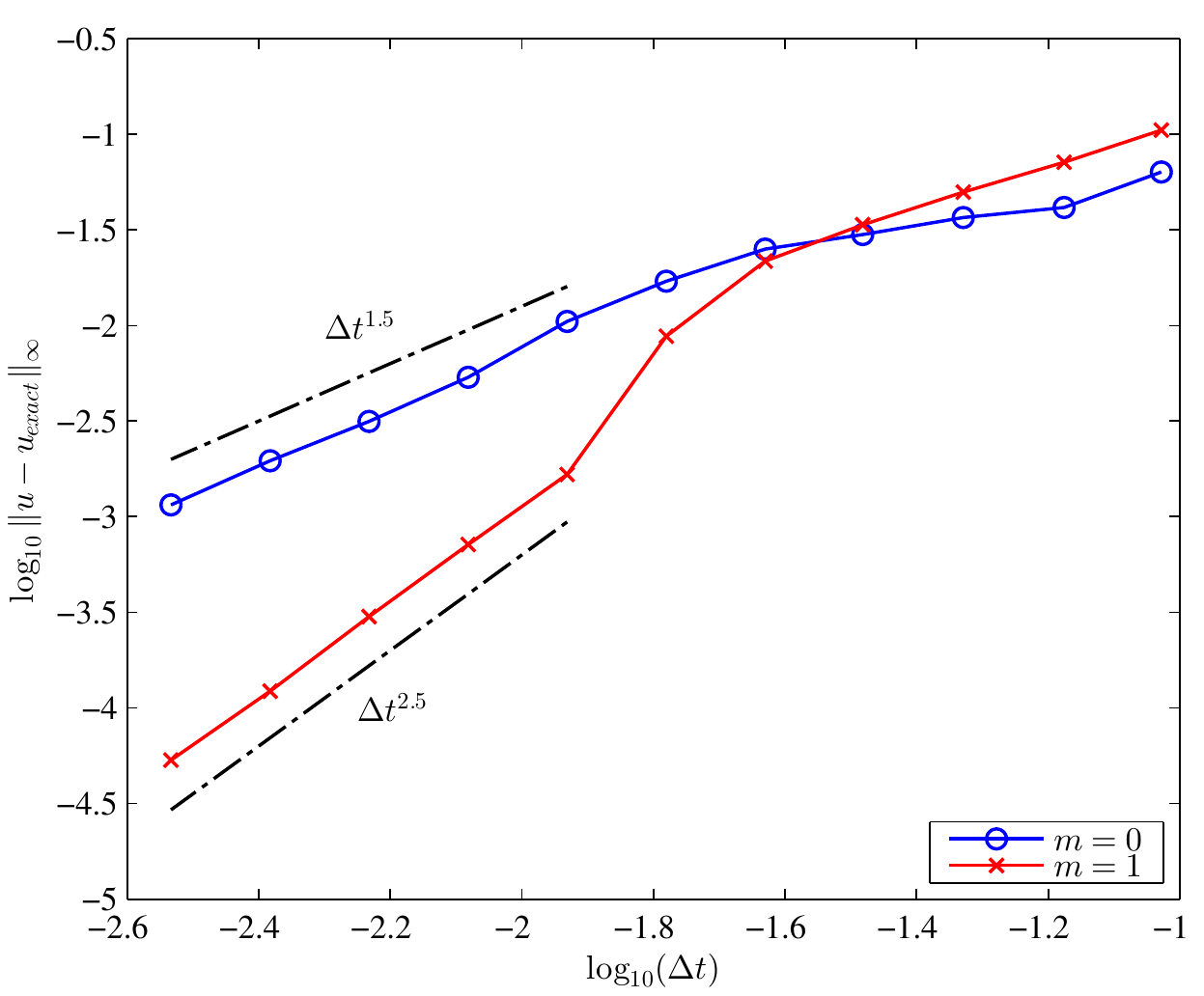}
\caption{Convergence study for the $\textrm{TM}_{z}$ mode cylindrical PEC cavity problem. The global convergence rates in $L^{\infty}(\Omega_0)$ are approximately 1.5 and 2.5 for $m=0$ and $m=1$ respectively. Here, $u=\left(H_{x,\eta},H_{y,\eta},E_{z,\eta}\right)$, and a smoothing parameter $c_f = 16$ was used for the $m=1$ case. No smoothing was required for $m=0$.}
\label{Fig_eigenfunction_convergence}
\end{figure}

Figure \ref{Fig_eigenfunction_convergence} shows the convergence of the penalized
solution to the exact solution (\ref{TM_eigenfunction_exact_solution}).

\subsection{Test 4: Scattering off of a PEC cylinder} \label{section_cylinder_scattering}

Our last test involving the circular geometry is for a time-dependent scattering computation.  Specifically, we take the following modulated Gaussian wave packet as initial data for a $\textrm{TE}_z$ mode 
\begin{subequations}
\begin{align}
	E_x(x, y, 0) &= 0 \\
	E_y(x, y, 0) &= \frac{2}{\sigma^2} (x - x_0) e^{-(\frac{x - x_0}{\sigma})^2} \\ 
	H_z(x, y, 0) &= \frac{2}{\sigma^2} (x - x_0) e^{-(\frac{x - x_0}{\sigma})^2}, 
\end{align} 
\end{subequations}
and compute the scattered wave packet off of a cylinder using our Fourier penalty method.  Here we take the domain parameters to be $\Omega = [-1, 1.5]^2$ with the cylinder centered at $(x_0,y_0)=(0, 0)$ with radius $a = 0.2$.  The initial data is chosen to have $\sigma = 0.125$ and $x_0 = -0.6$.  

\begin{figure}[htb!]
\centering
\includegraphics[width = \textwidth]{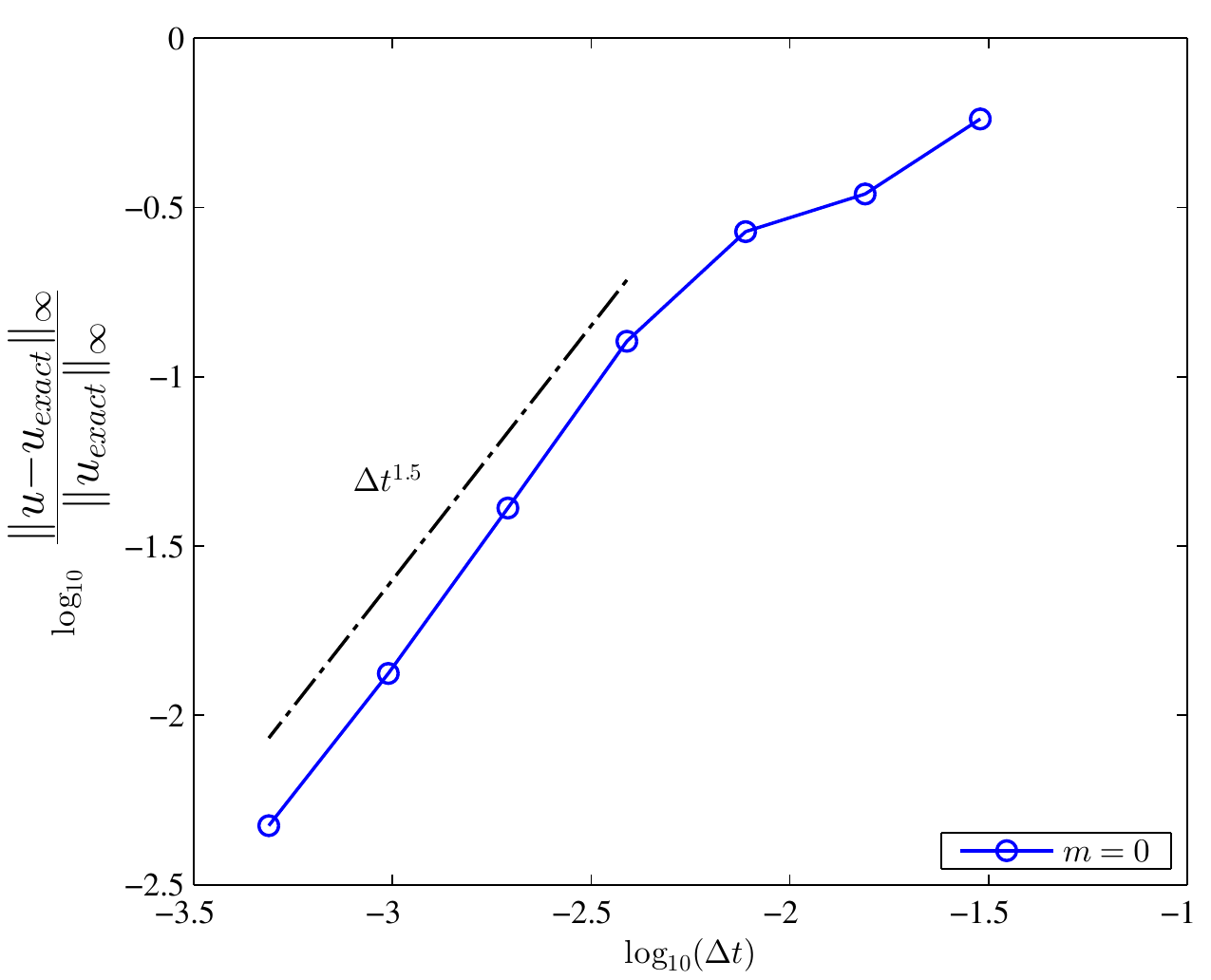}
\caption{Convergence study for the $\textrm{TE}_{z}$ mode cylindrical PEC
scattering problem. The global convergence rate in the relative error
is approximately 1.5 for $m=0$. Here, $u=\left(E_{x,\eta},E_{y,\eta},H_{z,\eta}\right)$.}
\label{Fig_scatter_convergence}
\end{figure}

\begin{figure}[htb!]
    \centering
    \includegraphics[width = .332\textwidth]{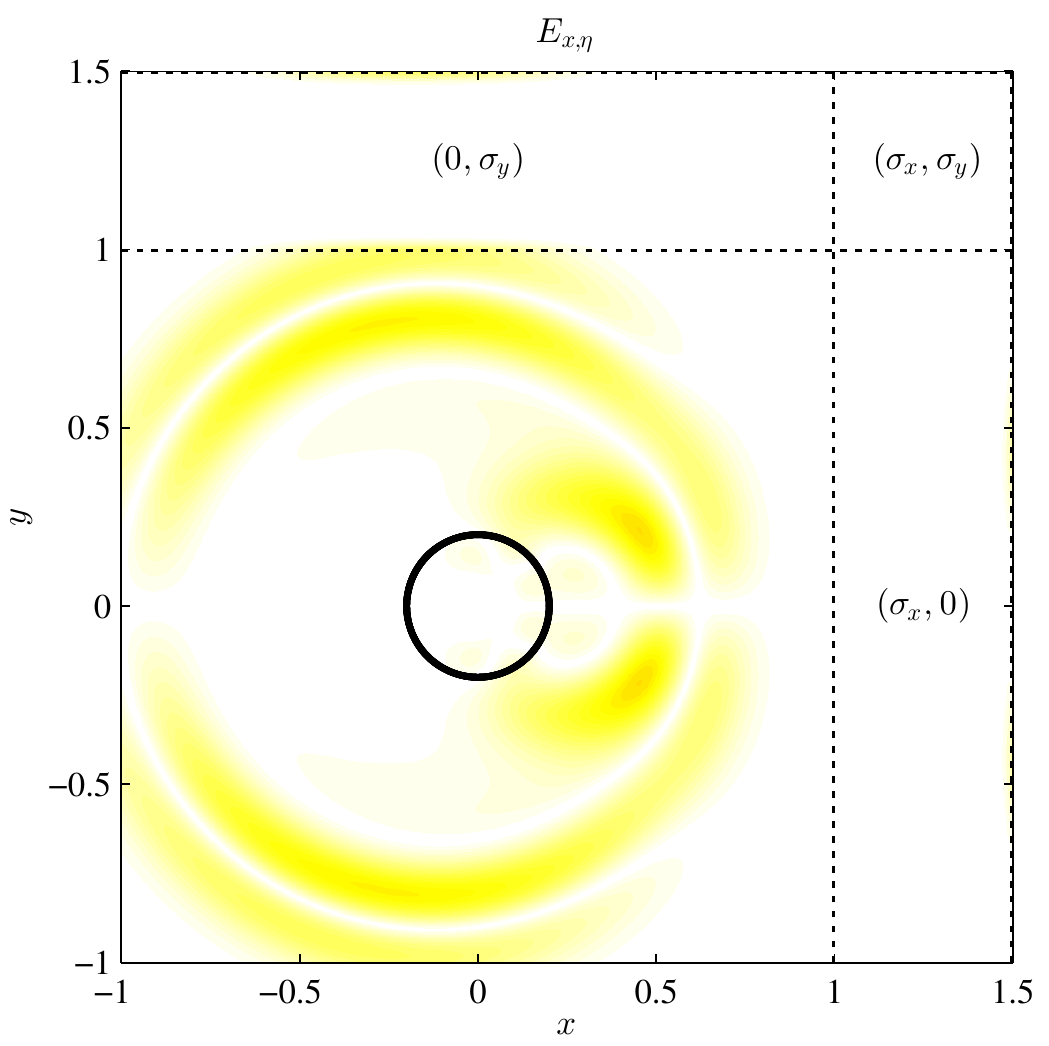} 
    \includegraphics[width = .325\textwidth]{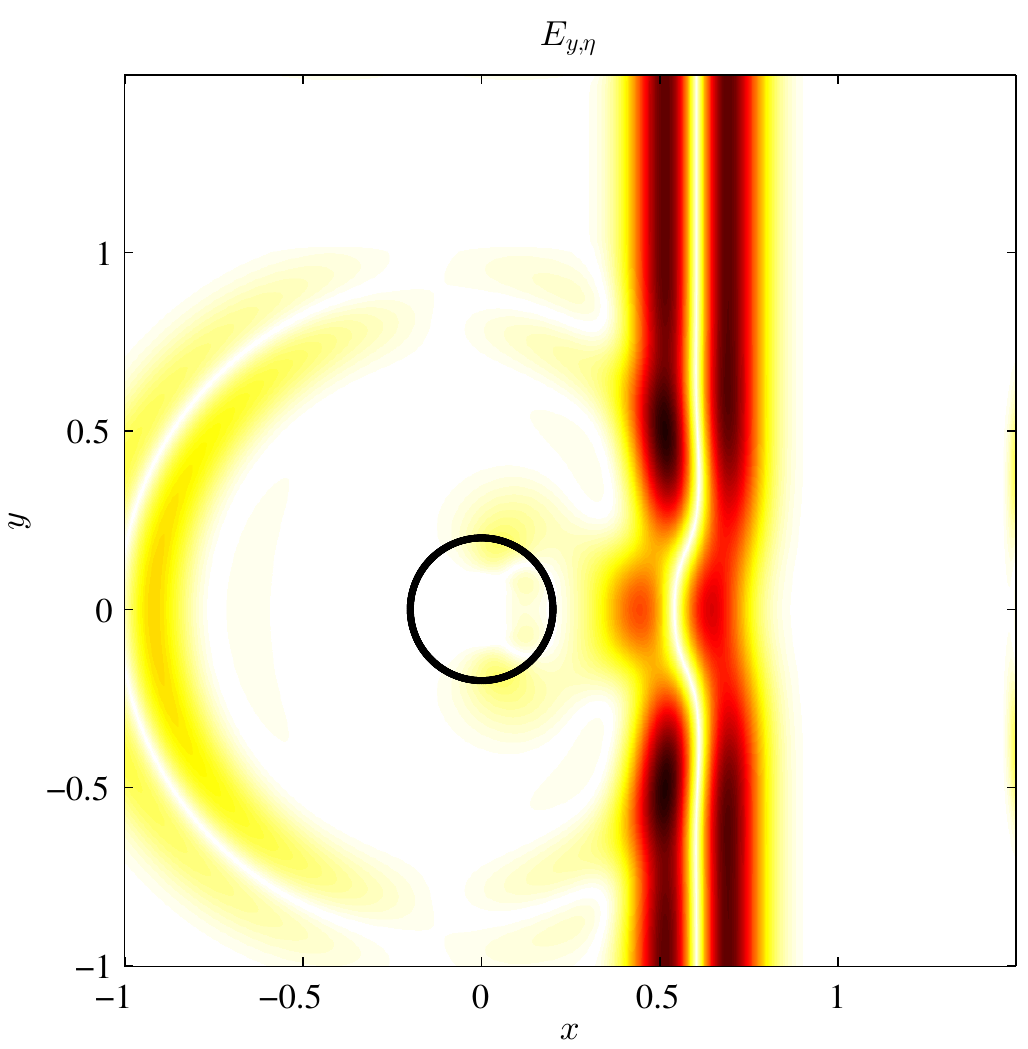}
    \includegraphics[width = .325\textwidth]{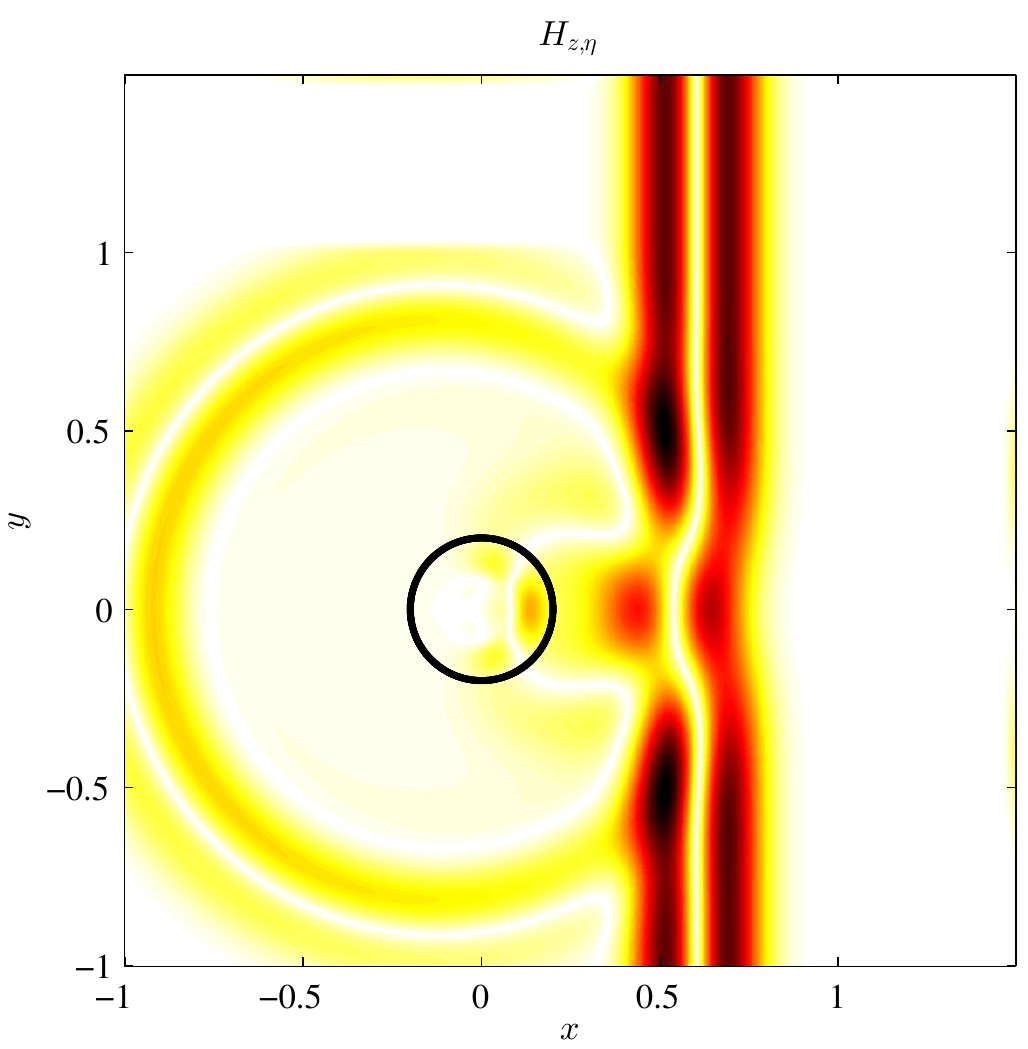} 
    \caption{Plots of the magnitudes of the (l-r) $E_{x,\eta}$, $E_{y,\eta}$ and $H_{z,\eta}$ components for a $\textrm{TE}_z$ wave scattering off of a cylinder with PEC boundary conditions. The approximate wavelength is $\lambda = 0.25$, so that the wavelength to radius ratio is $\lambda/a = 1.25$. The dashed lines in the $E_x$ plot show the region where we include a PML layer. Note that only one vertical and one horizontal PML strip are required.} \label{Fig_TE_scattering}
\end{figure}

To test the error, we perform the full time-dependent simulation of the scattered wave up to time $T = 1.5080$.  We then compare the penalized solution with the exact analytic solution throughout the entire domain. The Lorentz-Mie-Debye method for electromagnetic scattering off of a perfectly conducting infinite cylinder is used to compute the exact solution for scattering of a time-harmonic plane wave in the frequency domain (see, for example, \cite{Harrington2001}). We then compute the time-dependent scattered solution at each grid point by taking the inverse Fourier transform of the exact time-harmonic solution scaled by the Fourier transform of the envelope of the plane wave. Since the Fourier spectrum of the Gaussian envelope is band-limited in finite precision, we can perform this step via inverse FFT with high accuracy. Figure \ref{Fig_scatter_convergence} shows the convergence plot of the error, while Figure \ref{Fig_TE_scattering} shows a plot of the scattered wave.

\subsection{Test 5: Two-dimensional bent waveguide}

Next, to demonstrate some potential uses for the penalization scheme, we treat a bent waveguide problem. To define the waveguide geometry,
we first construct five segments of a piecewise parametric curve with
\begingroup
\allowdisplaybreaks
\begin{align*}
\mathbf{r}_{1}\left(\tau\right) & =\left[\begin{array}{c}
l_{0}\tau\\
y_{\mathrm{off}}
\end{array}\right]\\
\mathbf{r}_{2}\left(\tau\right) & =\left[\begin{array}{c}
l_{0}+r_{0}\cos\left(\frac{\pi}{2}\tau+\pi\right)\\
y_{\mathrm{off}}+r_{0}\left(1+\sin\left(\frac{\pi}{2}\tau+\pi\right)\right)
\end{array}\right]\\
\mathbf{r}_{3}\left(\tau\right) & =\left[\begin{array}{c}
l_{0}+r_{0}\left(2-\cos\left(\pi\tau\right)\right)\\
y_{\mathrm{off}}+r_{0}\left(1+\sin\left(\pi\tau\right)\right)
\end{array}\right]\\
\mathbf{r}_{4}\left(\tau\right) & =\left[\begin{array}{c}
l_{0}+r_{0}\left(4+\cos\left(\frac{\pi}{2}\tau-\frac{\pi}{2}\right)\right)\\
y_{\mathrm{off}}+r_{0}\left(1+\sin\left(\frac{\pi}{2}\tau-\frac{\pi}{2}\right)\right)
\end{array}\right]\\
\mathbf{r}_{5}\left(\tau\right) & =\left[\begin{array}{c}
l_{0}+4r_{0}+l_{0}\left(\tau-4\right)\\
y_{\mathrm{off}}
\end{array}\right]
\end{align*}
\endgroup
such that 
\begin{equation}
\mathbf{r}\left(\tau\right)=\begin{cases}
\mathbf{r}_{1}\left(\tau\right), & 0\le\tau<1\\
\mathbf{r}_{2}\left(\tau\right), & 1\le\tau<2\\
\mathbf{r}_{3}\left(\tau\right), & 2\le\tau<3\\
\mathbf{r}_{4}\left(\tau\right), & 3\le\tau<4\\
\mathbf{r}_{5}\left(\tau\right), & 4\le\tau<5
\end{cases}.
\end{equation}
The constants $l_{0}$, $r_{0}$, and $y_{\mathrm{off}}$ correspond to the
length of straight line segments, the radii of circular arcs, and
the $y$-offset for the parametrized curve respectively. We then compute
a parametrized two-dimensional surface $\mathbf{S}\left(\tau,v\right)=\left(
x\left(\tau,v\right), y\left(\tau,v\right), \psi\left(\tau,v\right)\right)$ given by 
\begin{align}
\left[\begin{array}{c}
x\\
y
\end{array}\right] & =\mathbf{r}\left(\tau\right)+v\hat{\mathbf{n}}\left(\tau\right),\\
\psi & =c-\left|v\right|
\end{align}
where $c$ is some positive constant and 
\begin{equation}
\mathbf{\hat{n}}\left(\tau\right)=\left[\begin{array}{cc}
0 & 1\\
-1 & 0
\end{array}\right]\mathbf{t}\left(\tau\right)
\end{equation}
with 
\begin{equation}
\mathbf{t}\left(\tau\right)=\frac{\frac{d\mathbf{r}}{d\tau}\left(\tau\right)}{\left\Vert \frac{d\mathbf{r}}{d\tau}\left(\tau\right)\right\Vert _{2}}.
\end{equation}
The boundary of the bent waveguide corresponds to the zero level set
of $\psi$. In our example, we set $r_{0}=1$, $l_{0}=\pi-2r_{0}$,
$y_{\mathrm{off}}=2$. To avoid producing a multivalued function, we sweep
$\tau$ from 0 to 5 and $v$ from -1 to 1 and take $c=0.5$.
This is sufficient for our purposes as we are only really interested
in the signed distance function $\psi$ in the vicinity of the zero level
set (in particular, only signed distances $h$ and $-L$ away from the boundary).

Next, we require an expression for the normal to the boundary. We
first compute two tangent directions on $\mathbf{S}$ by differentiating
with respect to $\tau$ and $v$, and take their cross product. Projecting the resulting normal into the $xy$-plane and normalizing to unit length yields 
\begin{equation}
\mathbf{n} = \mathrm{sign}\left(v\right)\hat{\mathbf{n}}\left(\tau\right).
\end{equation}

We note that for a given $\tau$ and $v$, we can evaluate the corresponding
location $\left(x,y\right)$, the level set value $\psi$, and its
corresponding normal $\mathbf{n}$. However, in our setting, we require the value
of $\psi$ and its corresponding normal at a set of known locations
(the grid points and a set of boundary points). Thus, we interpolate
from equally spaced data in the $\tau v$-plane (which is not equally
spaced in the $xy$-plane) to the grid points and necessary boundary
points. We do this once in the pre-processing stage of the algorithm
before we begin our time-stepping scheme.

\subsubsection{$TM_z$ mode manufactured solution}

\begin{figure}[htb!]
	\centering
    \includegraphics[width = \textwidth]{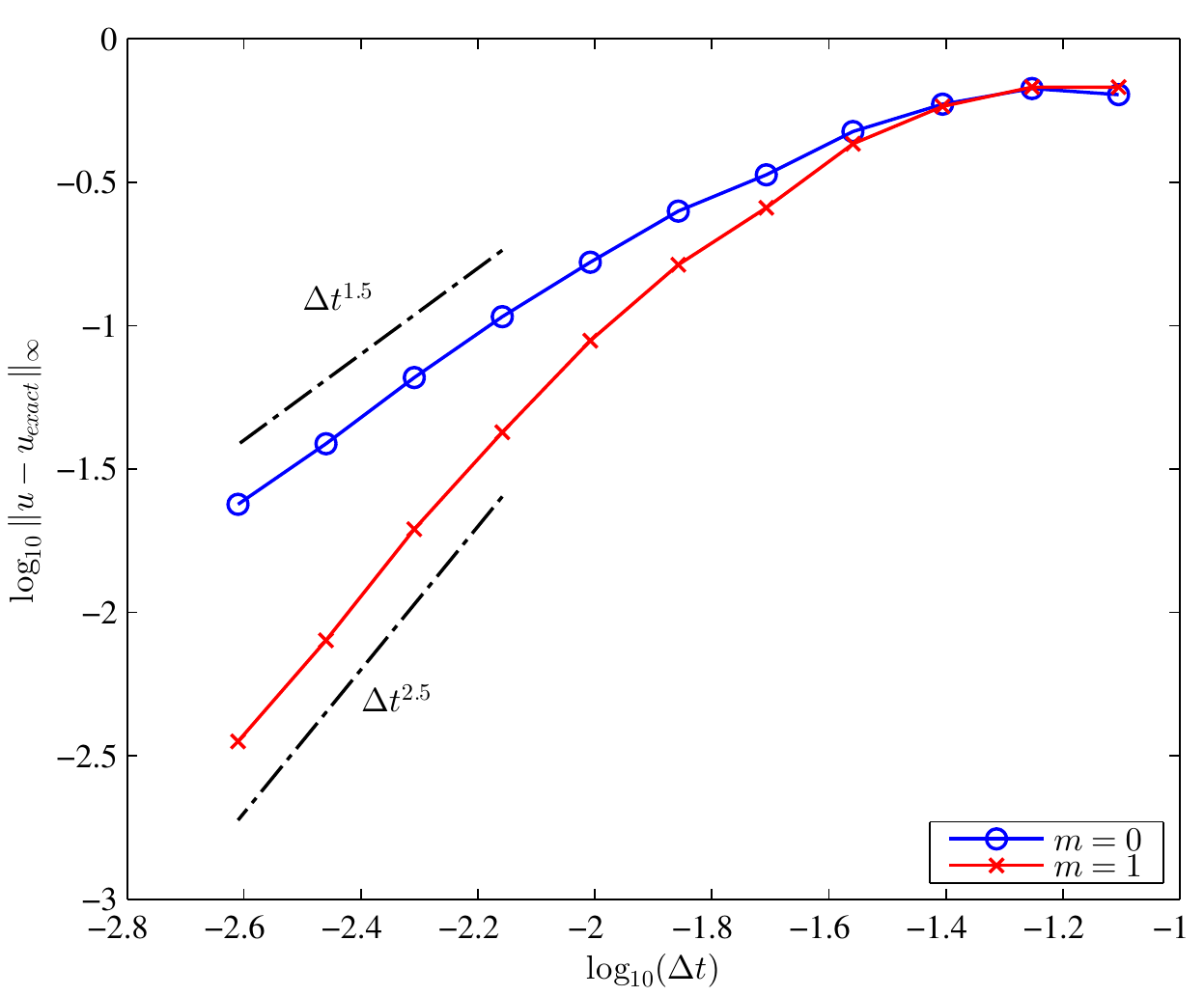} \\
    \caption{Convergence study for the manufactured $\textrm{TM}_z$ example applied to the bent waveguide geometry. The global convergence rates in $L^{\infty}(\Omega_0)$ are approximately 1.5 and 2.5 for $m=0$ and $m=1$ respectively. This agrees with our previous test case. Here, $u=(H_{x,\eta},H_{y,\eta},E_{z,\eta})$.} \label{TM_waveguide_convergence}
\end{figure}

As in Section \ref{section_TM_manufacture}, we first verify that our construction converges using a manufactured solution approach.
In fact, we use the same manufactured solution as before and only change the geometry of the obstacle. In addition, due to the curvature of the waveguide boundary, we are required to take a smaller decay length $L=0.4$ for $\tilde{\mathbf{g}}$. Otherwise, with the exception of $T=0.275\pi$ and $\eta = \Delta t$, all other parameters are left unchanged. Figure \ref{TM_waveguide_convergence} illustrates the same convergence rates as in the previous test and demonstrates the validity of the waveguide construction.

\subsubsection{$TE_z$ mode plane wave propagation}

\begin{figure}[p!] 
\centering $
\begin{array}{cc}
    \includegraphics[width = .45\textwidth]{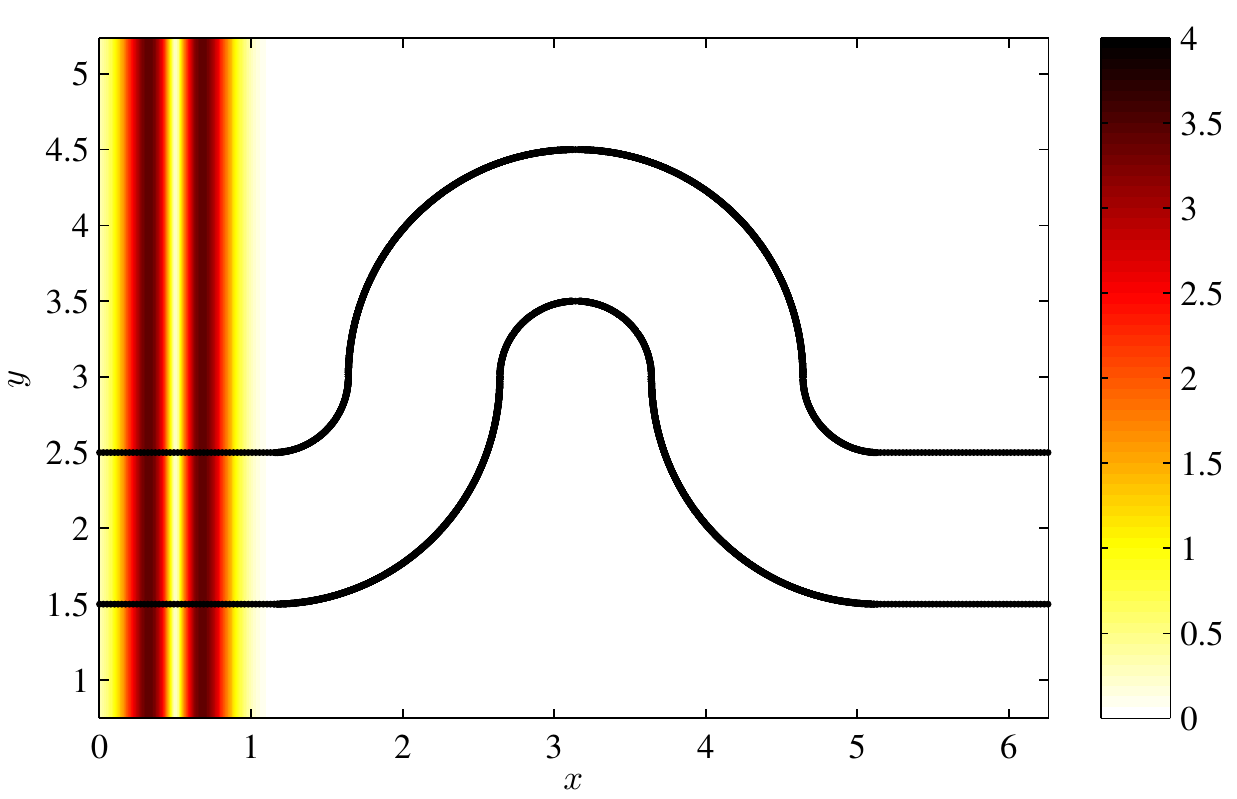} 
    \includegraphics[width = .45\textwidth]{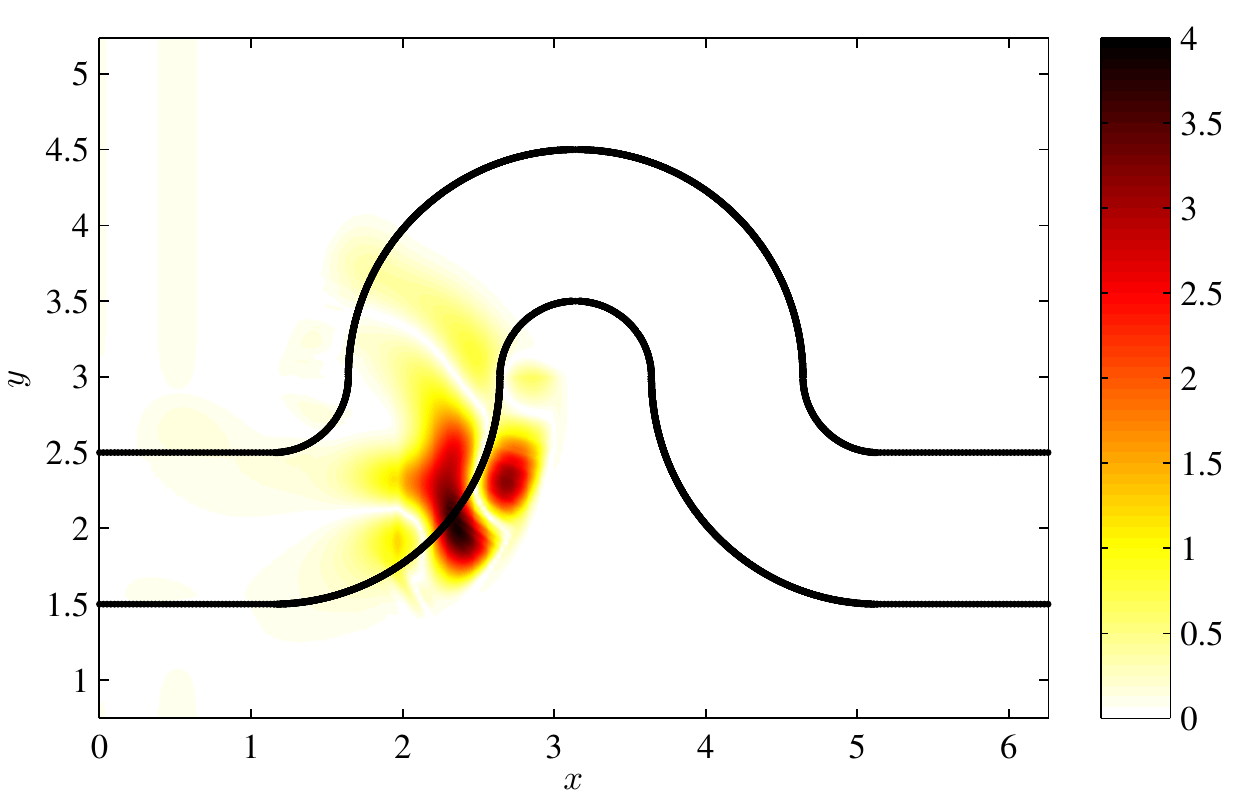} \\
    \includegraphics[width = .45\textwidth]{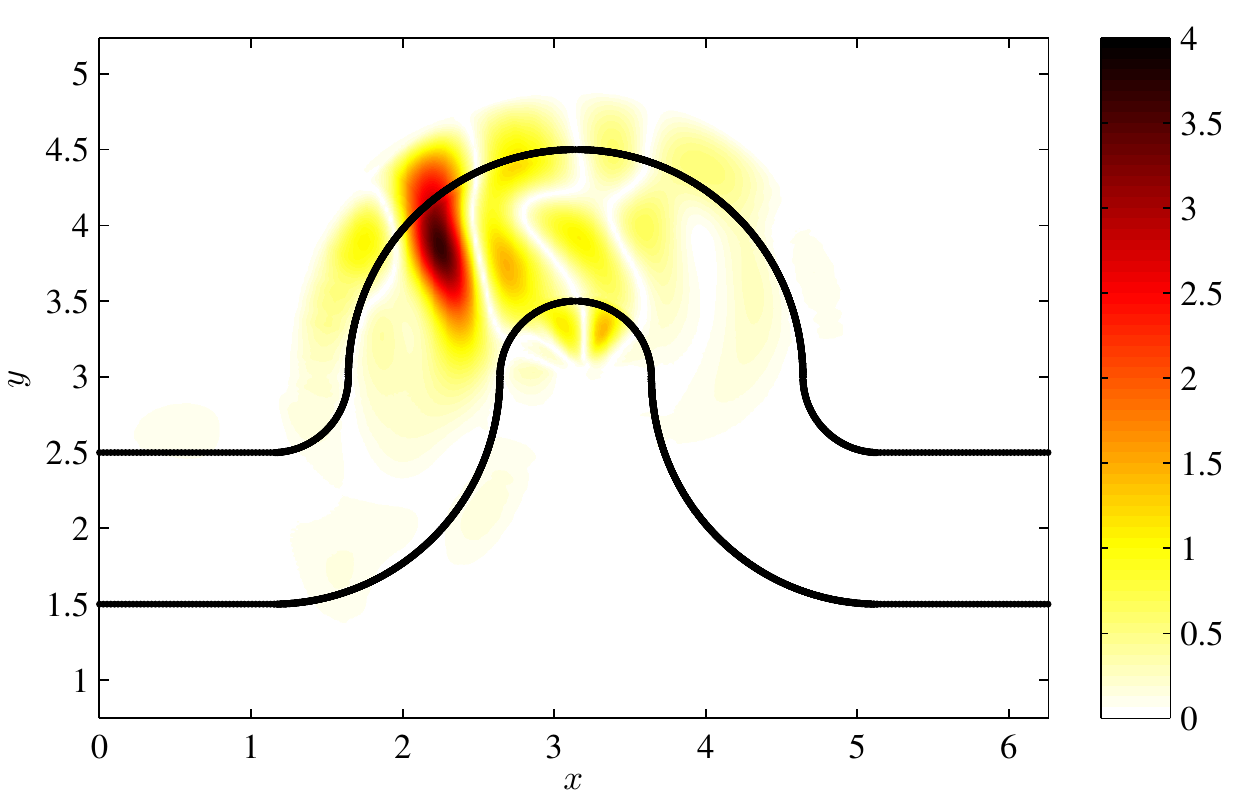} 
    \includegraphics[width = .45\textwidth]{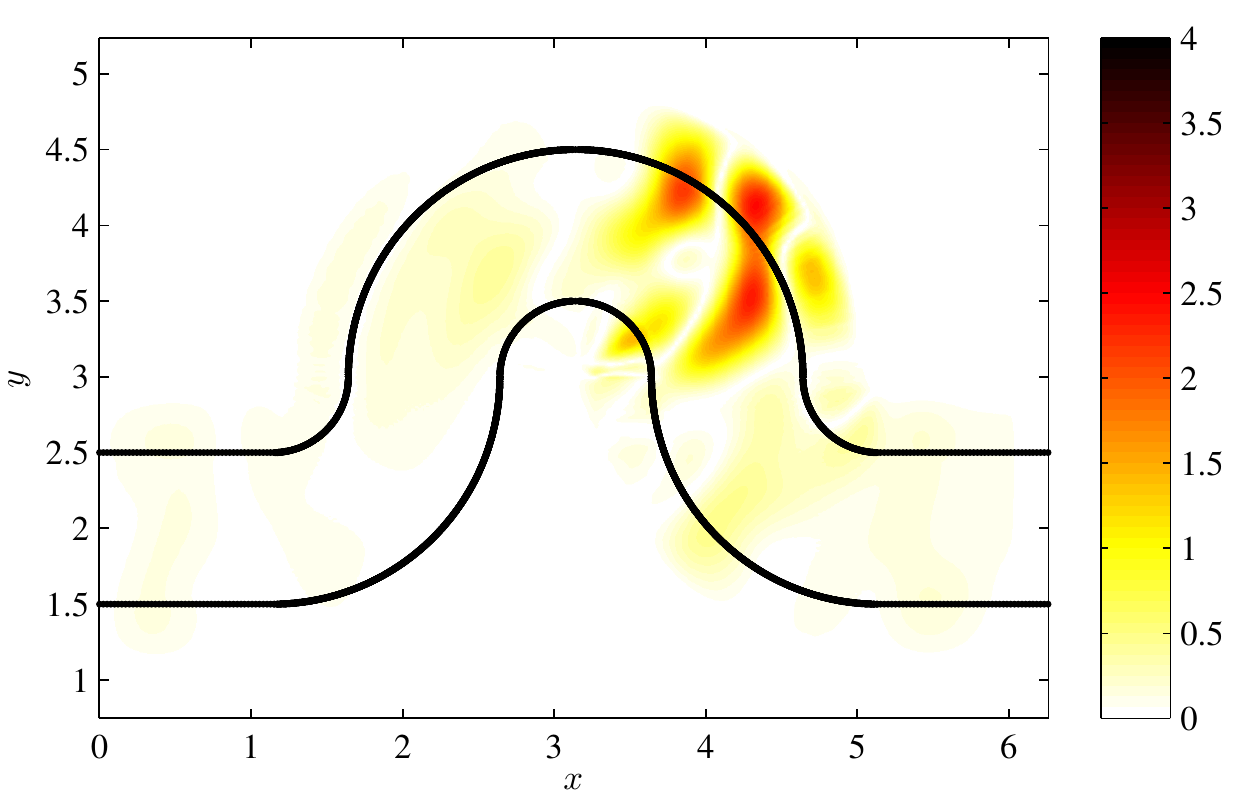} \\
    \includegraphics[width = .45\textwidth]{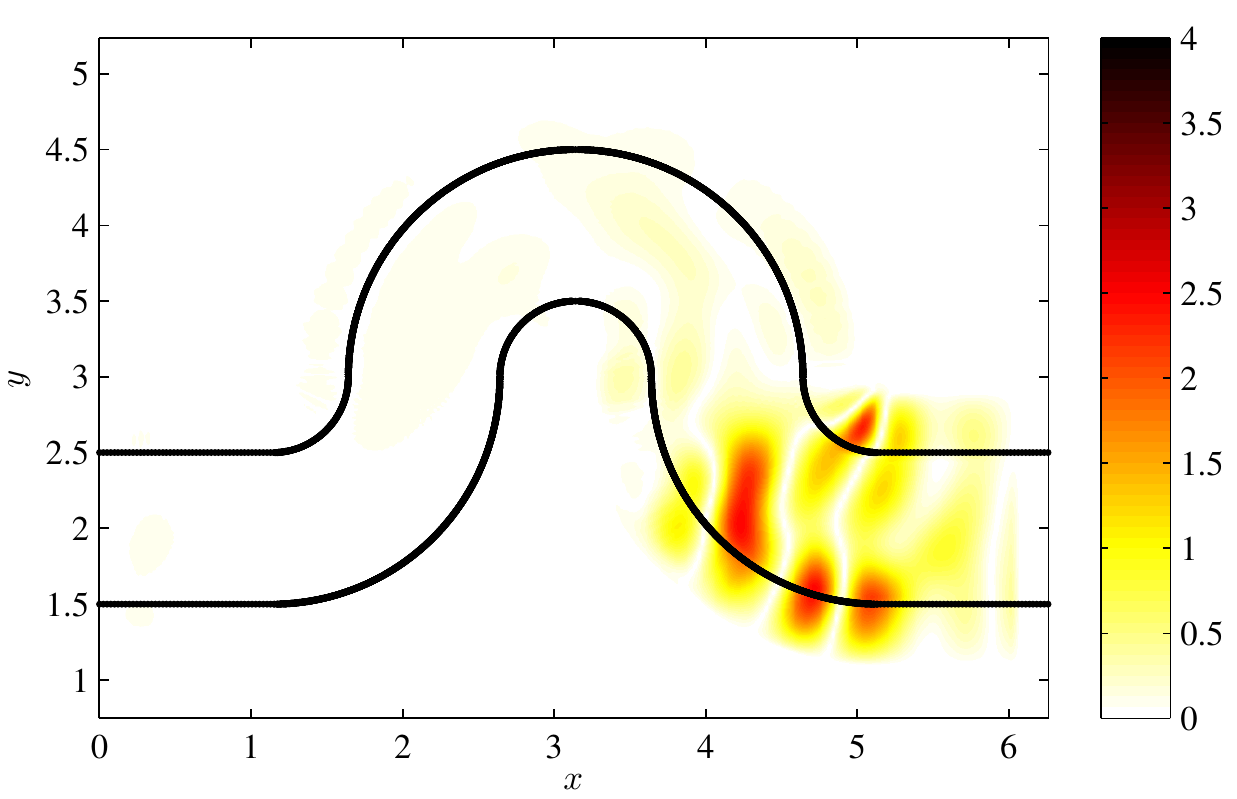} 
    \includegraphics[width = .45\textwidth]{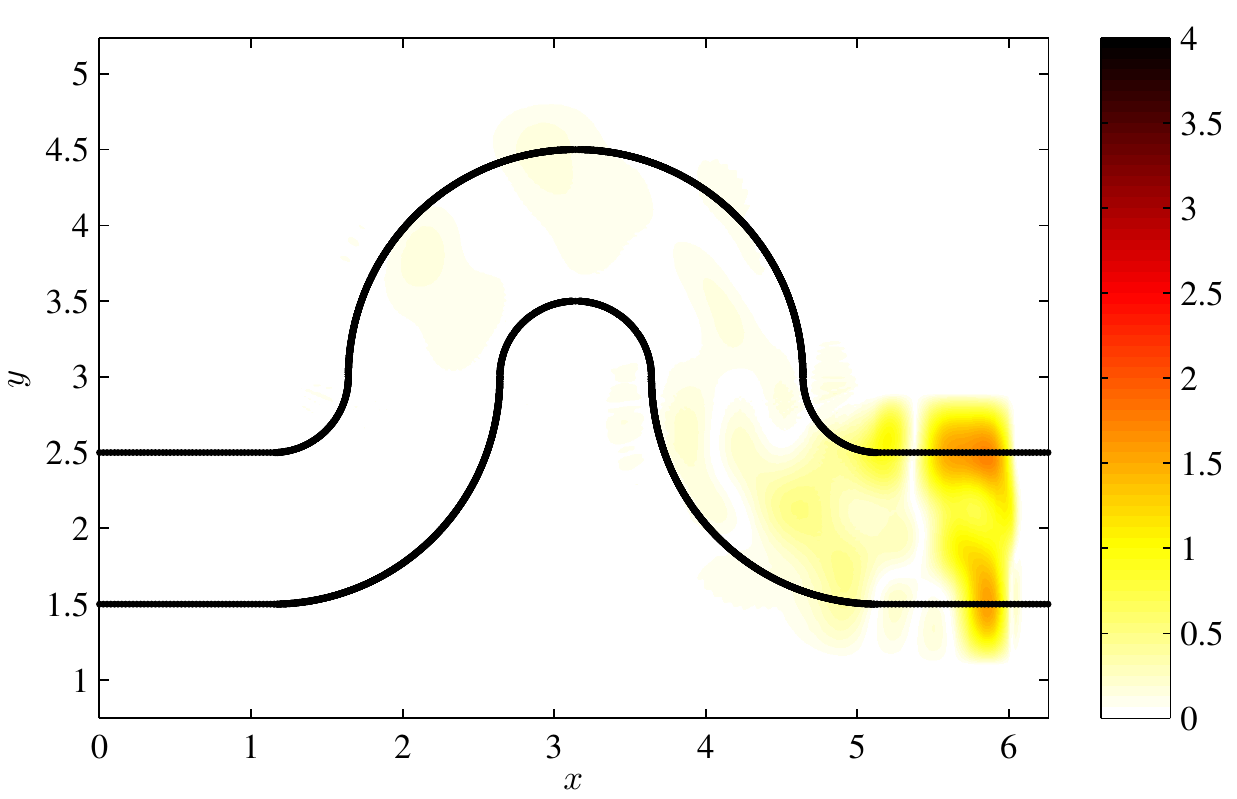}
\end{array}$
\caption{Six snapshots of the magnitude of the $E_{x,\eta}$ component of the plane wave propagating through the bent waveguide. The initial condition ($t=0$) is shown in the top left corner, with subsequent figures taken at times $t=2,4,6,8,10$.} \label{waveguide_snapshots}
\end{figure}

More practically, consider the same waveguide geometry, but with an initial condition corresponding to a pulsed Gaussian. Take, for example, the initial conditions
\begin{align}
E_{x}\left(x,y,0\right) & = \phantom{-}0\\
E_{y}\left(x,y,0\right) & =-\frac{2}{\sigma^2} \left( x - x_0 \right) e^{-\left( \frac{ x - x_0}{\sigma} \right)^2}\\
H_{z}\left(x,y,0\right) & =-\frac{2}{\sigma^2} \left( x - x_0 \right) e^{-\left( \frac{ x - x_0}{\sigma} \right)^2}
\end{align}
with $\sigma = 0.25$ and $x_0 = 0.5$.
We may then solve equations (\ref{perfectly_matched_layers}) with these initial conditions and the splitting $H_{zx}(x,y,0) = H_z(x,y,0)$ and $H_{zy}(x,y,0) = 0$. We take $\sigma_{x,max} = D / 2\Delta x$ for a slab of width $0.25$ and set $\sigma_y = 0$. In addition, $T=10$, $\Delta t =0.4 \Delta x$, $h=2 \Delta x$, and $\eta = \Delta t$ with $m=0$. In Figure \ref{waveguide_snapshots}, we plot a collection of snapshots of the behavior of the plane wave as it propagates down the waveguide.

\subsection{Test 6: Two-dimensional scattering off a windmill-like geometry}

To demonstrate scattering from objects that are not comprised of circular
boundaries, we consider a windmill-like geometry adapted from a rhodonea
curve given, in polar coordinates, by $r=a\sin\left(3\theta\right)$.
This trifolium is an algebraic curve corresponding to the zero level set of
\begin{equation}
\tilde{\psi}\left(x,y\right)=\left(x^{2}+y^{2}\right)^{2}-4ayx^{2}+ay\left(x^{2}+y^{2}\right).
\end{equation}
To obtain a single smooth boundary, we shift $\tilde{\psi}$ and work
with $\hat{\psi}\left(x,y\right)=\tilde{\psi}(x,y)-b$
with $a=3$ and $b=1$ fixed. This level set function is not a signed distance
function (one can check that $|\nabla\hat{\psi}|\ne1$)
so we must construct, as outlined in Section \ref{Sec_NumSolveCoordinates}, a signed distance function $\psi$ whose zero level set coincides with that of $\hat{\psi}$.

To illustrate one possible scattering solution in the vicinity of this windmill-like geometry, let us consider the $\textrm{TM}_z$ mode. We penalize the equations as described in Section \ref{Sec_PenaltyFunction_TMmode}. We then add
a PML to absorb outgoing scattered waves.
Unlike the PML discussed in Section \ref{sec_perfectly_matched_layers}, we use the complex coordinate stretching interpretation
of the PML \cite{ChewWeedon1994, Rappaport1995, TeixeiraChew1998} to avoid splitting the penalization term in our equations. We decompose solutions to Maxwell's equations
into terms of the form $\hat{E}_{z,\eta}\left(x,y\right)e^{\imath \omega t}$
(respectively $\hat{H}_{x,\eta}\left(x,y\right)e^{\imath \omega t}$ and $\hat{H}_{y,\eta}\left(x,y\right)e^{\imath \omega t}$)
and write the $\textrm{TM}_{z}$ mode equations in the frequency domain. We
then replace 
\begin{align}
\frac{\partial}{\partial x} & \rightarrow\frac{1}{1-\imath\sigma_{x}\omega^{-1}} \frac{\partial}{\partial x}\\
\frac{\partial}{\partial y} & \rightarrow\frac{1}{1-\imath\sigma_{y}\omega^{-1}} \frac{\partial}{\partial y},
\end{align}
multiply both sides by the denominators, ignore all terms containing
products of $\sigma_{x}$ or $\sigma_{y}$ with $\chi_h$ (this is valid
as long as the penalization and PML regions do not overlap), and finally
transform back to the time domain. The resulting equations are 
\begin{subequations}
\begin{align} \label{TM_PML_begin}
\frac{\partial H_{x,\eta}}{\partial t} & =-\frac{\partial E_{z,\eta}}{\partial y}-\sigma_{y}H_{x,\eta}\\
\frac{\partial H_{y,\eta}}{\partial t} & =\phantom{-}\frac{\partial E_{z,\eta}}{\partial x}-\sigma_{x}H_{y,\eta}\\
\frac{\partial E_{z,\eta}}{\partial t} & =-\frac{\partial H_{x,\eta}}{\partial y}+\frac{\partial H_{y,\eta}}{\partial x}-\eta^{-1}\chi_h\left(\mathbf{x}\right)\left(E_{z,\eta}-g_{z}\right)-\left(\sigma_{x}+\sigma_{y}\right)E_{z,\eta}+\Phi\\ \label{TM_PML_end}
\frac{\partial \Phi}{\partial t} & =-\sigma_{x}\sigma_{y}E_{z,\eta},
\end{align}
\end{subequations}
where $\Phi$ is an auxiliary variable (initialized to zero) added
to avoid integrals (terms of the form $-\imath\omega^{-1}$ in the
frequency domain) in the time domain representation.

\begin{figure}[htbp!]
\centering $
\begin{array}{cc}
    \includegraphics[width = .43\textwidth]{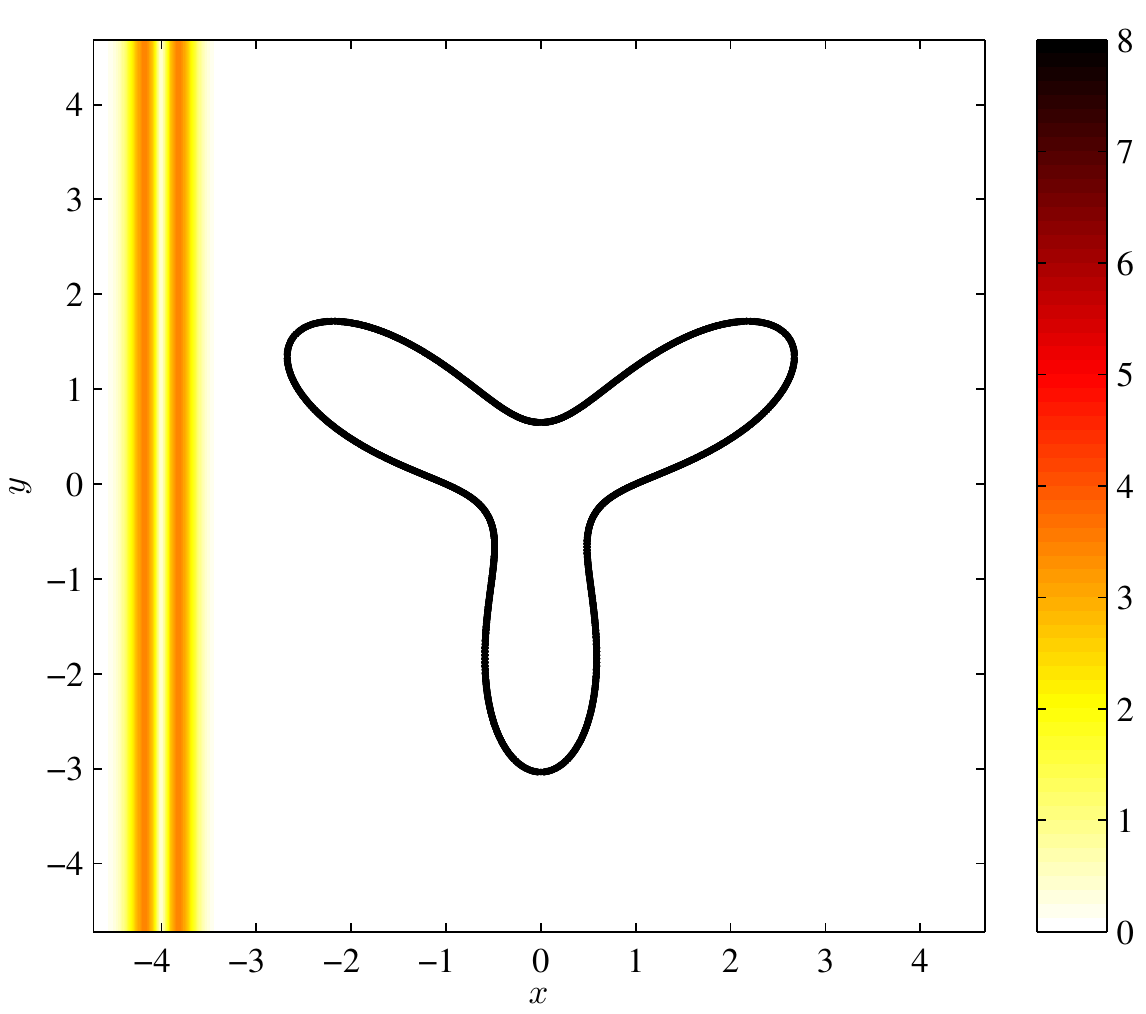} 
    \includegraphics[width = .43\textwidth]{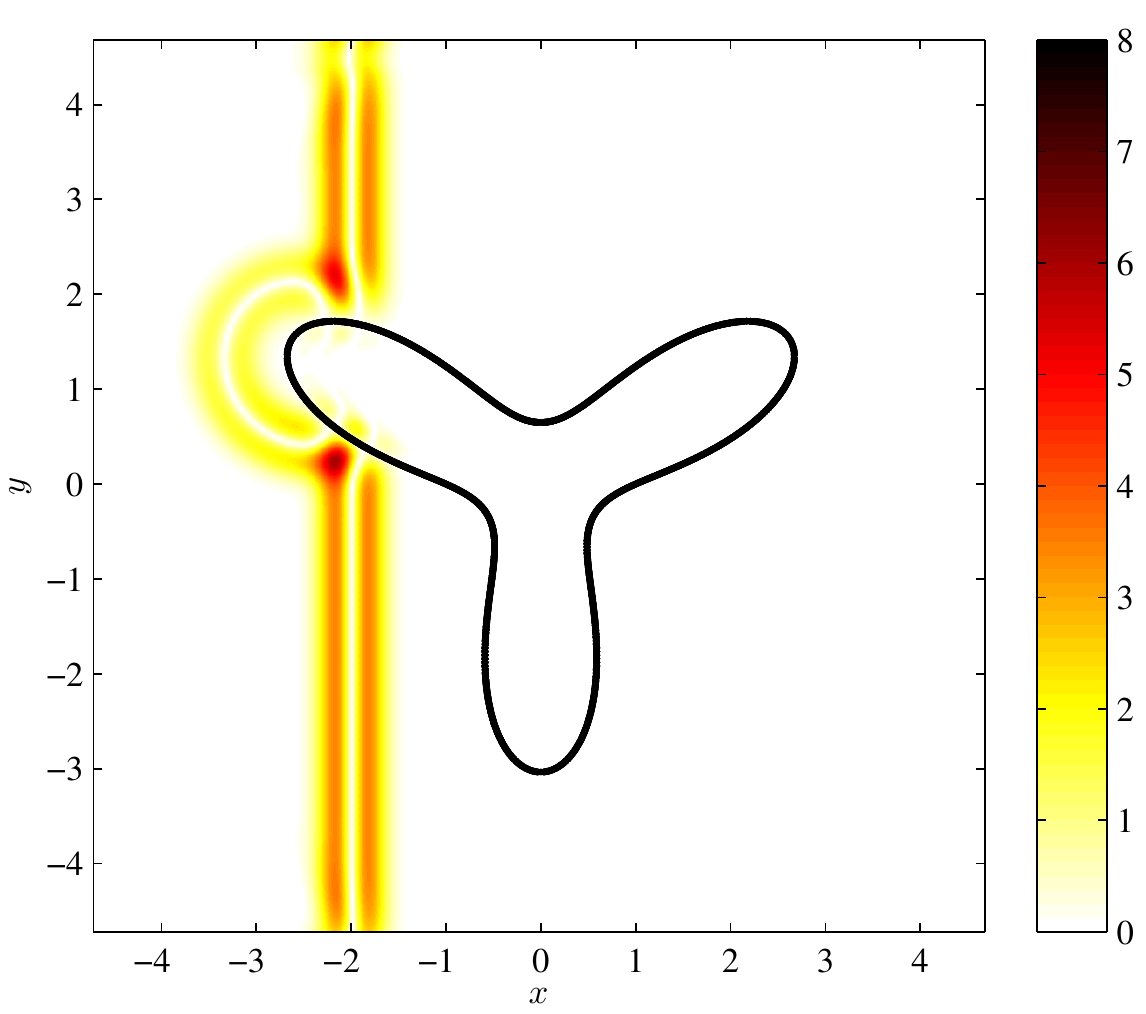} \\
    \includegraphics[width = .43\textwidth]{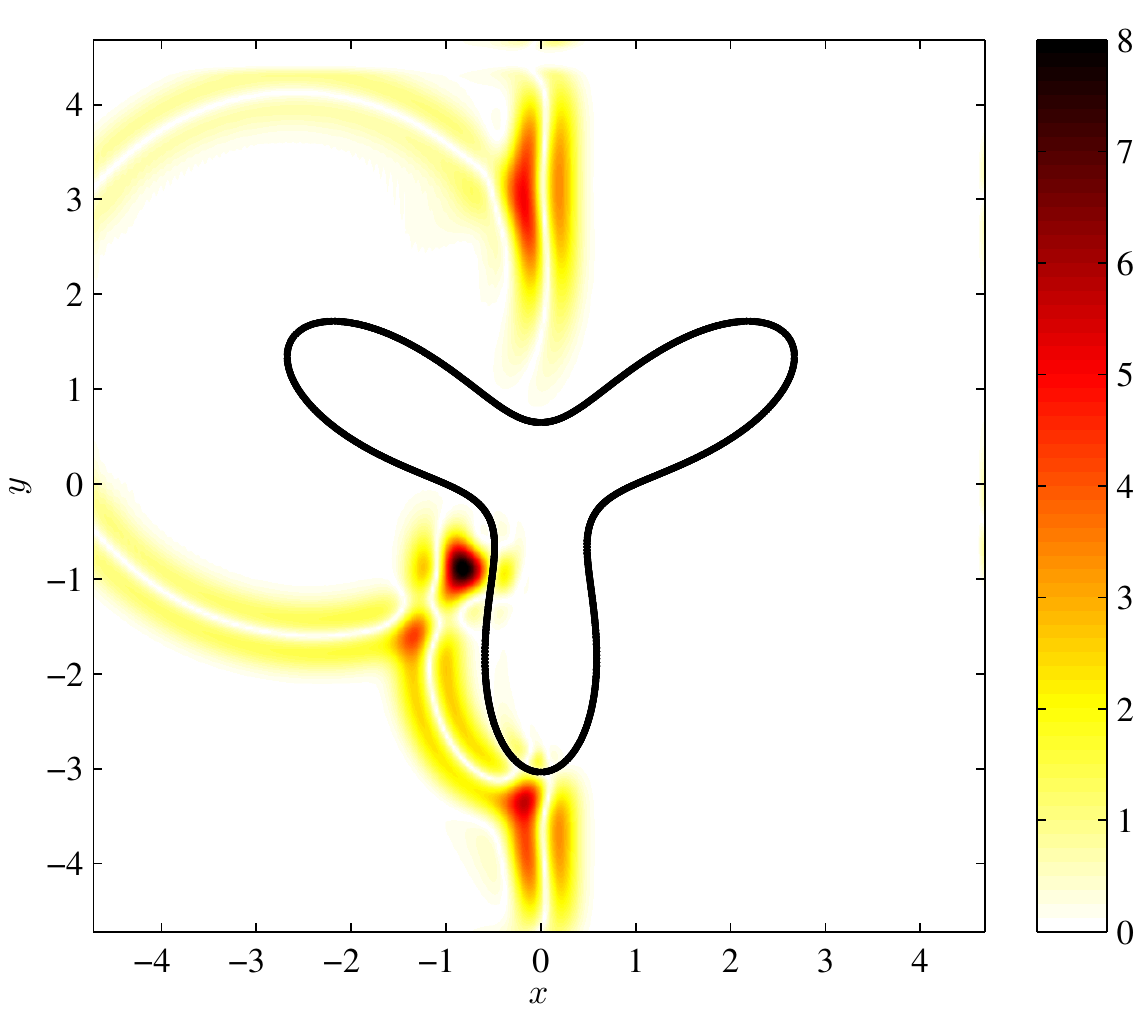} 
    \includegraphics[width = .43\textwidth]{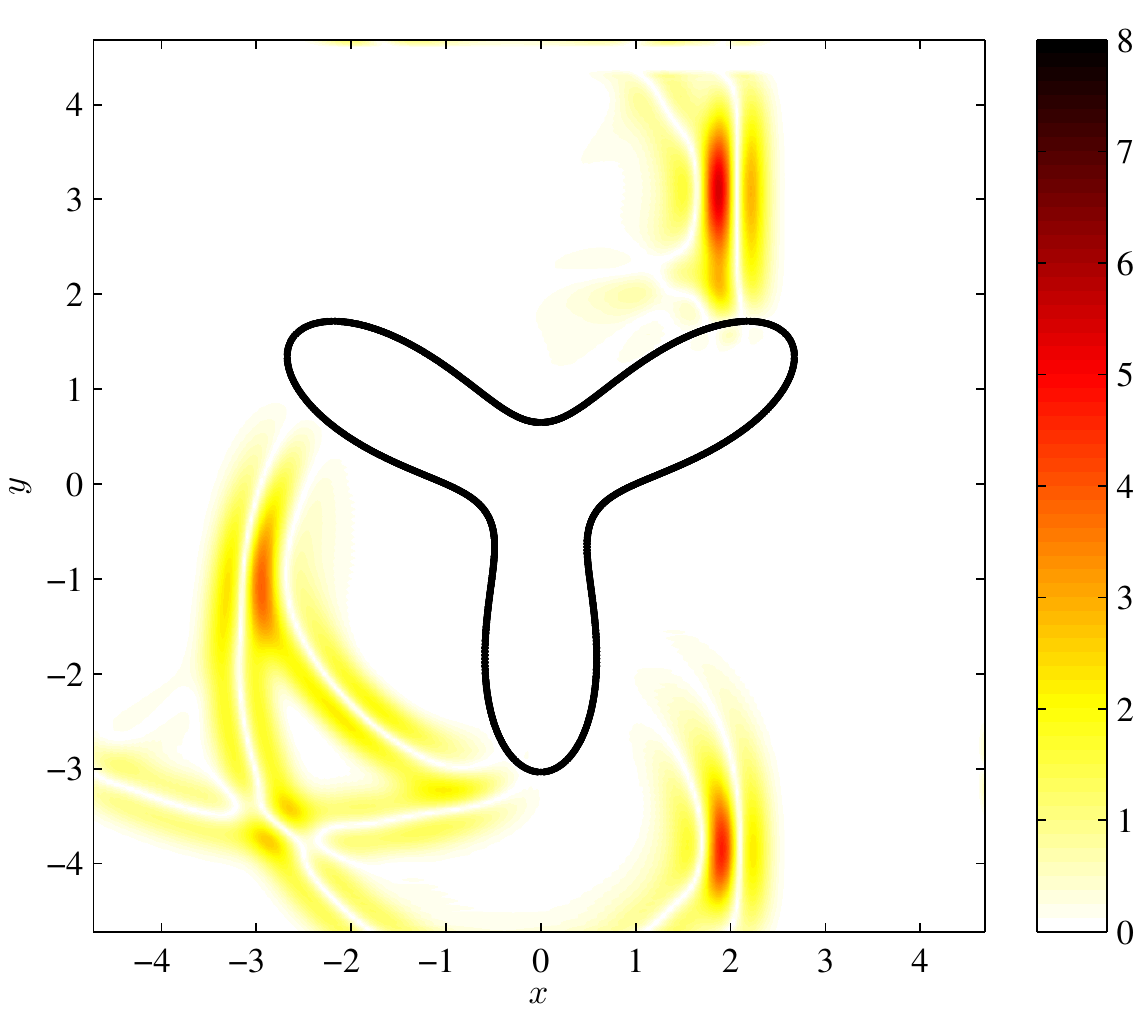} \\
    \includegraphics[width = .43\textwidth]{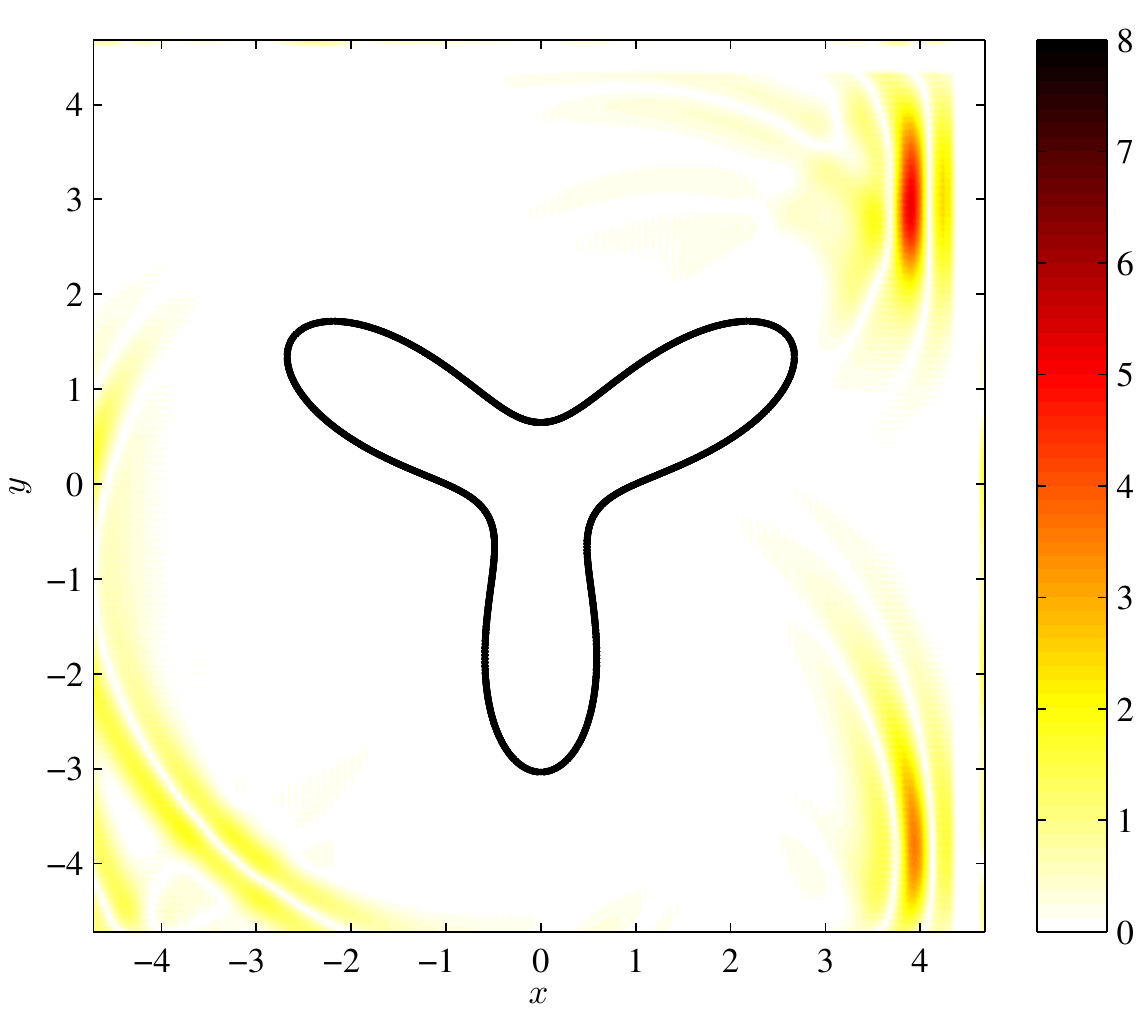} 
    \includegraphics[width = .43\textwidth]{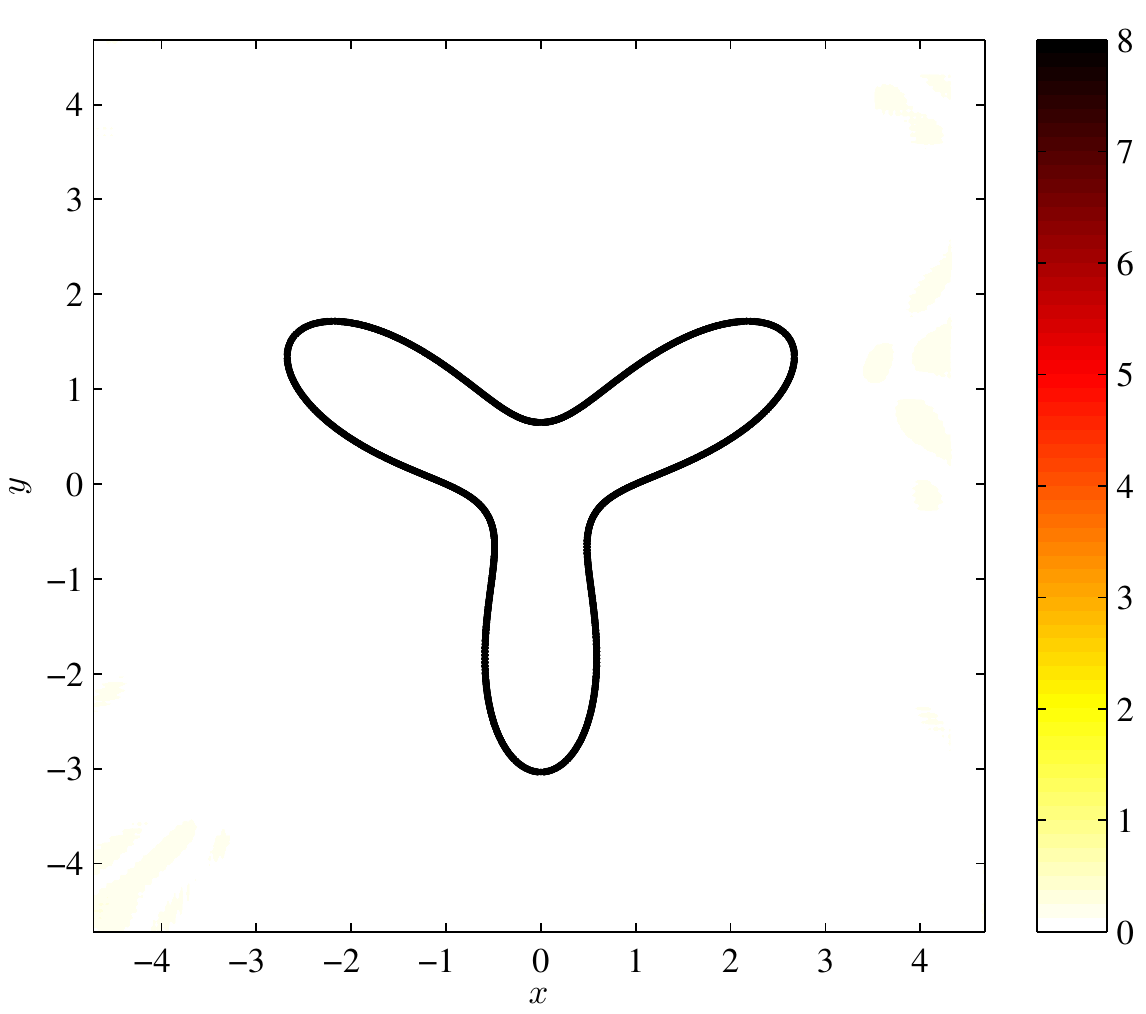}
\end{array} $
\caption{Six snapshots of the magnitude of the $E_{z,\eta}$ component of a plane wave scattering off the windmill-like geometry. The initial condition ($t=0$) is shown in the top left corner, with subsequent figures taken at times $t=2,4,6,8,10$. By the final snapshot, the reflected wave is on the verge of having completely left the computational domain.} \label{windmill_snapshots}
\end{figure}

Discretization in space yields the equations
\begin{subequations}
\begin{align}
\frac{\partial H_{x,\eta}}{\partial t} & =- \mathcal{F}^{-1}\left\{ \imath k_{y}\mathcal{F}\left\{ E_{z,\eta}\right\} \right\} -\sigma_{y}H_{x,\eta}\\
\frac{\partial H_{y,\eta}}{\partial t} & =\phantom{-} \mathcal{F}^{-1}\left\{ \imath k_{x}\mathcal{F}\left\{ E_{z,\eta}\right\} \right\} -\sigma_{x}H_{y,\eta}\\
\frac{\partial E_{z,\eta}}{\partial t} & =-\mathcal{F}^{-1}\left\{ \imath k_{y}\mathcal{F}\left\{ H_{x,\eta}\right\} \right\}+\mathcal{F}^{-1}\left\{ \imath k_{x}\mathcal{F}\left\{ H_{y,\eta}\right\} \right\} \\
&\phantom{=} \qquad -\eta^{-1}\chi_h\left(\mathbf{x}\right)\left(E_{z,\eta}-g_{z}\right)-\left(\sigma_{x}+\sigma_{y}\right)E_{z,\eta}+\Phi \notag \\
\frac{\partial \Phi}{\partial t} & =-\sigma_{x}\sigma_{y}E_{z,\eta}
\end{align}
\end{subequations}
which are integrated forward in time using RK4. We solve the problem on $\Omega = [-\frac{3\pi}{2},\frac{3\pi}{2}]^2$ with initial conditions
\begin{subequations} \label{Gaussian_pulse}
\begin{align}
H_{x}\left(x,y,0\right) & =\phantom{-} 0\\
H_{y}\left(x,y,0\right) & =-\frac{2}{\sigma^{2}}\left(x-x_{0}\right)e^{-\left(\frac{x-x_{0}}{\sigma}\right)^{2}}\\
E_{z}\left(x,y,0\right) & =\phantom{-} \frac{2}{\sigma^{2}}\left(x-x_{0}\right)e^{-\left(\frac{x-x_{0}}{\sigma}\right)^{2}} 
\end{align}
\end{subequations}
which corresponds to a pulsed wave traveling in the positive $x$-direction. For our example, $x_0 = -4$ and $\sigma = 0.25$. We take $\sigma_{x,max} = \sigma_{y,max} = N/2$ in slabs of width $1/3$. Unlike the $\textrm{TE}_z$ mode, we begin our simulation with both PMLs set to zero, and increase their values over half the duration of the simulation (we do so cubically). This is done to avoid nonphysical reflections of the initial plane wave off of the PML. Finally, we set $T=10$, $\Delta t =0.4 \Delta x$, $h=2 \Delta x$, $\eta = \Delta t$, and $L=0.3$. Table \ref{windmill_snapshots} illustrates the plane wave scattering off the windmill-like obstacle using $m=1$ normal derivatives matched at the boundary.

\subsection{Test 7: Three-dimensional manufactured solution for a domain with a spherical hole}

To demonstrate the applicability of the method in three dimensions, we first test a manufactured standing wave solution. We seek solutions to
\begin{subequations} \label{Maxwell_3D}
\begin{align}
\frac{\partial H_{x}}{\partial t} & = \frac{\partial E_{y}}{\partial z} - \frac{\partial E_{z}}{\partial y} \\
\frac{\partial H_{y}}{\partial t} & = \frac{\partial E_{z}}{\partial x} - \frac{\partial E_{x}}{\partial z} \\
\frac{\partial H_{z}}{\partial t} & = \frac{\partial E_{x}}{\partial y} - \frac{\partial E_{y}}{\partial x} \\
\frac{\partial E_{x}}{\partial t} & = \frac{\partial H_{z}}{\partial y} - \frac{\partial H_{y}}{\partial z} \\
\frac{\partial E_{y}}{\partial t} & = \frac{\partial H_{x}}{\partial z} - \frac{\partial H_{z}}{\partial x} \\ 
\frac{\partial E_{z}}{\partial t} & = \frac{\partial H_{y}}{\partial x} - \frac{\partial H_{x}}{\partial y}
\end{align}
\end{subequations}
with initial conditions
\begin{subequations} \label{initial_conditions_3D}
\begin{align}
\mathbf{E}\left(\mathbf{x},\frac{\pi}{2\sqrt{3}} \right) & = 0 \\
\mathbf{H}\left(\mathbf{x},\frac{\pi}{2\sqrt{3}} \right) & = 2 \left(\mathbf{k} \times \mathbf{E}_0 \right) \sin{ ( \sqrt{3} \, \mathbf{k} \cdot \mathbf{x} )}.
\end{align}
\end{subequations}
The solution to (\ref{Maxwell_3D}) subject to initial conditions (\ref{initial_conditions_3D}) is given by
\begin{subequations}
\begin{align}
\mathbf{E}\left(\mathbf{x},t \right) & = 2 \mathbf{E}_0 \cos{( \sqrt{3} \, t )} \cos{( \sqrt{3} \, \mathbf{k} \cdot \mathbf{x} )}\\
\mathbf{H}\left(\mathbf{x},t \right) & = 2 \left(\mathbf{k} \times \mathbf{E}_0 \right) \sin{( \sqrt{3} \, t )} \sin{( \sqrt{3} \, \mathbf{k} \cdot \mathbf{x} )},
\end{align}
\end{subequations}
which one can verify satisfies the divergence-free criteria and is periodic when $\mathbf{k} = \frac{1}{\sqrt{3}}(1,1,1)$ and $\mathbf{E}_0 = (1,-2,1)$.

\begin{figure}[htb!]
	\centering
    \includegraphics[width = \textwidth]{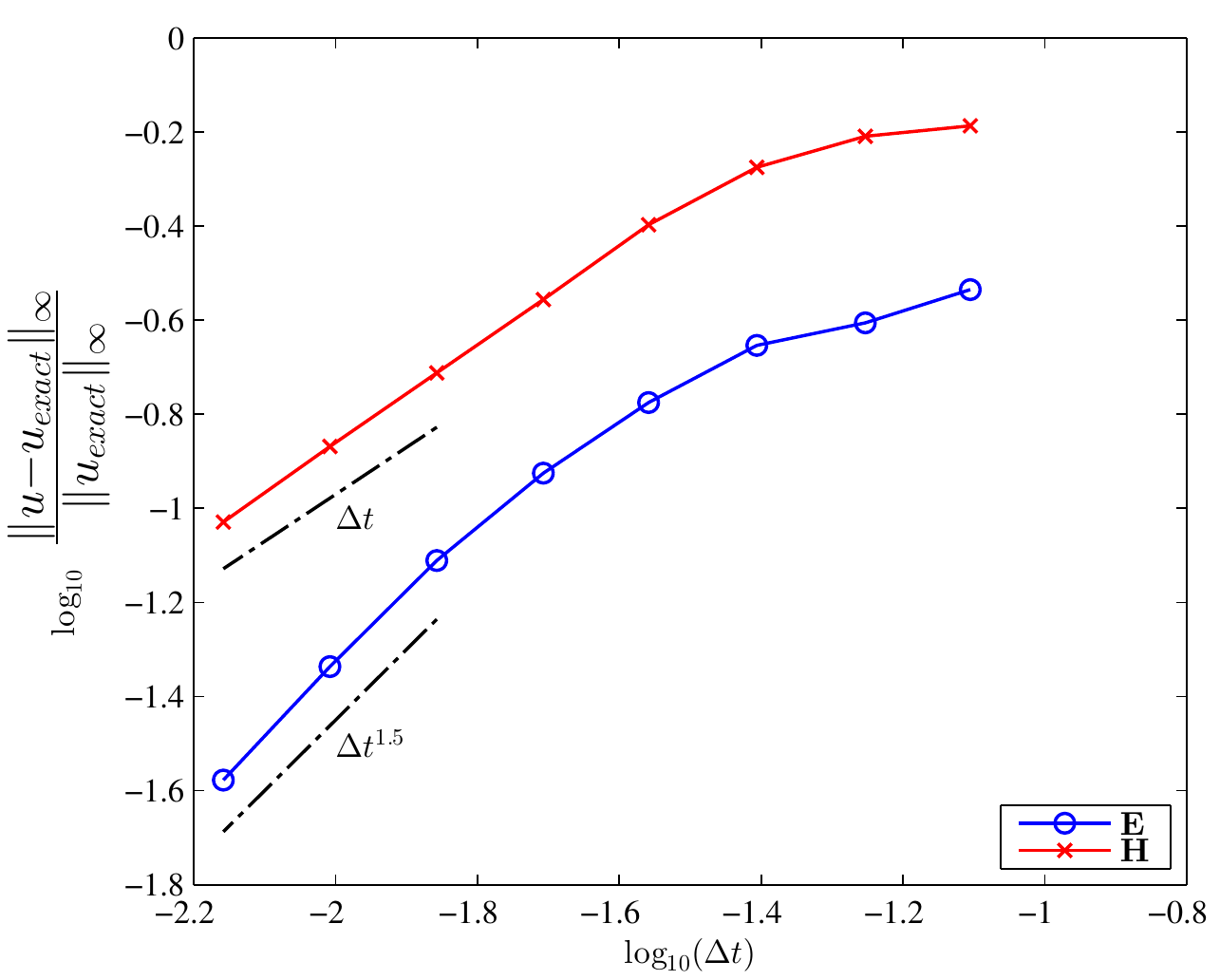} \\
    \caption{Convergence study for the full three-dimensional Maxwell equations in the presence of the spherical geometry. The global convergence rate in $L^{\infty}(\Omega_0)$ is $1.5$ for $\mathbf{E}$ and $1$ for $\mathbf{H}$. Here, $u$ is either $\mathbf{H}_{\eta}$ or $\mathbf{E}_{\eta}$.} \label{3D_convergence}
\end{figure}

We solve these equations on the periodic domain $\Omega=\left[0,2\pi\right]^3$ with a spherical hole removed.  Specifically, the boundary $\Gamma$ of the hole is given by the zero level set of the signed distance function
\begin{equation}
\psi\left(x,y,z\right)=\sqrt{\left(x-x_{0}\right)^{2}+\left(y-y_{0}\right)^{2}+\left(z-z_{0}\right)^{2}}-a \label{eq:levelset_3D}
\end{equation}
with radius $a$ and center $(x_0,y_0,z_0)$. In our tests, we fix $a=2$ with center $x_{0}=y_{0}=z_{0}=\pi$.
the penalized solution is computed by integration in time (using RK4) of 
\begin{subequations}
\begin{align}
\frac{\partial H_{x}}{\partial t} & = \mathcal{F}^{-1}\left\{ \imath k_{z}\mathcal{F}\left\{ E_{y,\eta}\right\} \right\} - \mathcal{F}^{-1}\left\{ \imath k_{y}\mathcal{F}\left\{ E_{z,\eta}\right\} \right\} \\
\frac{\partial H_{y}}{\partial t} & = \mathcal{F}^{-1}\left\{ \imath k_{x}\mathcal{F}\left\{ E_{z,\eta}\right\} \right\} - \mathcal{F}^{-1}\left\{ \imath k_{z}\mathcal{F}\left\{ E_{x,\eta}\right\} \right\} \\
\frac{\partial H_{z}}{\partial t} & = \mathcal{F}^{-1}\left\{ \imath k_{y}\mathcal{F}\left\{ E_{x,\eta}\right\} \right\} - \mathcal{F}^{-1}\left\{ \imath k_{x}\mathcal{F}\left\{ E_{y,\eta}\right\} \right\} \\
\frac{\partial E_{x}}{\partial t} & = \mathcal{F}^{-1}\left\{ \imath k_{y}\mathcal{F}\left\{ H_{z,\eta}\right\} \right\} - \mathcal{F}^{-1}\left\{ \imath k_{z}\mathcal{F}\left\{ H_{y,\eta}\right\} \right\} -\eta^{-1}\chi_{h}\left(x,y,z\right)\left(E_{x,\eta}-\tilde{g}_{x}\right) \\
\frac{\partial E_{y}}{\partial t} & = \mathcal{F}^{-1}\left\{ \imath k_{z}\mathcal{F}\left\{ H_{x,\eta}\right\} \right\} - \mathcal{F}^{-1}\left\{ \imath k_{x}\mathcal{F}\left\{ H_{z,\eta}\right\} \right\}  -\eta^{-1}\chi_{h}\left(x,y,z\right)\left(E_{y,\eta}-\tilde{g}_{y}\right) \\ 
\frac{\partial E_{z}}{\partial t} & = \mathcal{F}^{-1}\left\{ \imath k_{x}\mathcal{F}\left\{ H_{y,\eta}\right\} \right\} - \mathcal{F}^{-1}\left\{ \imath k_{y}\mathcal{F}\left\{ H_{x,\eta}\right\} \right\}  -\eta^{-1}\chi_{h}\left(x,y,z\right)\left(E_{z,\eta}-\tilde{g}_{z}\right)
\end{align}
\end{subequations}
where $\tilde{\mathbf{g}}$ is constructed using $\mathbf{g} = \mathbf{E}$. Figure \ref{3D_convergence} shows the convergence of the penalized solution to the exact solution at final time $T=\frac{\pi}{2\sqrt{3}} + 3$ with $\Delta t = 0.4 \Delta x$, $h=2 \Delta x$, $\eta = 4 \Delta t$, and $L = 1$. In three dimensions (as opposed to the $\textrm{TE}_z$ mode in two dimensions), $H_x$, $H_y$, and $E_z$ are not zero. As a result, we expect a slight decrease in accuracy for $\mathbf{H}$ as, in our approach, only $\mathbf{E}$ is penalized. This is indeed observed for the case $m=0$, where the convergence rate for $\mathbf{E}$ is $1.5$ and the convergence rate for $\mathbf{H}$ is $1$. One possible improvement could be to design a more complicated penalization involving $\mathbf{H}$ to reconcile its convergence rate with that of $\mathbf{E}$ in the three-dimensional case.


\subsection{Test 8: Three-dimensional scattering off a gyroid}

As a final example, we consider periodic scattering of a radiating dipole off of a gyroid. We solve the full three-dimensional Maxwell equations \eqref{Maxwell_3D} on the domain $\Omega = [0,1]^3$ with initial conditions corresponding to an ideal dipole \cite{Harrington2001} whose initial radial envelope in the azimuthal plane corresponds to a Gaussian pulse of the form \eqref{Gaussian_pulse}. The dipole is $z$-directed and lies at the point $(0.2704,0.4421,0.3902)$ with $Il=0.01$. To compute the initial condition, we use the approach described in Section \ref{section_cylinder_scattering}.

\begin{figure}[htbp!]
\centering $
\begin{array}{cc}
    \includegraphics[height = .31\textheight]{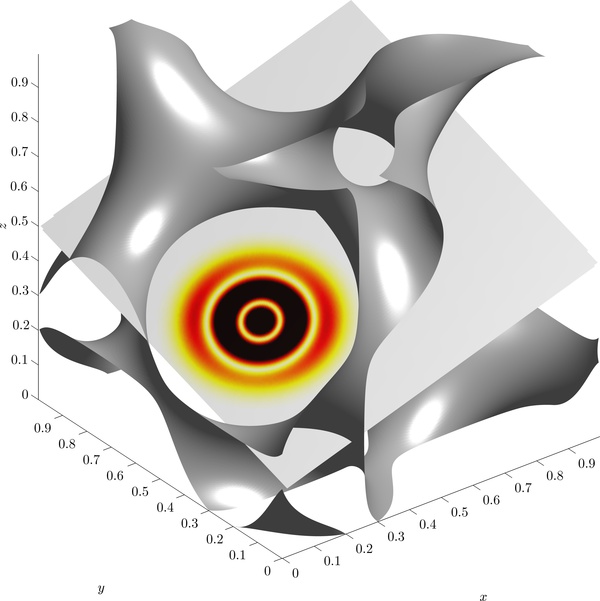} 
    \includegraphics[height = .31\textheight]{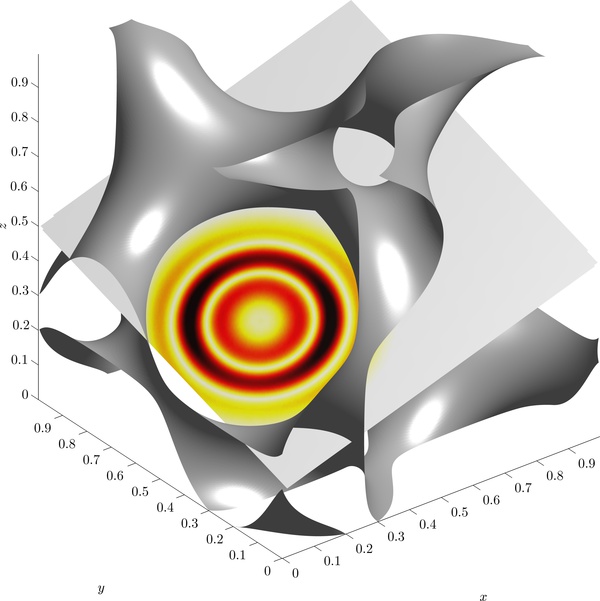} \\
    \includegraphics[height = .31\textheight]{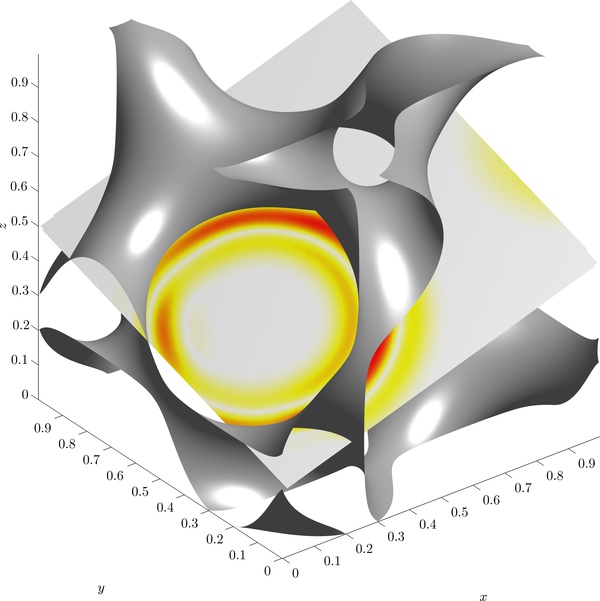} 
    \includegraphics[height = .31\textheight]{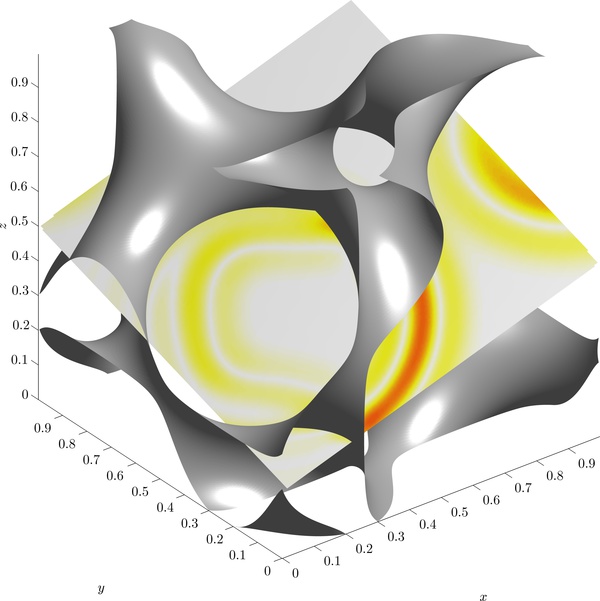} \\
    \includegraphics[height = .31\textheight]{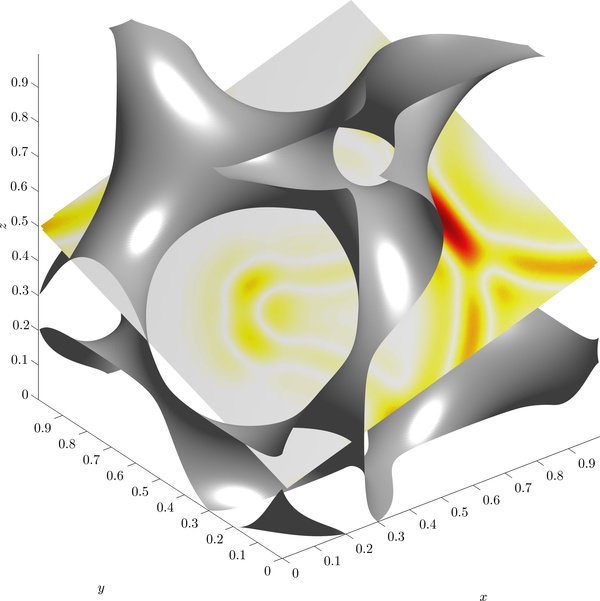} 
    \includegraphics[height = .31\textheight]{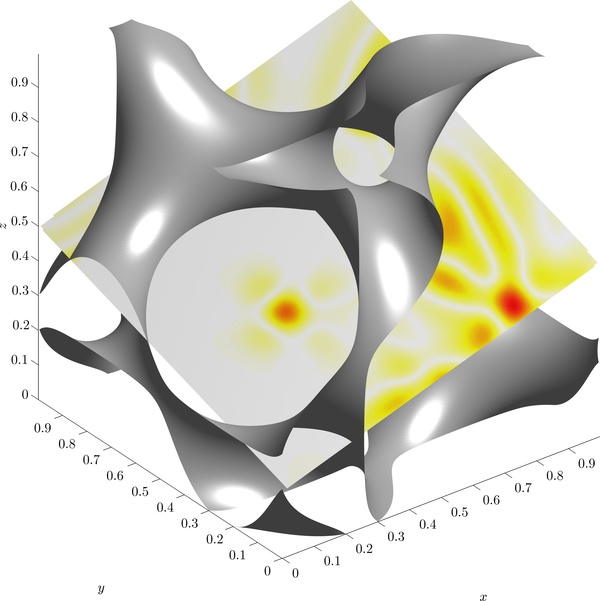}
\end{array} $
\caption{Six snapshots of a slice of the magnitude of the $E_{x,\eta}$ component of the Gaussian dipole wave scattering off the gyroid. From the top left corner, the figures are taken at times $t=0.0999,0.1844,0.3111,0.4096,0.5363,0.6348$. The gyroid surface is shown in dark grey.} \label{gyroid_Ex_snapshots}
\end{figure}

\begin{figure}[htbp!]
\centering $
\begin{array}{cc}
    \includegraphics[width = .42\textwidth]{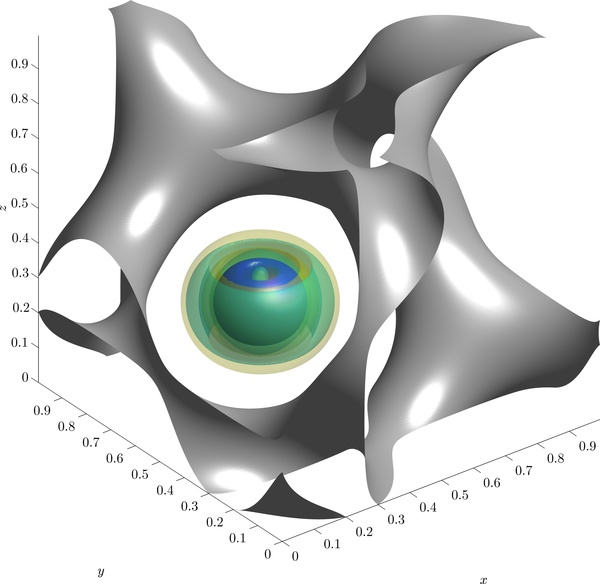} 
    \includegraphics[width = .42\textwidth]{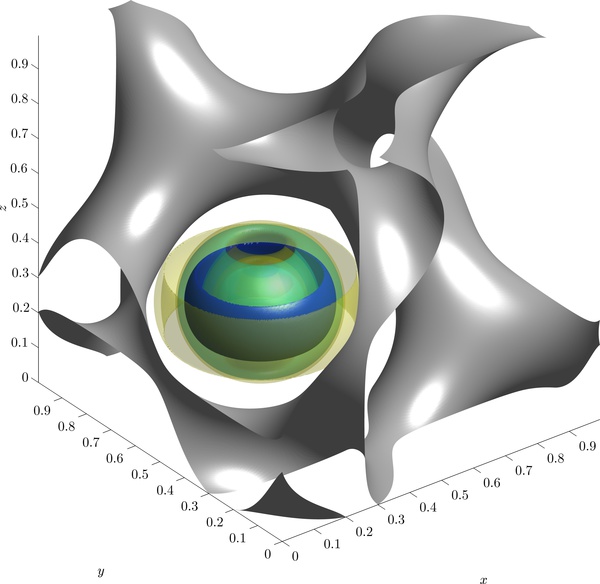} \\
    \includegraphics[width = .42\textwidth]{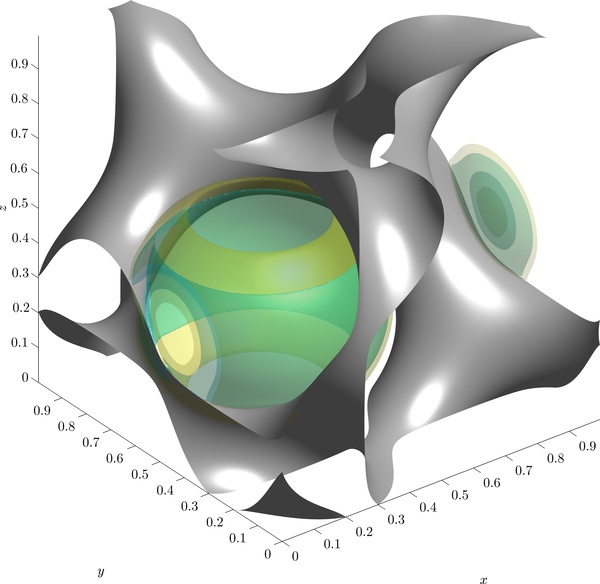} 
    \includegraphics[width = .42\textwidth]{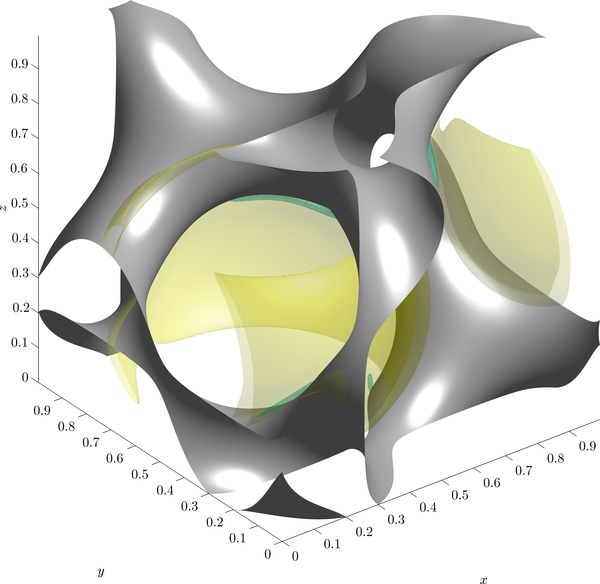} \\
    \includegraphics[width = .42\textwidth]{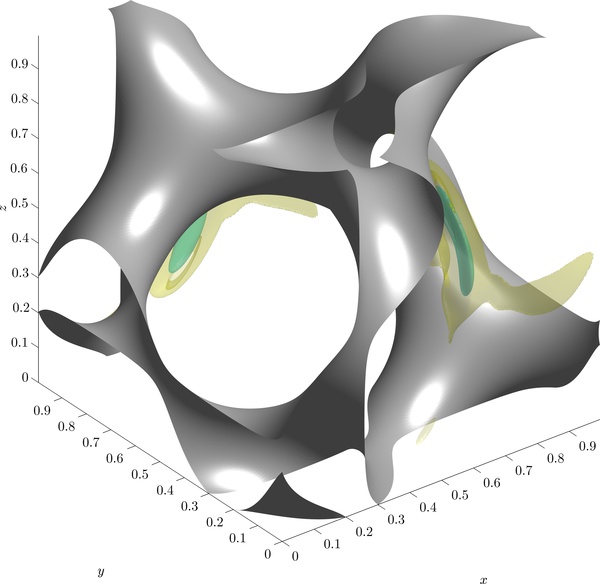} 
    \includegraphics[width = .42\textwidth]{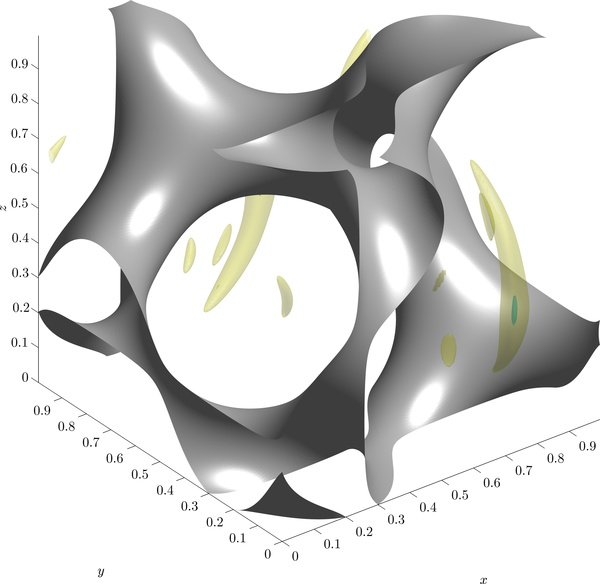}
\end{array} $
\caption{Six snapshots (matching those in Figure \ref{gyroid_Ex_snapshots}) of the energy density of the Gaussian dipole wave scattering off the gyroid. Three level sets ($c=1$ in opaque blue, $c=0.3$ in transparent cyan, and $c=0.15$ in transparent yellow) of the energy density are shown. The gyroid surface is shown in dark grey.} \label{gyroid_energy_snapshots}
\end{figure}

The gyroid is described by the zero level set of 
\begin{equation}
\hat\psi \left(x,y,z \right) = a - \sin \left( 2 \pi x \right)\cos \left( 2 \pi y \right) - \sin \left( 2 \pi y \right) \cos \left( 2 \pi z \right) - \sin \left( 2 \pi z \right) \cos \left( 2 \pi x \right)
\end{equation}
with $a=0.95$ which is not a signed distance function (see Section \ref{Sec_NumSolveCoordinates}). Figure \ref{gyroid_Ex_snapshots} illustrates a slice of the $E_{x,\eta}$ component of the radiating dipole for various times while Figure \ref{gyroid_energy_snapshots} illustrates three level sets of the corresponding energy density $\frac{1}{2}\left( |\mathbf{E}_{\eta}|^2 + |\mathbf{B}_{\eta}|^2 \right)$ of the computed wave solution. The figures were generated using $N=256$, $T=1$, $\Delta t = 0.35 \Delta x$, $\eta = 5 \Delta t$, $h=2 \Delta x$, and $L=0.025$.

\section{Conclusions}

In this paper, we have introduced a Fourier based penalty method for solving Maxwell's equations in domains with curved boundaries and perfect electric conductor boundary conditions.  The approach relied on embedding the physical domain in a larger periodic computational domain, followed by the introduction of a penalty forcing term.  We demonstrate that by constructing a penalty term that is a continuous extension of the electric field and that also satisfies the exact boundary condition, we may systematically improve the analytic convergence of the penalized PDE to the exact underlying PDE.  We show by analytic calculations in two dimensions that one achieves high order convergence for a $\textrm{TM}_z$ mode scattering off a straight wall. We also show through the direct computation of numerical eigenvalues that the scheme is numerically stable in dimension one (for $m = 0, 1, 2$) and dimension two (for $m = 0$). We conclude with several numerical examples of our Fourier based approach. Specifically, we show high order convergence in dimension one, as well as a lack of dispersion errors which typically result when solving for wave propagation at high frequencies.  We demonstrate the approach with several more practical examples including propagation in a waveguide geometry and scattering off a windmill-like geometry. Finally, we confirm that the method extends to three dimensions.

Despite the simplicity of the approach, several issues can still be improved. Future work aims to further improve the efficiency and simplicity of constructing the extension $\mathbf{\tilde{g}}$ through the formulation of a minimization problem.  In doing so, one can likely avoid the added step of solving (\ref{LocalCoordinates}) to compute the local coordinates.  Secondly, additional stability details arise in dimensions two and three that are not present in dimension one and that currently limit the accuracy of the method to either $2.5$ for $\textrm{TM}_z$ modes, or $1.5$ for $\textrm{TE}_z$ modes.  These issues appear due to the conditioning of the current construction for $\mathbf{\tilde{g}}$, which relies on building smooth extensions along rays. The conditioning may potentially be improved by taking an alternative, basis based, approach to the construction of the extension $\mathbf{\tilde{g}}$.  We leave the investigation of alternative constructions of $\mathbf{\tilde{g}}$ for future work. Finally, one may consider a full $\left( \mathbf{E},\mathbf{H} \right)$ penalization. This could raise the convergence rate of $\mathbf{H}$ by half an order (to match that of $\mathbf{E}$) at the expense of a more complicated scheme.


\section{Acknowledgments}
The authors would like to thank Mark Lyon, Dmitry Kolomenskiy and Kai Schneider for numerous enlightening conversations.

This research was partly supported through the NSERC Discovery and Discovery Accelerator Supplement grants of the third author.

This work was supported by a grant from the Simons Foundation ($\#359610$, David Shirokoff).

\bibliography{main}
\bibliographystyle{plain}
\end{document}